\documentclass[a4paper,10pt,final]{article}
\usepackage{a4wide}
\usepackage{amsmath,amssymb,amsfonts}
\usepackage[utf8]{inputenc} 
\usepackage{algorithm}
\usepackage[T1]{fontenc}
\usepackage{enumerate}
\usepackage{pdfpages}
\usepackage{graphicx}
\usepackage{hyperref}
\usepackage{tabularx}
\usepackage{xcolor}
\usepackage{url}
\newcommand{\MyParagraph}[1]{\smallskip\noindent\textbf{#1.}\;}
\newcommand{\Myparagraph}[1]{\smallskip\noindent\emph{#1.}\;}

\newcommand{\NN}{\mathbb{N}}
\newcommand{\ZZ}{\mathbb{Z}}
\newcommand{\RR}{\mathbb{R}}
\newcommand{\CC}{\mathbb{C}}
\newcommand{\dd}{\mathrm{d}}
\newcommand{\ee}{\mathrm{e}}
\newcommand{\ii}{\mathrm{i}}
\newcommand{\abs}[1]{\vert#1\vert}
\newcommand{\nE}{\mathcal{E}}
\newcommand{\nF}{\mathcal{F}}
\newcommand{\nI}{\mathcal{I}}
\newcommand{\nJ}{\mathcal{J}}
\newcommand{\nL}{\mathcal{L}}
\newcommand{\nM}{\mathcal{M}}
\newcommand{\nQ}{\mathcal{Q}}
\newcommand{\nQc}{\nQ^{(c)}}
\newcommand{\MyGrad}[1]{\nabla_{\!#1}}
\newcommand{\nIm}{\nI^{(m)}}
\newcommand{\nIpsi}{\nI^{(\psi)}}
\newcommand{\nJpsi}{\nJ^{(\psi)}}
\newcommand{\nIphipsi}{\nI^{(\varphi - \psi)}}
\newcommand{\nJphipsi}{\nJ^{(\varphi - \psi)}}
\newcommand{\nImpsi}{\nI^{(m,\psi)}}
\newcommand{\nImphipsi}{\nI^{(m,\varphi - \psi)}}
\newcommand{\nEb}{\nE^{(b)}}
\newcommand{\nFb}{\nF^{(b)}}
\newcommand{\mub}{\mu^{(b)}}
\newcommand{\Omb}{\Omega^{(b)}}
\newcommand{\Omo}{\Omega^{(0)}}
\begin{document}
%%%%%%%%%%%%%%%%%%%%%%%%%%%%%%%%%%%%%%%%%%%%%%%%%%%%%%%%%%%%%%%%%%%%%%%%%%%%%%%%%%%%%%%%%%%%%%%%%%%%%%%%%%%%%%%%%%%%%%%%%%%%%%%%%%%%%%%%%%%%%%%%%% 
%%%%%%%%%%%%%%%%%%%%%%%%%%%%%%%%%%%%%%%%%%%%%%%%%%%%%%%%%%%%%%%%%%%%%%%%%%%%%%%%%%%%%%%%%%%%%%%%%%%%%%%%%%%%%%%%%%%%%%%%%%%%%%%%%%%%%%%%%%%%%%%%%% 
\title{Novel approaches for the reliable and efficient \\ numerical evaluation of Landau-type operators}
\author{Jos{\'e} Antonio Carrillo, Mechthild Thalhammer
  \footnote{Addresses: J.\,A.~Carrillo, Mathematical Institute, University of Oxford, Andrew Wiles Building, Radcliffe Observatory Quarter, Woodstock Road, OX2 6GG Oxford, United Kingdom.
M.~Thalhammer, Department of Mathematics, University of Innsbruck, Techniker\-strasse~13/7, 6020~Innsbruck, Austria.
Emails: \url{jose.carrillo@maths.ox.ac.uk}, \url{mechthild.thalhammer@uibk.ac.at}. 
Websites: \url{www.maths.ox.ac.uk/people/jose.carrillodelaplata}, \url{techmath.uibk.ac.at/mecht}. }}
%\author[Carrillo J.\,A. and Thalhammer M.]{Jos{\'e} Antonio Carrillo\affil{1}, Mechthild %Thalhammer\affil{2}\comma\corrauth}
%\address{\affilnum{1}Mathematical Institute, University of Oxford, Andrew Wiles %Building, Radcliffe Observatory Quarter, Woodstock Road, OX2 6GG Oxford, United %Kingdom. \\
%\affilnum{2}Department of Mathematics, University of Innsbruck, %Techniker\-strasse~13/7, 6020~Innsbruck, Austria.}  
%\emails{{\tt jose.carrillo@maths.ox.ac.uk} (J.\,A.~Carrillo), {\tt %mechthild.thalhammer@uibk.ac.at} (M.~Thalhammer)}
% Websites: \url{www.maths.ox.ac.uk/people/jose.carrillodelaplata}, %\url{techmath.uibk.ac.at/mecht}. 
%%%%%%%%%%%%%%%%%%%%%%%%%%%%%%%%%%%%%%%%%%%%%%%%%%%%%%%%%%%%%%%%%%%%%%%%%%%%%%%%%%%%%%%%%%%%%%%%%%%%%%%%%%%%%%%%%%%%%%%%%%%%%%%%%%%%%%%%%%%%%%%%%% 
\maketitle
%%%%%%%%%%%%%%%%%%%%%%%%%%%%%%%%%%%%%%%%%%%%%%%%%%%%%%%%%%%%%%%%%%%%%%%%%%%%%%%%%%%%%%%%%%%%%%%%%%%%%%%%%%%%%%%%%%%%%%%%%%%%%%%%%%%%%%%%%%%%%%%%%% 
\subsubsection*{Abstract}
Numerical approximations of Landau-type operators represent fundamental components of time integration methods for demanding problems such as inhomogeneous Vlasov--Landau-type equations.
Substantial computational issues arise from the treatment of the physically most relevant three-dimensional case with Coulomb-type interaction.
This work is concerned with the introduction and numerical comparison of novel approaches for the reliable and efficient evaluation of Landau-type collision operators, where the focus is on the treatment of integral operators involving general singular kernels.
In the spirit of collocation, common tools are the identification of fundamental integrals, series expansions of the integral kernel and the density function on the main part of the velocity domain, and interpolation as well as quadrature approximation nearby the singularity of the kernel. 
Focusing on the favourable choice of the Fourier spectral method, their practical implementation uses the reduction to basic integrals, fast Fourier techniques, and summations along certain directions.
Moreover, an important observation is that a significant percentage of the overall computational effort can be transferred to precomputations which are independent of the density function. 
For the purpose of exposition and numerical validation, the cases of constant, regular, and singular integral kernels are distinguished, and the procedure is adapted accordingly to the increasing complexity of the problem.
%\ams{65D30,65M12,65M20,65M70}%52B10, 65D18, 68U05, 68U07
%\keywords{Vlasov--Landau-type system, Landau-type equation, Landau-type %collision operator, Numerical approximation, Fourier spectral method, Operator %splitting, Reliability, Efficiency, Mass conservation.}
%%%%%%%%%%%%%%%%%%%%%%%%%%%%%%%%%%%%%%%%%%%%%%%%%%%%%%%%%%%%%%%%%%%%%%%%%%%%%%%%%%%%%%%%%%%%%%%%%%%%%%%%%%%%%%%%%%%%%%%%%%%%%%%%%%%%%%%%%%%%%%%%%% 
%%%%%%%%%%%%%%%%%%%%%%%%%%%%%%%%%%%%%%%%%%%%%%%%%%%%%%%%%%%%%%%%%%%%%%%%%%%%%%%%%%%%%%%%%%%%%%%%%%%%%%%%%%%%%%%%%%%%%%%%%%%%%%%%%%%%%%%%%%%%%%%%%% 
%%%%%%%%%%%%%%%%%%%%%%%%%%%%%%%%%%%%%%%%%%%%%%%%%%%%%%%%%%%%%%%%%%%%%%%%%%%%%%%%%%%%%%%%%%%%%%%%%%%%%%%%%%%%%%%%%%%%%%%%%%%%%%%%%%%%%%%%%%%%%%%%%% 
\section{Introduction}
\label{sec:Introduction}
%%%%%%%%%%%%%%%%%%%%%%%%%%%%%%%%%%%%%%%%%%%%%%%%%%%%%%%%%%%%%%%%%%%%%%%%%%%%%%%%%%%%%%%%%%%%%%%%%%%%%%%%%%%%%%%%%%%%%%%%%%%%%%%%%%%%%%%%%%%%%%%%%% 
%%%%%%%%%%%%%%%%%%%%%%%%%%%%%%%%%%%%%%%%%%%%%%%%%%%%%%%%%%%%%%%%%%%%%%%%%%%%%%%%%%%%%%%%%%%%%%%%%%%%%%%%%%%%%%%%%%%%%%%%%%%%%%%%%%%%%%%%%%%%%%%%%% 
\MyParagraph{Scope of applications} 
The present work is inspired by various contributions on the numerical simulation of kinetic equations modelling the distribution of charged particles in a collisional plasma, see in particular \cite{bailo2024collisional,BCM2023,CarrilloEtAl2020,DLP2015,DP2014,FilbetPareschi2002,Gamba2017,HuQi2020,LiEtAl2019,MangeneyEtAl2002,PR22,PareschiRusso2000,PareschiEtAl2000,PG2019,Wang2015,ZhangGamba2017} and the references given therein.
We intend to lay the foundation for a future application of operator splitting methods to Vlasov--Maxwell--Landau-type and Vlasov--Poisson--Landau-type equations, where the efficient time integration of a Landau-type equation via spectral methods represents a fundamental component of the entire algorithm.

\MyParagraph{Vlasov--Landau-type equation}
Solving an inhomogeneous Vlasov--Landau-type equation~\cite{Landau1936} is still one of the most computationally costly kinetic problems in the field and of paramount importance in plasma physics.
With the variables $x \in \Omega^{(x)} \subseteq \RR^d$, $v \in \Omega^{(v)} \subseteq \RR^d$, and $t \in [t_0, T] \subset \RR$ representing position, velocity, and time, the functions $f: \Omega^{(x)} \times \Omega^{(v)} \times [t_0, T] \rightarrow \RR$ and $F: \Omega^{(x)} \times \Omega^{(v)} \times [t_0, T] \rightarrow \RR$ describing the distribution of charged particles and a given or self-consistent force field including electromagnetic effects, and a Landau-type operator~$\nQ(f, f)$ capturing collisions between particles, the associated Vlasov--Landau-type equation reads as  
\begin{equation*}
\partial_t f + v \cdot \nabla_x f - F \cdot \nabla_v f = \nQ(f, f)\,.
\end{equation*}
Throughout, for notational simplicity, we neglect dependences on variables when appropriate and no confusion arises.

\MyParagraph{Landau-type operator}
Our main concern is the efficient numerical evaluation of a Landau-type collision operator
\begin{subequations}
\label{eq:LandauOperatorxvt}
\begin{equation}
\nQ(f, f) = \MyGrad{v} \cdot \nQc(f, f)\,, 
\end{equation}
which is given by the divergence of the integral operator 
\begin{equation}
\label{eq:FormulaQc1}
\begin{split}
\nQc(f, f)(v) 
&= C \int_{\Omega^{(v)}} A(v-w) \, \big(f(w) \, \MyGrad{v} f(v) - \MyGrad{w} f(w) \, f(v)\big) \; \dd w \\
&= \nI(f)(v) \, \MyGrad{v} f(v) + \nJ(f)(v) \, f(v)\,, \quad v \in \Omega^{(v)}\,.
\end{split}
\end{equation}
Here, we set
\begin{equation}
A(z) = \varphi(z) \, (\abs{z}^2 \, I_d - z \otimes z)\,, \quad z \in \RR^d\,, 
\end{equation}
\end{subequations}
and denote by $\abs{z} = \sqrt{z^T z} \in \RR_{\geq 0}$ and $z \otimes z = z \, z^T \in \RR^{d \times d}$ the Euclidean norm and the outer product of a column vector $z \in \RR^d$, by $I_d \in \RR^{d \times d}$ the identity matrix, and by $C > 0$ a positive constant.
For the gradient operator comprising the partial derivatives with respect to the velocity components $v = (v_1, \dots, v_d)^T \in \RR^d$, we employ the notation $\MyGrad{v} = (\partial_{v_1}, \dots, \partial_{v_d})^T$.
The purpose of efficiently evaluating the Landau-type operator~\eqref{eq:LandauOperatorxvt} is closely related to the numerical approximation of the associated Landau-type equation 
\begin{equation}
\label{eq:LandauEquation}
\partial_t f = \MyGrad{v} \cdot \nQc(f, f)\,.
\end{equation}

\MyParagraph{Maxwellian molecules and Coulomb interaction}
The choices $\Omega^{(v)} = \RR^d$ with $d \in \{2, 3\}$ and $\varphi(z) = C$ with some $C \in \RR$, referred to as Maxwellian molecules cases, serve as basic test examples, since the specification of density functions like 
\begin{equation*}
\begin{gathered}
d = 2\,, \quad f(v) = \ee^{- \abs{v}^2} \, \abs{v}^2\,, \quad v \in \RR^2 \,, \\
d = 3\,, \quad f(v) = \ee^{- \abs{v}^2} \, \big(2 \, \abs{v}^2 - 1\big)\,, \quad v \in \RR^3 \,, 
\end{gathered}
\end{equation*}
permits to determine the integrals inside the Landau operators.
As another well-established test problem, we also consider the physically relevant and numerically challenging case of Coulomb interaction in three dimensions
\begin{equation*}
d = 3\,, \quad \Omega^{(v)} = \RR^3\,, \quad \varphi(z) = C \, \abs{z}^{- \, 3}\,, \quad z \in \RR^3\,, 
\end{equation*}
which leads to a strong singularity in the integral operator. Remarkably, explicit expressions of the Fourier expansion of the Landau operator are known in this particular case,  leading to fast spectral methods~\cite{PareschiRusso2000,ZhangGamba2017}.

\MyParagraph{General setting}
In contrast to the above situations, where $\varphi(z) = C \, \abs{z}^{\beta}$ for some constant $C > 0$ and exponent $\beta \geq - \, 3$, we here admit more general kernels.
Relevant instances
\begin{equation}
\label{eq:KernelCoulombGeneral}
\varphi(z) = C \, \abs{z}^{\beta} \, \ee^{\gamma \, \abs{z}}\,, \quad \beta, \gamma < 0\,,  \quad z \in \RR^d\,, 
\end{equation}
 again include the case of an isolated strong singularity in the integral operator.
Due to the fact that a coupling of all velocity directions is unavoidable in such settings , a further development of approaches previously studied in \cite{FilbetPareschi2002,PR22,PareschiRusso2000,PareschiEtAl2000,ZhangGamba2017} is required. However, as explicit formulas for the Fourier expansion of the Landau operator are no longer available, the canonical approach relies on the implementation of accurate quadrature approximations for the arising integrals. 
In the following, we detail an alternative strategy that significantly improves accuracy and reduces the computational effort. 
More general integral kernels are features common to other Landau-type operators such as the Lenard--Balescu operator, see \cite{DW23} and the references therein, or simplified related equations, see for instance \cite{scullard2016numerical,scullard2020adaptive}, where our approach could be potentially used as well.

\MyParagraph{Simplifications}
Henceforth, we tacitly assume that the density function satisfies suitable regularity and integrability requirements such that~\eqref{eq:LandauOperatorxvt} is well-defined in the classical sense.
For the sake of concreteness, we assume that a single isolated singularity of the kernel arises at the origin.
For simplicity, we refer to the Landau-type operator~\eqref{eq:LandauOperatorxvt} as Landau operator and to the Landau-type equation~\eqref{eq:LandauEquation} as Landau equation in the rest of this work.

\MyParagraph{Novel strategy}
In this work, we propose a strategy in the spirit of collocation methods for the reliable and efficient numerical evaluation of the Landau operator~\eqref{eq:LandauOperatorxvt}. 
For illustrative purposes, we apply it in the time integration of the Landau equation~\eqref{eq:LandauEquation}, keeping in mind a future implementation for more demanding problems such as Vlasov--Landau equations based on operator splitting. 
Our alternative spectral method avoids quadrature approximations on the entire velocity domain by Fourier expanding not only the unknown density distribution but also the integral kernel. Quadrature approximations are  needed nearby the singularity of the kernel, where a Fourier series representation is not suitable and inaccurate.
Let us emphasise once again that 
the value of our strategy relies on the fact that we deal with more general singular kernels. 
On the one hand, 
we gain spectral accuracy instead of polynomial quadrature order. 
On the other hand, by advances in the precomputational stage, the computational effort for the time evolution is reduced significantly. 
Substantial  simplifications are feasible, provided that explicit expressions of the arising integrals are available, typically in the case of Coulomb interaction~\cite{PareschiRusso2000,ZhangGamba2017}.
When we employ a representation of the Landau-type operators in divergence form, conservation of mass is ensured. However, in comparison with~\cite{PR22}, the 
preservation of momenta is not an intrinsic property of our approach.

\MyParagraph{Basic concepts}    
We next sketch basic concepts for the numerical approximation of \eqref{eq:LandauOperatorxvt}-\eqref{eq:LandauEquation}  with~\eqref{eq:KernelCoulombGeneral} and essential components of the resulting algorithms.
Our common starting point is the representation of the Landau equation as a nonlocal drift-diffusion equation. 
We identify fundamental integrals of the form
\begin{subequations}
\begin{equation*}
\int_{\RR^3} \varphi(v - w) \, p(v - w) \, g(w) \; \dd w
\end{equation*}
that involve the singular integral kernel $\varphi: \RR^3 \to \RR$, a polynomial of degree two $p: \RR^3 \to \RR$, and a regular function $g: \RR^3 \to \RR$ reflecting the values of the density function or derivatives thereof, respectively.
In addition, we make use of the decomposition 
\begin{align}
\label{eq:Integral1}
\int_{\RR^3} \varphi(v - w) \, p(v - w) \, g(w) \; \dd w  = &\,\int_{\RR^3} \psi(v - w) \, p(v - w) \, g(w) \; \dd w \\
&\, + \int_{\RR^3} (\varphi - \psi)(v - w) \, p(v - w) \, g(w) \; \dd w\,, \nonumber
\end{align}
where~$\psi$ is a suitable regularisation of the kernel, obtained by interpolation nearby the isolated singularity, such that the difference~$\varphi - \psi$ vanishes on the main part of the velocity domain, see Figure~\ref{fig:Kernel}.
In order to numerically compute these integrals, we employ series expansions of~$\psi$ and~$g$ as well as quadrature approximations. 
Due to the particular properties of Fourier functions 
\begin{equation}
\label{eq:Fourier1}
\begin{gathered}
\nF_{\!\kappa}: \RR \longrightarrow \RR: \xi \longmapsto \tfrac{1}{\sqrt{2 \, b}} \; \ee^{\, \mu_{\kappa} (\xi + b)}\,, \\
\nF_{\!m}: \RR^3 \longrightarrow \RR: v \longmapsto \nF_{\!m_1}(v_1) \, \nF_{\!m_2}(v_2) \, \nF_{\!m_3}(v_3)\,, \\
\mu_{\kappa} = \tfrac{\pi \, \ii \, \kappa}{b} \in \CC\,, \quad \mu_m = \big(\mu_{m_1}, \mu_{m_2}, \mu_{m_3}\big) \in \CC^{1 \times 3}\,, \\
b > 0\,, \quad \kappa \in \ZZ\,, \quad m = (m_1, m_2, m_3) \in \ZZ^3\,, 
\end{gathered}
\end{equation}
and in view of the highly efficient practical implementation of Fourier series expansions by fast transforms, i.e.
\begin{equation}
\label{eq:Expansion1}  
\sum_{m \in \nM} g_m \, \nF_m \approx g\,, \quad \sum_{m \in \nM} \psi_m \, \nF_m \approx \psi\,, \quad \nM = \big\{- \tfrac{M}{2}\,, \dots, \tfrac{M}{2} - 1\big\}^3\,,
\end{equation}
for an even positive integer number $M \in \NN$, we favour the Fourier spectral method.
Under the reasonable presumption of a localised density function, we may replace the unbounded velocity domain by a Cartesian product of intervals, characterised by a sufficiently large positive real number $b > 0$. 
Accordingly, we choose uniform grid points that cover the truncated domain
\begin{equation}
\label{eq:Domain1}  
v_{\ell} \in [- \, b, b]^3\,, \quad \ell = (\ell_1, \ell_2, \ell_3) \in \nL = \{1\,, \dots, M\}^3\,, 
\end{equation}
Likewise, the relatively small neighbourhood of the origin, where the interpolant of the singular kernel is devised
\begin{equation}
\label{eq:DomainSmall1}  
(\varphi - \psi)\big\vert_{[- \, b, b]^3 \setminus [- \, b_0, b_0]^3} = 0\,, 
\end{equation}
is defined by a positive real number $0 < b_0 <\!\! < b$, which we adjust in such a way that the point $(b_0, 0, 0)$ coincides with a grid point.
Thus, non-zero values 
\begin{equation}
\label{eq:SetSmall1}  
(\varphi - \psi)(v_{\ell}) \neq 0\,, \quad \ell \in \widetilde{\nL} \subset \nL\,, 
\end{equation}
occur, but only for a small subset $\widetilde{\nL} \subset \nL$.
Because of the multiplicativity of Fourier functions, the concrete tasks for the numerical computation of the fundamental integrals~\eqref{eq:Integral1} amount to determine one-dimensional integrals of the form 
\begin{equation}
\label{eq:Integral1a}  
\int_{- b}^{b} \xi^i \, \nF_{\!\kappa}(\xi) \; \dd \xi\,, \quad i \in \{0, 1, 2\}\,, \quad \kappa \in \big\{- \tfrac{M}{2}\,, \dots, \tfrac{M}{2} - 1\big\}\,,
\end{equation}
to approximate three-dimensional integrals by quadrature 
\begin{equation}
\label{eq:Integral1b}  
\int_{[- \, b_0, b_0]^3} (\varphi - \psi)(w) \, p(w) \, \nF_{\!m}(w) \; \dd w\,, \quad m \in \nM\,, 
\end{equation}
\end{subequations}
and to apply fast Fourier transforms or summations along certain directions, respectively.
A benefit of our strategy to decompose  the kernel and accordingly the domain of integration is that we are improving the polynomial order to spectral order on the main part of the domain. The width $b_0$ has to be determined in a precomputation, since it is affected by the properties of the kernel. 
We propose to choose it as small as possible such that the kernel and its Fourier expansion coincide outside $[-b_0,b_0]^3$ up to high accuracy.

\MyParagraph{Comparison with alternative spectral methods}
The efficient implementation of our approaches strongly relies on spectral methods, in particular, on fast Fourier techniques. 
However, compared to~\cite{FilbetPareschi2002,PR22,PareschiRusso2000,PareschiEtAl2000,ZhangGamba2017}, rather than approximating the weak formulation of the Landau operator, we approximate the distribution function in physical space by collocation. 
One of the advantages of such a Fourier collocation method is the localisation of the singular kernel in velocity space and the accurate representation of the distribution function by Fourier series outside the truncated domain.
As a consequence, in connection with the future application of splitting methods for Vlasov--Landau-type systems, e.g., we can hand over the values of the density function to the solvers of the other subproblems. 
We propose different approaches using or avoiding numerical differentiation, depending on whether mass conservation or higher accuracy is preferable, respectively, see Table~\ref{tab:Table1} for an overview. 
Since the kernel is well represented by Fourier series except on a small neighbourhood of the singularity, restricting the quadrature approximation to this neighbourhood has the advantage of significantly reducing the precomputation times as well as the computation times during evolution, see Table~\ref{tab:Table2}.

\MyParagraph{Outline}
This manuscript has the following structure. 
In Section~\ref{sec:Section2}, we introduce our numerical methods for the evaluation of Landau operators involving  general singular integral kernels.  Comments on simplifications for regular kernels and auxiliary results are collected in the appendix. 
Numerical comparisons for different test problems are discussed in Section~\ref{sec:Section3}. 
%%%%%%%%%%%%%%%%%%%%%%%%%%%%%%%%%%%%%%%%%%%%%%%%%%%%%%%%%%%%%%%%%%%%%%%%%%%%%%%%%%%%%%%%%%%%%%%%%%%%%%%%%%%%%%%%%%%%%%%%%%%%%%%%%%%%%%%%%%%%%%%%%% 
%%%%%%%%%%%%%%%%%%%%%%%%%%%%%%%%%%%%%%%%%%%%%%%%%%%%%%%%%%%%%%%%%%%%%%%%%%%%%%%%%%%%%%%%%%%%%%%%%%%%%%%%%%%%%%%%%%%%%%%%%%%%%%%%%%%%%%%%%%%%%%%%%% 
\section{Numerical methods}
\label{sec:Section2}
%%%%%%%%%%%%%%%%%%%%%%%%%%%%%%%%%%%%%%%%%%%%%%%%%%%%%%%%%%%%%%%%%%%%%%%%%%%%%%%%%%%%%%%%%%%%%%%%%%%%%%%%%%%%%%%%%%%%%%%%%%%%%%%%%%%%%%%%%%%%%%%%%% 
%%%%%%%%%%%%%%%%%%%%%%%%%%%%%%%%%%%%%%%%%%%%%%%%%%%%%%%%%%%%%%%%%%%%%%%%%%%%%%%%%%%%%%%%%%%%%%%%%%%%%%%%%%%%%%%%%%%%%%%%%%%%%%%%%%%%%%%%%%%%%%%%%% 
In this section, we develop the proposed spectral collocation methods for the evaluation of the Landau operator and the time integration of the Landau equation.
We begin with an overview of the fundamental approach and then describe in full detail a slight generalisation and related approaches.
%%%%%%%%%%%%%%%%%%%%%%%%%%%%%%%%%%%%%%%%%%%%%%%%%%%%%%%%%%%%%%%%%%%%%%%%%%%%%%%%%%%%%%%%%%%%%%%%%%%%%%%%%%%%%%%%%%%%%%%%%%%%%%%%%%%%%%%%%%%%%%%%%% 
\subsection{Overview}
%%%%%%%%%%%%%%%%%%%%%%%%%%%%%%%%%%%%%%%%%%%%%%%%%%%%%%%%%%%%%%%%%%%%%%%%%%%%%%%%%%%%%%%%%%%%%%%%%%%%%%%%%%%%%%%%%%%%%%%%%%%%%%%%%%%%%%%%%%%%%%%%%% 
\MyParagraph{Preliminary remark}
We recall that our strategy in the spirit of collocation methods is designed for a future time integration of demanding problems such as Vlasov--Landau equations by operator splitting. 
That is, the numerical evaluation of the Landau operator and the short time integration of the Landau equation will be  fundamental components of the entire algorithm. 
Notably, a significant percentage of the overall computational effort can be transferred to precomputations which are independent of the density function.
For this work, for illustrative purposes, we perform the evolution of the Landau equation over longer time frames. 

\MyParagraph{Underlying time integration method}
In the following, we specify details for a time integration of~\eqref{eq:LandauEquation} by means of the explicit Euler method. It is  the simplest representative for standard classes of explicit and implicit, one-step and multi-step time integrators. 
Evidently, the actual choice of the time integration method will depend on the complexity of the problem and the processing capabilities. 
As a general rule, to enhance reliability and efficiency, it will be worthwhile to use an automatic stepsize selection based on a local error control. Specifically, due to the fact that the Landau equation is a drift-diffusion equation, it is expedient to adapt (increase) the temporal increments over longer time frames. 
Our pragmatic and practical realisation  is based on an embedded pair of explicit Runge--Kutta methods, even though the stiffness of the test problems in two and three dimensions imposes considerable stability restrictions. 
An alternative class of methods is provided by Adams--Moulton predictor-corrector methods, which have larger stability regions in comparison with their explicit analogues, the Adams--Bashforth methods.   
Another practicable alternative, but without a straightforward option for local error control, is the use of a stabilised semi-implicit variant of the explicit Euler method, which allows for larger time increments and potentially a resolution of the linear systems by Fourier techniques. 

\MyParagraph{Guide line}
We next  illustrate the key aspects of our algorithms.  For a sequence of positive  stepsizes 
$(\tau_n)_{n=0}^{N-1}$ such that $\tau_0 + \dots + \tau_{N-1} = T$ and a given initial approximation, the explicit Euler solution is determined by the recurrence 
\begin{equation}
\label{eq:ExplicitEuler}
\begin{cases}
&f^{(n)} = f^{(n-1)} + \tau_{n-1} \, \MyGrad{v} \cdot \nQc\big(f^{(n-1)}, f^{(n-1)}\big)\,, \quad n \in \{1, \dots, N\}\,, \\
&f^{(0)} \; \text{given}\,.
\end{cases}
\end{equation}
Evidently, the main computational issue in each time step is related to the numerical evaluation of the integral operator. As an example for a singular kernel, we recall~\eqref{eq:KernelCoulombGeneral} and that we assume that a single isolated singularity at the origin is allowed.
Briefly summarised, we have the following guide line.

\Myparagraph{Input (Problem data)} 
\begin{enumerate}[(i)]
\item
The function~$f^{(0)}$ defines the inital state of the density function, see~\eqref{eq:ExplicitEuler}. 
\item
The function~$\varphi$ defines the integral kernel with isolated singularity at the origin, see~\eqref{eq:LandauOperatorxvt}. 
\end{enumerate}

\begin{figure}[t!]
\begin{center}
\includegraphics[width=4.2cm]{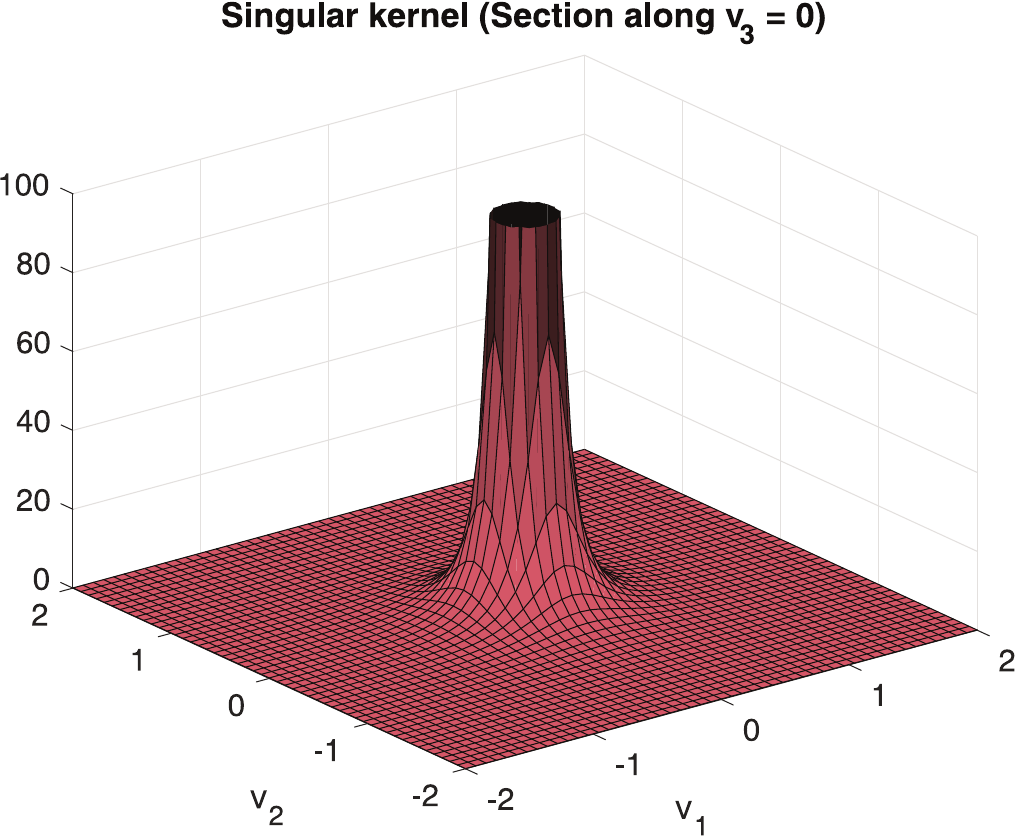} \quad
\includegraphics[width=4.2cm]{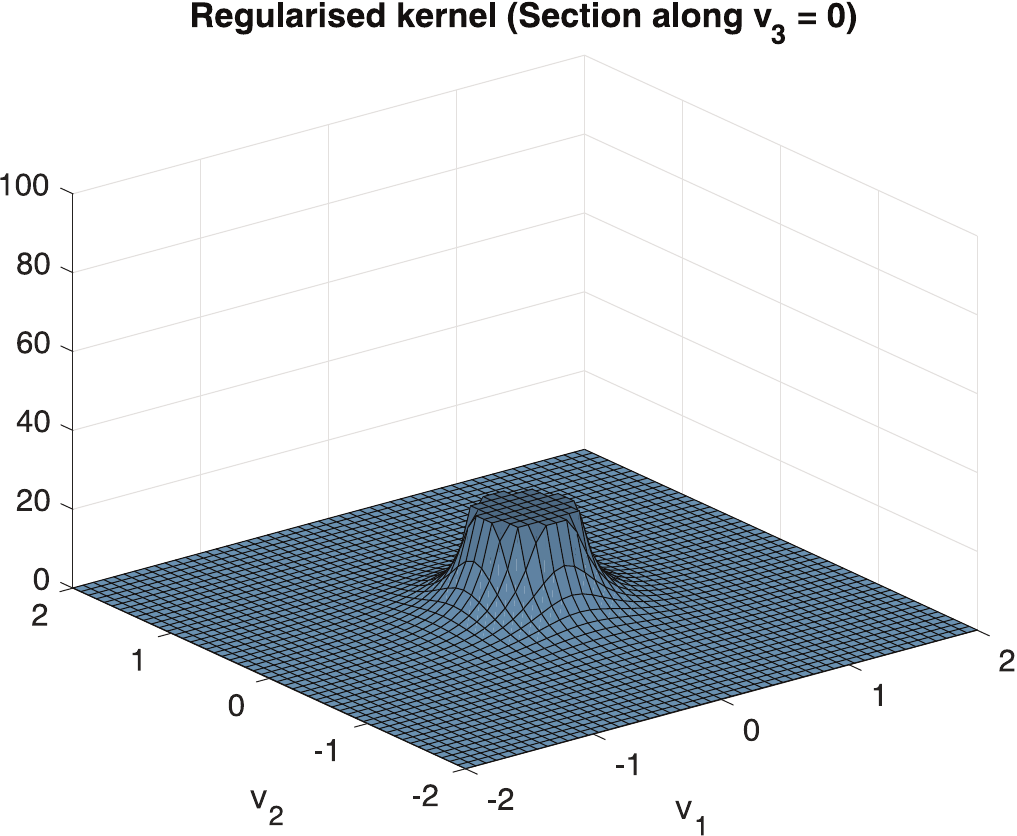} \quad
\includegraphics[width=4.2cm]{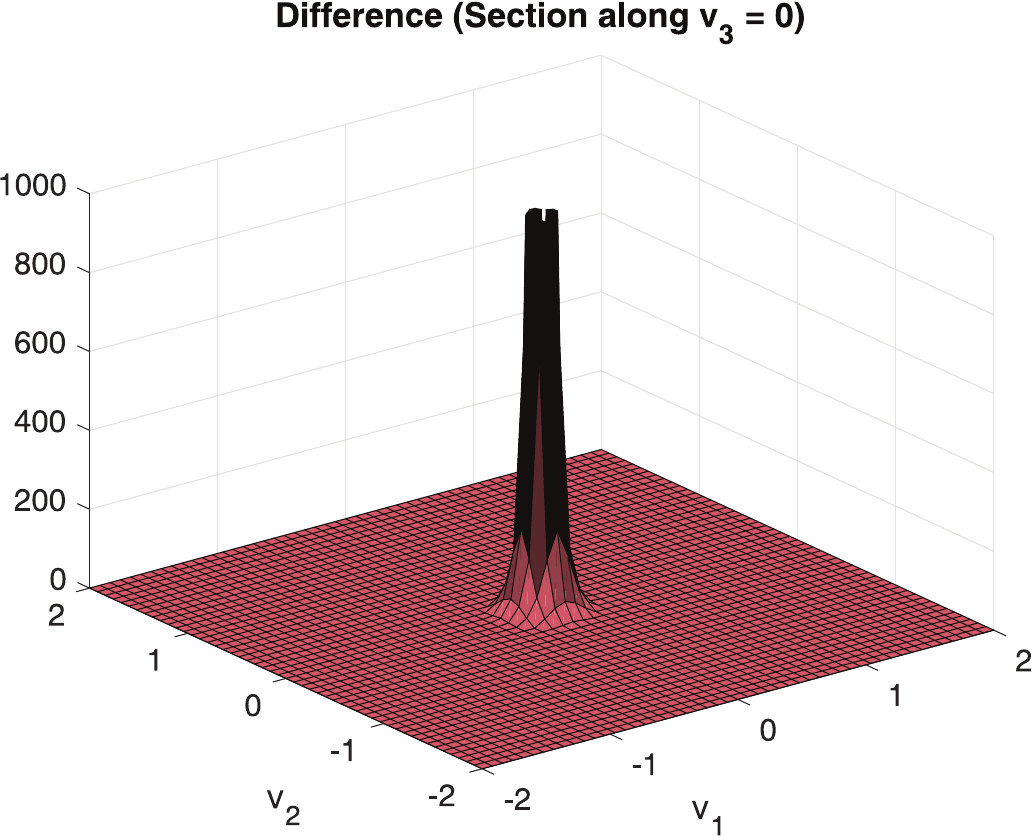} \\[2mm]
\includegraphics[width=4.2cm]{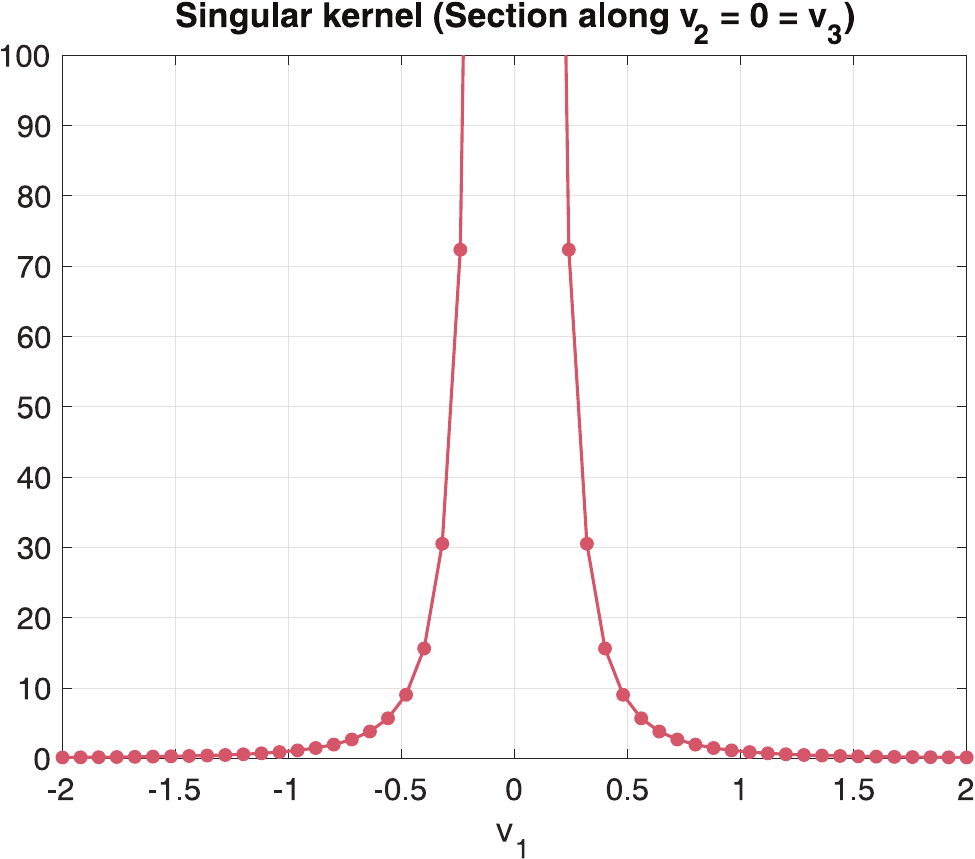} \quad
\includegraphics[width=4.2cm]{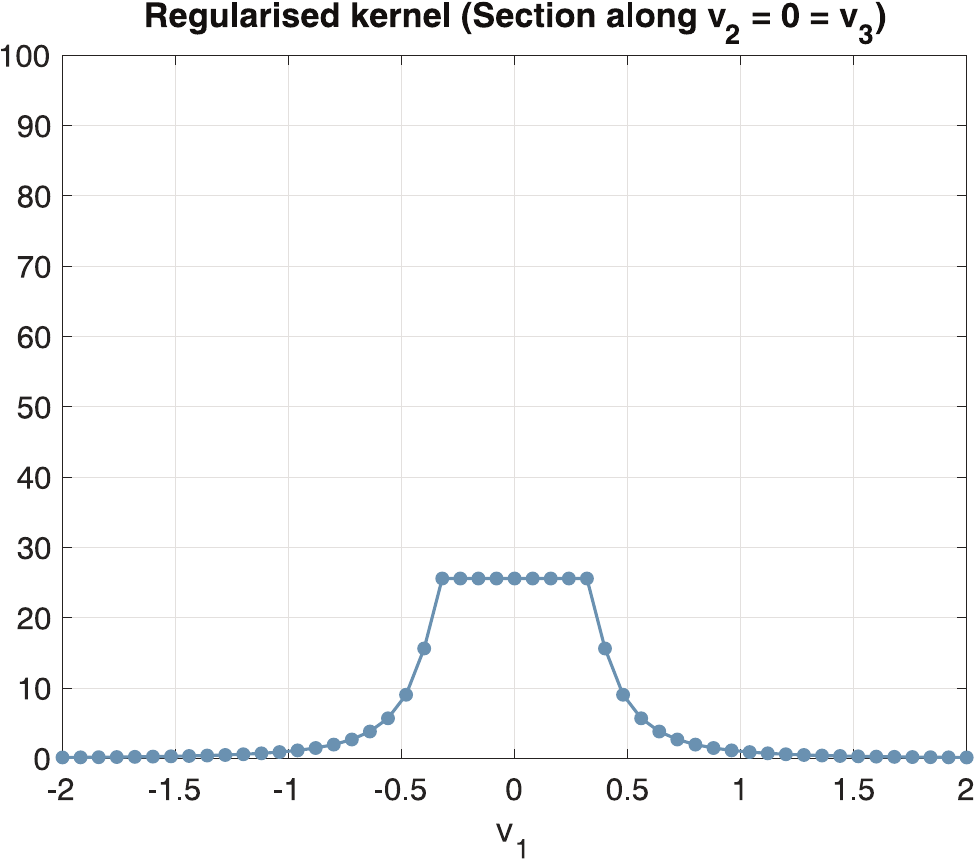} \quad
\includegraphics[width=4.2cm]{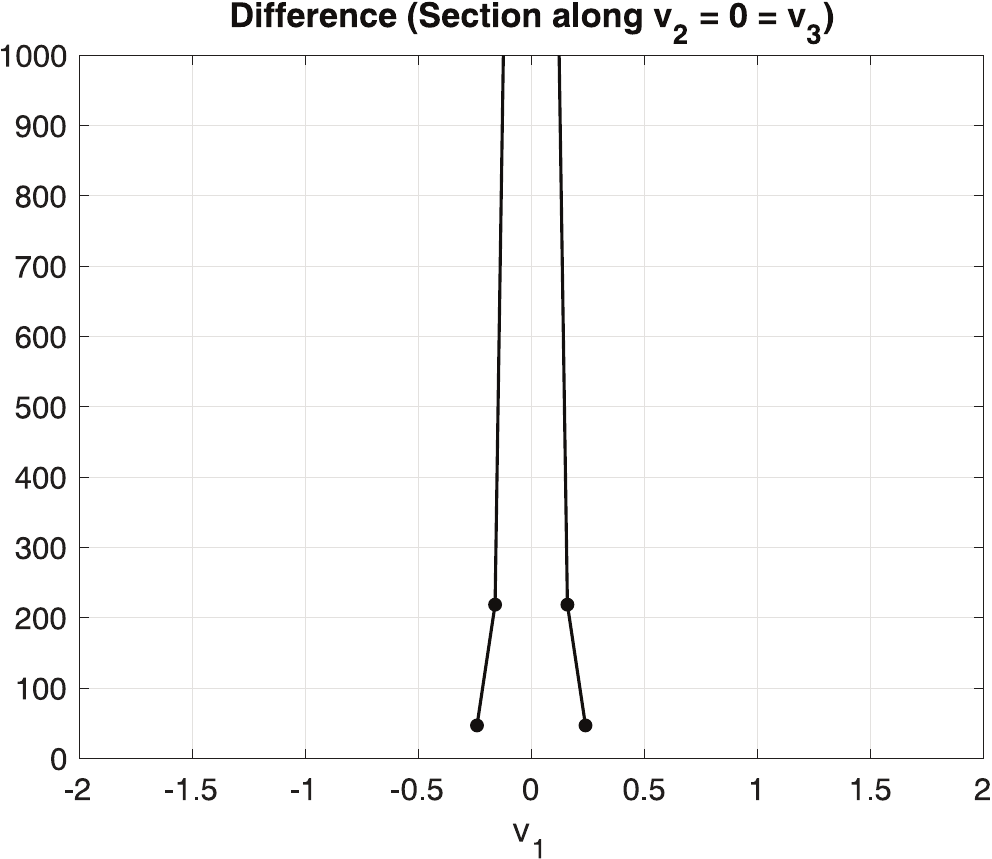} 
\caption{Illustration of singular integral kernels such as $\varphi: \RR^3 \setminus \{0\} \to \RR: z \mapsto \abs{z}^{- \, 3}$ arising in the case of Coulomb interaction and generalisations, see~\eqref{eq:LandauOperatorxvt} and~\eqref{eq:KernelCoulombGeneral}.
A regularised kernel is obtained by interpolation on a small neighbourhood of the origin.
The remaining difference vanishes on the main part of the velocity domain.}
\label{fig:Kernel}
\end{center}
\end{figure}

\Myparagraph{Input (Discretisation)} 
\begin{enumerate}[(i)]
\item
The positive real number $b > 0$ defines the truncated velocity domain, see~\eqref{eq:Domain1}.
\item
The even positive integer number $M \in \NN$ reflects the total numbers of Fourier functions and uniform grid points covering the truncated domain, see~\eqref{eq:Expansion1} and~\eqref{eq:Domain1}.
\item
The small positive real number $b_0 > 0$ is chosen such that $(b_0, 0, 0)$ coincides with a grid point and defines a neighbourhood of the origin, where a regularisation of the singular integral kernel is devised by interpolation, see~\eqref{eq:DomainSmall1} and Figure~\ref{fig:Kernel}. 
\item
A sequence of positive time stepsizes $(\tau_n)_{n=0}^{N-1}$ is defined, and, if required for stability and accuracy, suitably adapted during time integration, see~\eqref{eq:ExplicitEuler}. 
\end{enumerate}

\Myparagraph{Precomputations}
\begin{enumerate}[(i)]
\item
The uniform grid points $(v_{\ell})_{\ell \in \nL}$ are computed and stored, see~\eqref{eq:Domain1}. 
\item
The complex eigenvalues $(\mu_m)_{m \in \nM}$ associated with the included Fourier functions $(\nF_m)_{m \in \nM}$ are computed and stored, see~\eqref{eq:Fourier1} and~\eqref{eq:Expansion1}.
\item
The basic integrals~\eqref{eq:Integral1a} are computed and stored.  
\item
The values of the integral kernel at the grid points are computed. 
Based on them, the values of a suitable regularisation on the previously defined neighbourhood of the origin are determined by interpolation.  
On the one hand, the associated spectral coefficients $(\psi_m)_{m \in \nM}$ are computed through a fast Fourier transform, see~\eqref{eq:Expansion1}.
Auxiliary quantities involving the basic integrals~\eqref{eq:Integral1a} and complex exponentials~\eqref{eq:Fourier1} are computed through summations along certain directions, in total three single sums and three double sums. 
On the other hand, the non-zero values corresponding to the difference between the singular integral kernel and its interpolant 
\begin{equation*}
\big(\varphi(v_{\ell}) - \psi(v_{\ell})\big)_{\ell \in \widetilde{\nL}}\,, 
\end{equation*}
are computed, see~\eqref{eq:SetSmall1}.   
Based on these quantities, quadrature approximations to the in total $6 \, M^3$ integrals
\begin{equation*}
\begin{gathered}
\int_{[- \, b_0, b_0]^3} w_i \, w_j \, (\varphi - \psi)(w) \, \nF_{\!m}(w) \; \dd w\,, \\
(i, j) \in \{(1, 1), (1, 2), (1, 3), (2, 2), (2, 3), (3, 3)\}\,, \quad m \in \nM\,, 
\end{gathered}
\end{equation*}
are computed and stored, see~\eqref{eq:Integral1b}.
\end{enumerate}

\Myparagraph{Computations} 
From the values of the initially prescribed density function at the velocity grid points, the associated spectral coefficients are computed through a fast Fourier transform 
\begin{equation*}
n = 1\,, \quad \big(f^{(n-1)}(v_{\ell})\big)_{\ell \in \nL}\,, \quad \big(f^{(n-1)}_m\big)_{m \in \nM}\,, 
\end{equation*}
yielding the approximation
\begin{equation*}
n = 1\,, \quad \sum_{m \in \nM} f^{(n-1)}_m \, \nF_m \approx f^{(n-1)}\,,
\end{equation*}
see also~\eqref{eq:Expansion1}.
In each substep of the time integration, based for instance on the explicit Euler method~\eqref{eq:ExplicitEuler}, the following computations are carried out, in order to determine the values and spectral coefficients of the discrete solution.
\begin{enumerate}[(i)]
\item
Numerical approximations to the values of the gradient at the velocity grid points are computed through pointwise multiplications by the corresponding eigenvalues and three inverse fast Fourier transforms
\begin{equation*}
\sum_{m \in \nM} \mu_m \, f^{(n-1)}_m \, \nF_m(v_{\ell}) \approx \MyGrad{v} f^{(n-1)}(v_{\ell})\,, \quad \ell \in \nL\,, 
\end{equation*}
see~\eqref{eq:FormulaQc1}.  
\item 
The main computational cost for the numerical approximation of fundamental integrals, which involve the regularised kernel, amounts to pointwise multiplications and in total 18 inverse fast Fourier transforms, see~\eqref{eq:Integral1}.
\item 
The representation~\eqref{eq:FormulaQc1} is used to determine the values of the integral operator 
\begin{equation*}
\nQc\big(f^{(n-1)}, f^{(n-1)}\big)(v_{\ell})\,, \quad \ell \in \nL\,.
\end{equation*}
In each case, the divergence is computed by three fast Fourier transforms and an additional inverse fast Fourier transform
\begin{equation}
\label{eq:GradQc1}
\MyGrad{v} \cdot \nQc\big(f^{(n-1)}, f^{(n-1)}\big)(v_{\ell})\,, \quad \ell \in \nL\,.
\end{equation}
Simple summations finally yield the values and spectral coefficients of the new approximation 
\begin{equation*}
\big(f^{(n)}(v_{\ell})\big)_{\ell \in \nL}\,, \quad \big(f^{(n)}_m\big)_{m \in \nM}\,,
\end{equation*}
see also~\eqref{eq:ExplicitEuler}.
\end{enumerate}

\Myparagraph{Output}
\begin{enumerate}[(i)]
\item 
By means of the above procedure, we obtain the values and the spectral coefficients of the discrete solution 
\begin{equation*}
\big(f^{(n)}(v_{\ell})\big)_{\ell \in \nL}\,, \quad \big(f^{(n)}_m\big)_{m \in \nM}\,, \quad n \in \{0, 1, \dots, N\}\,. 
\end{equation*}
\item 
If desired, approximations at intermediate points can be computed through the relation 
\begin{equation*}
f^{(n-1)}(v) \approx \sum_{m \in \nM} f^{(n-1)}_m \, \nF_m(v)\,, \quad v \in [- \, b, b]^3\,, \quad n \in \{0, 1, \dots, N\}\,. 
\end{equation*}
\end{enumerate}

\MyParagraph{Computational cost}
In the previously described situation, we may consider the 26 Fourier transforms for the computation of the Landau operator in a substep of the time integration as the computationally most elaborate components.
Other processes such as summations and pointwise multiplications can be optimised by parallelisation. 
We point out that the choice of the neighbourhood, where a regularisation of the singular integral kernel is determined, is crucial, see Figure~\ref{fig:Kernel}.  
Compared to a quadrature approximation on the whole domain, a suitable adjustment of the relatively small subset makes it possible to significantly reduce the precomputation time for the same accuracy, see Table~\ref{tab:Table2}.

\MyParagraph{Implementation and improvements}
In order to perform numerical comparisons and to design graphical illustrations, we found it convenient to implement our approaches in \textsc{Matlab}.
An elementary code that has the purpose to illustrate the practical implementation and reproduces numerical results discussed in Section~\ref{subsec:NumericalResults} is available through~\cite{CM2024Code}.
As we prioritised readability and the validation of common components, our code comprises several subroutines, which lower speed.
We thus consider the observed computation times as rough indicators and see opportunities for improvements based on efficient software packages.
Detailed derivations in a more general setting and additional numerical comparisons for different test examples with known solutions are described in the subsequent sections, see also 
Table~\ref{tab:Table1} for an overview.

\MyParagraph{Conservation of mass}
A characteristical property of the solution to the Landau equation~\eqref{eq:LandauEquation} is that it conserves mass, momentum, and energy and that entropy decays over time 
\begin{equation}
\label{eq:MassMomentumEnergyEntropy}
\begin{gathered}
\int_{\Omega^{(v)}} f(v, t) \; \dd v = \int_{\Omega^{(v)}} f(v, t_0) \; \dd v\,, \\
\int_{\Omega^{(v)}} v \, f(v, t) \; \dd v = \int_{\Omega^{(v)}} v \, f(v, t_0) \; \dd v\,, \\
\int_{\Omega^{(v)}} \abs{v}^2 \, f(v, t) \; \dd v = \int_{\Omega^{(v)}} \abs{v}^2 \, f(v, t_0) \; \dd v\,, \\
\int_{\Omega^{(v)}} f(v, t) \, \ln\big(f(v, t)\big) \; \dd v \leq \int_{\Omega^{(v)}} f(v, t_0) \, \ln\big(f(v, t_0)\big) \; \dd v\,, \\
t \in [t_0, T]\,.
\end{gathered}
\end{equation}
Due to the employed representation of the Landau operator in divergence form, the conservation of mass is ensured for the discrete solution as well, provided that the mass is computed subsequently to~\eqref{eq:GradQc1} through the trivial identity for the associated spectral coefficients
\begin{equation*}
f^{(n)}_0 = f^{(n-1)}_0\,, \quad n \in \{1, \dots, N\}. 
\end{equation*}
Any additional application of a fast Fourier transform or inverse fast Fourier transform, however, will cause numerical perturbations, since the values of a regular function are not retained, in general, i.e.
\begin{equation*}
\text{IFFT}\big(\text{FFT}(g)\big) \neq g\,,
\end{equation*}
but we may expect highly accurate approximations.
Likewise, we cannot ensure the preservation of momentum and energy as well as the strict positivity of the density function and the decay of entropy in case the reconstruction is positive, but we may expect highly accurate numerical approximations for suitable discretisations based on the Fourier spectral method and geometric time integrators, see for instance \cite{HairerLubichWanner2006,Iserles2008,McLachlanQuispel2002,SanzSernaCalvo2018}.
%%%%%%%%%%%%%%%%%%%%%%%%%%%%%%%%%%%%%%%%%%%%%%%%%%%%%%%%%%%%%%%%%%%%%%%%%%%%%%%%%%%%%%%%%%%%%%%%%%%%%%%%%%%%%%%%%%%%%%%%%%%%%%%%%%%%%%%%%%%%%%%%%% 
\subsection{Detailed description}
%%%%%%%%%%%%%%%%%%%%%%%%%%%%%%%%%%%%%%%%%%%%%%%%%%%%%%%%%%%%%%%%%%%%%%%%%%%%%%%%%%%%%%%%%%%%%%%%%%%%%%%%%%%%%%%%%%%%%%%%%%%%%%%%%%%%%%%%%%%%%%%%%% 
\MyParagraph{Compact exposition}
For the benefit of compact representations, we employ convenient abbreviations for the integral kernel and the term involving the outer product
\begin{subequations}
\label{eq:LandauOperator}
\begin{equation}
\begin{gathered}
\varphi: \Omega^{(\varphi)} \longrightarrow \RR: v \longmapsto C \, \abs{v}^{\beta} \, \widetilde{\varphi}(v)\,, \quad
\Omega^{(\varphi)} = \begin{cases} \RR^d & \!\!\!\text{if } \, \beta \geq 0\,, \\ \RR^d \setminus \{0\} & \!\!\!\text{if } \, \beta < 0\,, \end{cases} \\
P: \RR^d \longrightarrow \RR^{d \times d}: v \longmapsto \abs{v}^2 \, I_d - v \otimes v\,.
\end{gathered}
\end{equation}
In order to reflect the occurence of singular kernels and to cover instances such as~\eqref{eq:KernelCoulombGeneral} with $\beta = - \, 1$ and $\widetilde{\varphi}(v) = \ee^{\gamma \, \abs{v}}$, we here introduce auxiliary exponents $\beta \in \RR$, $\gamma \leq 0$, and a regular function $\widetilde{\varphi}: \RR \to \RR$.
When applicable, we omit the dependence of the density function, the integral operator, and the Landau operator on the time variable, i.e., we write 
\begin{equation}
\begin{gathered}
\nQc(f, f)(v) = \int_{\Omega} (\varphi \, P)(v - w) \, \big(\MyGrad{v} \, f(v) \, f(w) - f(v) \, \MyGrad{w} \, f(w)\big) \; \dd w\,, \\
\nQ(f, f)(v) = \MyGrad{v} \cdot \nQc(f, f)(v)\,, \\
\quad v \in \Omega = \Omega^{(v)}\,, 
\end{gathered}
\end{equation}
\end{subequations}
for short.
For the purpose of numerical validation and comparison, we distinguish the cases of constant, regular, and singular integral kernels and adapt our procedure accordingly to the increasing complexity of the problem. 
We exemplify calculations for two dimensions and state the analogous results for the significantly more involved three-dimensional case.
%%%%%%%%%%%%%%%%%%%%%%%%%%%%%%%%%%%%%%%%%%%%%%%%%%%%%%%%%%%%%%%%%%%%%%%%%%%%%%%%%%%%%%%%%%%%%%%%%%%%%%%%%%%%%%%%%%%%%%%%%%%%%%%%%%%%%%%%%%%%%%%%%% 
\subsection{Approaches}
\label{sec:Approaches}
%%%%%%%%%%%%%%%%%%%%%%%%%%%%%%%%%%%%%%%%%%%%%%%%%%%%%%%%%%%%%%%%%%%%%%%%%%%%%%%%%%%%%%%%%%%%%%%%%%%%%%%%%%%%%%%%%%%%%%%%%%%%%%%%%%%%%%%%%%%%%%%%%% 
In this section and the appendix, we outline different approaches for the numerical evaluation of the Landau operator~\eqref{eq:LandauOperator}, see Table~\ref{tab:Table1}.
We begin with detailed considerations for a natural procedure and then indicate possible alternatives. 
Useful means are suitable reformulations of the integral operator obtained by Fourier series expansions of the density function and a regularised integral kernel, the application of a linear integral transform and a quadrature rule on a relatively small neighbourhood of the origin, and subsequent differentiation.
In order to underline the sources of errors due to the truncation of the unbounded velocity domain, the replacement of infinite by finite sums, and quadrature, we first state the employed representations of the Landau operator on the basis of formal series expansions and then sketch the corresponding algorithms. 
We recall that auxiliary abbreviations and results related to the Fourier spectral method are found in Section~\ref{sec:Auxiliaries}.
For the benefit of shorter formulas, we henceforth set
\begin{equation}
\label{eq:nM123}
\begin{gathered}
d = 2\,, \\
\nM_1 = \{(1, 1), (1, 2), (2, 2)\}\,, \\
\nM_2 = \{(1, 1, 2), (1, 2, 1), (1, 2, 2), (2, 2, 1)\}\,, \\
\nM_3 = \{(0, 0), (1, 0), (0, 1), (2, 0), (1, 1), (0, 2)\}\,, \\  
d = 3\,, \\
\nM_1 = \{(1, 1), (1, 2), (1, 3), (2, 2), (2, 3), (3, 3)\}\,, \\
\nM_2 = \{(1, 1, 2), (1, 1, 3), (1, 2, 1), (1, 2, 2), (1, 3, 1), (1, 3, 3), \\
\qquad\qquad\; (2, 2, 1), (2, 2, 3), (2, 3, 2), (2, 3, 3), (3, 3, 1), (3, 3, 2)\}\,, \\
\nM_3 = \{(0, 0, 0), (1, 0, 0), (0, 1, 0), (0, 0, 1), (2, 0, 0), (1, 1, 0), \\
(1, 0, 1), (0, 2, 0), (0, 1, 1), (0, 0, 2)\}\,, \qquad\quad  
\end{gathered}
\end{equation}
%%%%%%%%%%%%%%%%%%%%%%%%%%%%%%%%%%%%%%%%%%%%%%%%%%%%%%%%%%%%%%%%%%%%%%%%%%%%%%%%%%%%%%%%%%%%%%%%%%%%%%%%%%%%%%%%%%%%%%%%%%%%%%%%%%%%%%%%%%%%%%%%%% 
\subsection{Approach~CST2 (Conservative form Singular kernel Transform twice)}
\label{sec:SectionCST2}
%%%%%%%%%%%%%%%%%%%%%%%%%%%%%%%%%%%%%%%%%%%%%%%%%%%%%%%%%%%%%%%%%%%%%%%%%%%%%%%%%%%%%%%%%%%%%%%%%%%%%%%%%%%%%%%%%%%%%%%%%%%%%%%%%%%%%%%%%%%%%%%%%% 
\MyParagraph{Integral operator}
Our starting point is the following representation of the integral operator 
\begin{equation*}
\nQc(f, f)(v) = \nI(f)(v) \, \MyGrad{v} f(v) + \nJ(f)(v) \, f(v)\,, \quad v \in \Omega\,.
\end{equation*}
The arising matrix- and vector-valued operators
\begin{equation*}
\begin{gathered}
d = 2\,, \\
\nI = \begin{pmatrix} \nI_{22} & - \, \nI_{12} \\ - \, \nI_{12} & \nI_{11} \end{pmatrix}\,, \quad 
\nJ = \begin{pmatrix} - \, \nJ_{221} + \nJ_{122} \\ \nJ_{121} - \nJ_{112} \end{pmatrix}\,, \\
d = 3\,, \\
\nI =
\begin{pmatrix}
\nI_{22} + \nI_{33} & - \, \nI_{12} & - \, \nI_{13} \\
- \, \nI_{12} & \nI_{11} + \nI_{33} & - \, \nI_{23} \\
- \, \nI_{13} & - \, \nI_{23} & \nI_{11} + \nI_{22}
\end{pmatrix}\,, \\
\nJ =
\begin{pmatrix}
- \, \nJ_{221} - \nJ_{331} + \nJ_{122} + \nJ_{133} \\ \nJ_{121} - \nJ_{112} - \nJ_{332} + \nJ_{233} \\ \nJ_{131} + \nJ_{232} - \nJ_{113} - \nJ_{223} \end{pmatrix}\,, \\
\end{gathered}
\end{equation*}
are defined by fundamental integrals, which involve polynomials of degree two, the singular kernel, and the density function  
\begin{equation*}
\begin{gathered}
\nI_{ij}(f)(v) = \int_{\Omega} (v_i - w_i) \, (v_j - w_j) \, \varphi(v - w) \, f(w) \; \dd w\,, \quad (i, j) \in \nM_1\,, \\
\nJ_{ijk}(f)(v) = \int_{\Omega} (v_i - w_i) \, (v_j - w_j) \, \varphi(v - w) \, \partial_{w_k} f(w) \; \dd w\,, \quad (i, j, k) \in \nM_2\,, \\
v \in \Omega\,.
\end{gathered}
\end{equation*}
For the sake of concreteness, we next explain our main strategy on the basis of the two instances 
\begin{equation}
\label{eq:I11}  
\begin{gathered}
\nI_{11}(f)(v) = \int_{\Omega} (v_1 - w_1)^2 \, \varphi(v - w) \, f(w) \; \dd w\,, \\
\nJ_{112}(f)(v) = \int_{\Omega} (v_1 - w_1)^2 \, \varphi(v - w) \, \partial_{w_2} f(w) \; \dd w\,, \\
v \in \Omega\,, 
\end{gathered}
\end{equation}
and then state the corresponding representations obtained in the general case. 

\MyParagraph{Singular integral kernel}
Our basic idea is to exploit the evident identity 
\begin{subequations}
\label{eq:KernelDecomposition} 
\begin{equation}
\varphi(v) = \psi(v) + (\varphi - \psi)(v)\,, \quad v \in \Omega\,,
\end{equation}  
and to adjust the regular function~$\psi$ such that the remainder $\varphi - \psi$ vanishes on the main part of the velocity domain.
More concretely, we make use of the fact that the singular integral kernel is regular outside a relatively small neighbourhood of the origin
\begin{equation}
\Omo = [b_{11}^{(0)}, b_{12}^{(0)}] \times \cdots \times [b_{d1}^{(0)}, b_{d2}^{(0)}] \subset \Omb \subset \Omega\,,
\end{equation}
see also~\eqref{eq:Omb}. 
Hence, excluding this set, the kernel defines a regular function
\begin{equation}
\psi(v) = \varphi(v)\,, \quad v \in \Omega \setminus \Omo\,.
\end{equation}
\end{subequations}
Nearby the origin, straightforward interpolation of~$\varphi$ is applied to determine~$\psi$.

\MyParagraph{Fourier series expansions}
For the regularised integral kernel and the density function, we may assume that favourable approximations are provided by Fourier series expansions.
We meanwhile employ the formal representations 
\begin{equation}
\label{eq:FourierSeriesExpansions}  
\begin{gathered}
\psi(v) = \sum_{\ell \in \ZZ^d} \psi_{\ell} \, \nFb_{\ell}(v)\,, \\
f(v) = \sum_{m \in \ZZ^d} f_m \, \nFb_m(v)\,, \quad
\partial_{v_k} f(v) = \sum_{m \in \ZZ^d} \mu^{(b_{k1}, \, b_{k2})}_{m_k} \, f_m \, \nFb_m(v)\,, \\ 
v \in \Omega\,, \quad k \in \{1, \dots, d\}\,.
\end{gathered}
\end{equation}
Their practical implementation presupposes suitable truncations of the unbounded velocity domain to avoid significant aliasing effects, see for instance~\cite{PareschiRusso2000,PareschiEtAl2000}.
Moreover, it requires truncations of the infinite sums as well as the application of the trapezoidal rule for the numerical computation of the spectral coefficients.

\MyParagraph{Decisive integrals}
The above stated Fourier series expansions of the density function and its partial derivatives imply representations such as 
\begin{equation*}
\begin{gathered}
\nI_{11}(f)(v) = \sum_{m \in \ZZ^d} f_m \; \nIm_{11}(v)\,, \\
\nJ_{112}(f)(v) = \sum_{m \in \ZZ^d} \mu^{(b_{21}, \, b_{22})}_{m_2} \, f_m \; \nIm_{11}(v)\,, \\
\nIm_{11}(v) = \int_{\Omega} (v_1 - w_1)^2 \, \varphi(v - w) \, \nFb_m(w) \; \dd w\,, \quad m \in \ZZ^d\,, \\
v \in \Omega\,, 
\end{gathered}
\end{equation*}
see~\eqref{eq:I11} and~\eqref{eq:FourierSeriesExpansions}.
Due to the decomposition of the kernel into a regular function and a singular function that vanishes on the main part of the domain, the decisive integrals comprise the two contributions 
\begin{equation}
\label{eq:Ipsi_Iphipsi}
\begin{gathered}
\nIm_{11}(v) = \nImpsi_{11}(v) + \nImphipsi_{11}(v)\,, \\
\nImpsi_{11}(v) = \int_{\Omega} (v_1 - w_1)^2 \, \psi(v - w) \, \nFb_m(w) \; \dd w\,, \\
\nImphipsi_{11}(v) = \int_{\Omega} (v_1 - w_1)^2 \, (\varphi - \psi)(v - w) \, \nFb_m(w) \; \dd w\,, \\
m \in \ZZ^d\,, \quad v \in \Omega\,, 
\end{gathered}
\end{equation}
see also~\eqref{eq:KernelDecomposition}.

\MyParagraph{Linear integral transform}
Applying in both cases the linear integral transform $u = v - w$, yields the equivalent reformulations 
\begin{equation*}
\begin{gathered}
\nImpsi_{11}(v) = \int_{v - \Omega} u_1^2 \; \psi(u) \, \nFb_m(v - u) \; \dd u\,, \\
\nImphipsi_{11}(v) = \int_{v - \Omega} u_1^2 \, (\varphi - \psi)(u) \, \nFb_m(v - u) \; \dd u\,, \\
m \in \ZZ^d\,, \quad v \in \Omega\,, 
\end{gathered}
\end{equation*}
where we employ a symbolic notation for the shifted domain comprising an additional sign owing to $\dd u_k = - \, \dd w_k$ for $k \in \{1, \dots, d\}$.
By means of the Fourier series expansion of the regularised kernel~\eqref{eq:FourierSeriesExpansions} and the identity~\eqref{eq:FourierIdentity2}, the first multiple integral reduces to one-dimensional integrals 
\begin{equation*}
\nImpsi_{11}(v) = \nEb_{-m}(b_1) \, \nFb_{m}(v) \sum_{\ell \in \ZZ^d} \psi_{\ell} \int_{v - \Omega} u_1^2 \; \nFb_{\ell - m}(u) \; \dd u\,, \quad
m \in \ZZ^d\,, \quad v \in \Omega\,.
\end{equation*}
For the second multiple integral, we instead apply the relation~\eqref{eq:FourierIdentity1} and have 
\begin{equation*}
\nImphipsi_{11}(v) = \nFb_m(v) \int_{v - \Omega} u_1^2 \; (\varphi - \psi)(u) \, \nEb_{- m}(u) \; \dd u\,, \quad
m \in \ZZ^d\,, \quad v \in \Omega\,.
\end{equation*}

\MyParagraph{Approach~CST2}
Summarising the above procedure, we obtain the following relations for the derivatives of the density function
\begin{subequations}
\label{eq:ApproachCST2}
\begin{equation}
\partial_{v_k} f(v) = \sum_{m \in \ZZ^d} \mu^{(b_{k1}, \, b_{k2})}_{m_k} \, f_m \, \nFb_m(v)\,, \quad k \in \{1, \dots, d\}\,, \quad v \in \Omega\,,
\end{equation}
the fundamental integrals 
\begin{equation}
\begin{gathered}
\nImpsi_{ij}(v) = \nEb_{-m}(b_1) \, \nFb_{m}(v) \sum_{\ell \in \ZZ^d} \psi_{\ell} \int_{v - \Omega} u_i \, u_j \; \nFb_{\ell - m}(u) \; \dd u\,, \\
\nImphipsi_{ij}(v) = \nFb_m(v) \int_{v - \Omega} u_i \, u_j \; (\varphi - \psi)(u) \, \nEb_{- m}(u) \; \dd u\,, \\
\nIm_{ij}(v) = \nImpsi_{ij}(v) + \nImphipsi_{ij}(v)\,, \\
\nI_{ij}(f)(v) = \sum_{m \in \ZZ^d} f_m \; \nIm_{ij}(v)\,, \quad (i, j) \in \nM_1\,, \\
\nJ_{ijk}(f)(v) = \sum_{m \in \ZZ^d} \mu^{(b_{k1}, \, b_{k2})}_{m_k} \, f_m \; \nIm_{ij}(v)\,, \quad (i, j, k) \in \nM_2\,, \\
v \in \Omega\,,
\end{gathered}
\end{equation}
and the integral operator 
\begin{equation}
\begin{gathered}
\nQc(f, f)(v) = \nI(f)(v) \, \MyGrad{v} f(v) + \nJ(f)(v) \, f(v)\,, \\
d = 2\,, \\
\nI = \begin{pmatrix} \nI_{22} & - \, \nI_{12} \\ - \, \nI_{12} & \nI_{11} \end{pmatrix}\,, \quad 
\nJ = \begin{pmatrix} - \, \nJ_{221} + \nJ_{122} \\ \nJ_{121} - \nJ_{112} \end{pmatrix}\,, \\
d = 3\,, \\
\nI =
\begin{pmatrix}
\nI_{22} + \nI_{33} & - \, \nI_{12} & - \, \nI_{13} \\
- \, \nI_{12} & \nI_{11} + \nI_{33} & - \, \nI_{23} \\
- \, \nI_{13} & - \, \nI_{23} & \nI_{11} + \nI_{22}
\end{pmatrix}\,, \\
\nJ =
\begin{pmatrix}
- \, \nJ_{221} - \nJ_{331} + \nJ_{122} + \nJ_{133} \\ \nJ_{121} - \nJ_{112} - \nJ_{332} + \nJ_{233} \\ \nJ_{131} + \nJ_{232} - \nJ_{113} - \nJ_{223} \end{pmatrix}\,, \\
v \in \Omega\,.
\end{gathered}
\end{equation}
Furthermore, in order to determine the Landau operator, we use the associated Fourier series representation 
\begin{equation}
\begin{gathered}
\nQc(f, f)(v) = \sum_{m \in \ZZ^d} \nQc_m \, \nFb_m(v)\,, \\
\nQ(f, f)(v) = \MyGrad{v} \cdot \nQc(f, f)(v) = \sum_{m \in \ZZ^d} \mub_m \, \nQc_m \, \nFb_m(v)\,, \\
\quad v \in \Omega\,.
\end{gathered}
\end{equation}
\end{subequations}

\MyParagraph{Practical implementation}
The practical implementation of~\eqref{eq:ApproachCST2} requires suitable truncations of the velocity domain and the infinite sums.
We point out that replacing the shifted domain $v - \Omega$ by the truncated integral domain~$\Omb$ permits the application of fast Fourier techniques.
We recall that the essential steps of the algorithm, distinguishing between precomputations that are independent of the density function and computations that are carried out repeatedly in course of the time integration of the Landau equation, are described in Section~\ref{sec:Introduction}.
Moreover, a link to an elementary \textsc{Matlab} code is provided there.
%%%%%%%%%%%%%%%%%%%%%%%%%%%%%%%%%%%%%%%%%%%%%%%%%%%%%%%%%%%%%%%%%%%%%%%%%%%%%%%%%%%%%%%%%%%%%%%%%%%%%%%%%%%%%%%%%%%%%%%%%%%%%%%%%%%%%%%%%%%%%%%%%% 
\subsection{Approach~CST1 (Conservative form Singular kernel Transform once)}
\label{sec:SectionCST1}
%%%%%%%%%%%%%%%%%%%%%%%%%%%%%%%%%%%%%%%%%%%%%%%%%%%%%%%%%%%%%%%%%%%%%%%%%%%%%%%%%%%%%%%%%%%%%%%%%%%%%%%%%%%%%%%%%%%%%%%%%%%%%%%%%%%%%%%%%%%%%%%%%% 
\MyParagraph{Modification}
Regarding a possible modification of our first approach, we reconsider the quantity 
\begin{equation*}
\begin{gathered}
\nImpsi_{ij}(v) = \int_{\Omega} (v_i - w_i) \, (v_j - w_j) \, \psi(v - w) \, \nFb_m(w) \; \dd w\,, \\
(i, j) \in \nM_1\,, \quad m \in \ZZ^d\,, \quad v \in \Omega\,, 
\end{gathered}
\end{equation*}
see also~\eqref{eq:Ipsi_Iphipsi}.
Expanding the integrand, inserting the Fouries series representation of the regularised kernel~\eqref{eq:FourierSeriesExpansions}, and applying the identity~\eqref{eq:FourierIdentity3} yields the alternative reformulation
\begin{subequations}
\label{eq:nImpsi}  
\begin{equation}
\begin{gathered}
\widetilde{I}_{ij}^{(m - \ell)}(v) = \int_{\Omega} (v_i \, v_j + v_i \, w_j + v_j \, w_i + w_i \, w_j) \, \nFb_{m - \ell}(w) \; \dd w\,, \\
\nImpsi_{ij}(v) = \sum_{\ell \in \ZZ^d} \psi_{\ell} \, \nEb_{- \ell}(b_1) \, \nFb_{\ell}(v) \, \widetilde{I}_{ij}^{(m - \ell)}(v)\,, \\
(i, j) \in \nM_1\,, \quad m \in \ZZ^d\,, \quad v \in \Omega\,.
\end{gathered}
\end{equation}
Employing the short notation
\begin{equation}
\widetilde{I}(k, m) = \int_{\Omega} w^k \, \nFb_{m}(w) \; \dd w\,, \quad k \in \nM_3\,, \quad m \in \ZZ^d\,, 
\end{equation}
the two-dimensional case takes the form 
\begin{equation}
\begin{gathered}
\widetilde{I}_{11}^{(m)}(v) = v_1^2 \, \widetilde{I}\big((0, 0), m\big) + 2 \, v_1 \, \widetilde{I}\big((1, 0), m\big) + \widetilde{I}\big((2, 0), m\big)\,, \\
\widetilde{I}_{22}^{(m)}(v) = v_2^2 \, \widetilde{I}\big((0, 0), m\big) + 2 \, v_2 \, \widetilde{I}\big((0, 1), m\big) + \widetilde{I}\big((0, 2), m\big)\,, \\
\widetilde{I}_{12}^{(m)}(v) = v_1 \, v_2 \, \widetilde{I}\big((0, 0), m\big) + v_1 \, \widetilde{I}\big((0, 1), m\big) + v_2 \, \widetilde{I}\big((1, 0), m\big) + \widetilde{I}\big((1, 1), m\big)\,.
\end{gathered}
\end{equation}
In three dimensions, we instead have 
\begin{equation}
\begin{gathered}
\widetilde{I}_{11}^{(m)}(v) = v_1^2 \, \widetilde{I}\big((0, 0, 0), m\big) + 2 \, v_1 \, \widetilde{I}\big((1, 0, 0), m\big) + \widetilde{I}\big((2, 0, 0), m\big)\,, \\
\widetilde{I}_{22}^{(m)}(v) = v_2^2 \, \widetilde{I}\big((0, 0, 0), m\big) + 2 \, v_2 \, \widetilde{I}\big((0, 1, 0), m\big) + \widetilde{I}\big((0, 2, 0), m\big)\,, \\
\widetilde{I}_{33}^{(m)}(v) = v_3^2 \, \widetilde{I}\big((0, 0, 0), m\big) + 2 \, v_3 \, \widetilde{I}\big((0, 0, 1), m\big) + \widetilde{I}\big((0, 0, 2), m\big)\,, \\
\widetilde{I}_{12}^{(m)}(v) = v_1 \, v_2 \, \widetilde{I}\big((0, 0, 0), m\big) + v_1 \, \widetilde{I}\big((0, 1, 0), m\big) + v_2 \, \widetilde{I}\big((1, 0, 0), m\big)
+ \widetilde{I}\big((1, 1, 0), m\big)\,, \\
\widetilde{I}_{13}^{(m)}(v) = v_1 \, v_3 \, \widetilde{I}\big((0, 0, 0), m\big) + v_1 \, \widetilde{I}\big((0, 0, 1), m\big) + v_3 \, \widetilde{I}\big((1, 0, 0), m\big) 
+ \widetilde{I}\big((1, 0, 1), m\big)\,, \\
\widetilde{I}_{23}^{(m)}(v) = v_2 \, v_3 \, \widetilde{I}\big((0, 0, 0), m\big) + v_2 \, \widetilde{I}\big((0, 0, 1), m\big) + v_3 \, \widetilde{I}\big((0, 1, 0), m\big) 
+ \widetilde{I}\big((0, 1, 1), m\big)\,.
\end{gathered}
\end{equation}
\end{subequations}

\MyParagraph{Approach~CST1}
For the sake of completeness, we recapitulate the resulting representations for the derivatives of the density function
\begin{subequations}
\label{eq:ApproachCST1}
\begin{equation}
\partial_{v_k} f(v) = \sum_{m \in \ZZ^d} \mu^{(b_{k1}, \, b_{k2})}_{m_k} \, f_m \, \nFb_m(v)\,, \quad k \in \{1, \dots, d\}\,, \quad v \in \Omega\,,
\end{equation}
the basic integrals given by polynomials and Fourier functions 
\begin{equation}
\widetilde{I}(k, m) = \int_{\Omega} w^k \, \nFb_{m}(w) \; \dd w\,, \quad k \in \nM_3\,, \quad m \in \ZZ^d\,, \\
\end{equation}
the fundamental integrals involving the regularised kernel 
\begin{equation}
\begin{gathered}
\widetilde{I}_{ij}^{(m - \ell)}(v) = v_i \, v_j \int_{\Omega} \nFb_{m - \ell}(w) \; \dd w + v_i \int_{\Omega} w_j \, \nFb_{m - \ell}(w) \; \dd w \\
\qquad\qquad\qquad\quad + \; v_j \int_{\Omega} w_i \, \nFb_{m - \ell}(w) \; \dd w + \int_{\Omega} w_i \, w_j \, \nFb_{m - \ell}(w) \; \dd w\,, \\
\nImpsi_{ij}(v) = \sum_{\ell \in \ZZ^d} \psi_{\ell} \, \nEb_{- \ell}(b_1) \, \nFb_{\ell}(v) \, \widetilde{I}_{ij}^{(m - \ell)}(v)\,, \\
\nIpsi_{ij}(f)(v) = \sum_{m \in \ZZ^d} f_m \; \nImpsi_{ij}(v)\,, \quad (i, j) \in \nM_1\,, \\
\nJpsi_{ijk}(f)(v) = \sum_{m \in \ZZ^d} \mu^{(b_{k1}, \, b_{k2})}_{m_k} \, f_m \; \nImpsi_{ij}(v)\,, \quad (i, j,k) \in \nM_2\,, \\
v \in \Omega\,,
\end{gathered}
\end{equation}
as well as the fundamental integrals involving the difference of the singular kernel and its regularisation 
\begin{equation}
\begin{gathered}
\nImphipsi_{ij}(v) = \nFb_m(v) \int_{v - \Omega} u_i \, u_j \; (\varphi - \psi)(u) \, \nEb_{- m}(u) \; \dd u\,, \\
\nIphipsi_{ij}(f)(v) = \sum_{m \in \ZZ^d} f_m \; \nImphipsi_{ij}(v)\,, \quad (i, j) \in \nM_1\,, \\
\nJphipsi_{ijk}(f)(v) = \sum_{m \in \ZZ^d} \mu^{(b_{k1}, \, b_{k2})}_{m_k} \, f_m \; \nImphipsi_{ij}(v)\,, \quad (i, j,k) \in \nM_2\,, \\
v \in \Omega\,,
\end{gathered}
\end{equation}
the corresponding sums 
\begin{equation}
\begin{gathered}
\nI_{ij}(f)(v) = \nIpsi_{ij}(f)(v) + \nIphipsi_{ij}(f)(v)\,, \quad (i, j) \in \nM_1\,, \\
\nJ_{ijk}(f)(v) = \nJpsi_{ijk}(f)(v) + \nJphipsi_{ijk}(f)(v)\,, \quad (i, j, k) \in \nM_2\,, \\
v \in \Omega\,,
\end{gathered}
\end{equation}
and the integral operator 
\begin{equation}
\begin{gathered}
\nQc(f, f)(v) = \nI(f)(v) \, \MyGrad{v} f(v) + \nJ(f)(v) \, f(v)\,, \\
d = 2\,, \\
\nI = \begin{pmatrix} \nI_{22} & - \, \nI_{12} \\ - \, \nI_{12} & \nI_{11} \end{pmatrix}\,, \quad 
\nJ = \begin{pmatrix} - \, \nJ_{221} + \nJ_{122} \\ \nJ_{121} - \nJ_{112} \end{pmatrix}\,, \\
d = 3\,, \\
\nI =
\begin{pmatrix}
\nI_{22} + \nI_{33} & - \, \nI_{12} & - \, \nI_{13} \\
- \, \nI_{12} & \nI_{11} + \nI_{33} & - \, \nI_{23} \\
- \, \nI_{13} & - \, \nI_{23} & \nI_{11} + \nI_{22}
\end{pmatrix}\,, \\
\nJ =
\begin{pmatrix}
- \, \nJ_{221} - \nJ_{331} + \nJ_{122} + \nJ_{133} \\ \nJ_{121} - \nJ_{112} - \nJ_{332} + \nJ_{233} \\ \nJ_{131} + \nJ_{232} - \nJ_{113} - \nJ_{223} \end{pmatrix}\,, \\
v \in \Omega\,.
\end{gathered}
\end{equation}
Subsequent differentiation is again based on the Fourier series expansion 
\begin{equation}
\begin{gathered}
\nQc(f, f)(v) = \sum_{m \in \ZZ^d} \nQc_m \, \nFb_m(v)\,, \\
\nQ(f, f)(v) = \MyGrad{v} \cdot \nQc(f, f)(v) = \sum_{m \in \ZZ^d} \mub_m \, \nQc_m \, \nFb_m(v)\,, \\
\quad v \in \Omega\,.
\end{gathered}
\end{equation}
\end{subequations}

\MyParagraph{Practical implementation}
Concerning the practical implementation of~\eqref{eq:ApproachCST1}, we prescribe the dimension~$d$ of the velocity domain, the integral kernel~$\varphi$, and the value of the density function~$f$ at the initial time and the considered space grid points.
We replace the original integral domain~$\Omega$ by the truncated domain~$\Omb$ and make use of the fact that the integrals in~\eqref{eq:nImpsi} are related to the basic integrals given in~\eqref{eq:BasicIntegrals}.
It is again expedient to distinguish between precomputations and computations in course of the time integration, see also~\eqref{eq:ExplicitEuler}.

\Myparagraph{Input (Discretisation)}
We first specify 
\begin{enumerate}[(i)]
\item
the real numbers $(b_{i1}, b_{i2})_{i=1}^d$ with $b_{i1} < b_{i2}$ for $i \in \{1, \dots, d\}$ defining the truncated domain~$\Omb$, see~\eqref{eq:Omb}, 
\item
even positive integers to define the set~$\nM_M$ reflecting the total number of Fourier functions, see~\eqref{eq:nM},
\item
the real numbers $(b_{i1}^{(0)}, b_{i2}^{(0)})_{i=1}^d$ with $b_{i1}^{(0)} < b_{i2}^{(0)}$ for $i \in \{1, \dots, d\}$ defining the small restricted domain~$\Omo$, see~\eqref{eq:KernelDecomposition}.
\end{enumerate}

\Myparagraph{Precomputations}
Due to the fact that certain quantities do not depend on the values of the density function, it is advantageous to compute them in advance. 
This in particular concerns 
\begin{enumerate}[(i)]
\item
the in total $M_1 \cdots M_d$ equidistant grid points covering the truncated domain~$\Omb$,
\item
the eigenvalues $(\mub_m)_{m \in \nM_M}$ associated with the Fourier functions, see also~\eqref{eq:Fourier},
\item
the basic integrals~$I(k,m)$ for $k \in \nM_3$ and~$m \in \nM_M$, see~\eqref{eq:BasicIntegrals}, 
\item
the values of the integral kernel~$\varphi$ on the equidistant grid, excluding the singularity at the origin, 
\item
the values of the regularised kernel~$\psi$ on the grid points that are contained in the small domain~$\Omo$, obtained by interpolation, 
\item
the associated spectral coefficients~$(\psi_{\ell})_{\ell \in \nM_M}$ through a fast Fourier transform as well as the products $\psi_{\ell} \, \nEb_{- \ell}(b_1)$ for $\ell \in \nM_M$, see also~\eqref{eq:FourierSeriesExpansions},
\item
and quadrature approximations to the integrals
\begin{equation*}
\int_{\Omo} u_i \, u_j \, (\varphi - \psi)(u) \, \nEb_{- m}(u) \; \dd u\,, \quad (i, j) \in \nM_1\,, \quad m \in \nM_m\,, 
\end{equation*}
see~\eqref{eq:nM123} and~\eqref{eq:ApproachCST1}.
\end{enumerate}

\Myparagraph{Computations}
Regarding the evaluation of the Landau operator, we compute approximations to 
\begin{enumerate}[(i)]
\item
the spectral coefficients~$(f_m)_{m \in \nM_M}$ associated with the density function through a fast Fourier transform,
\item
the values of the derivatives~$\partial_{v_k} f$ for $k \in \{1, \dots, d\}$ at the prescribed equidistant grid points by means of the representations in~\eqref{eq:FourierSeriesExpansions}, requiring componentwise multiplications by eigenvalues and the application of fast inverse Fourier transforms,  
\item
on the one hand, the values of the following fundamental quantities at the grid points 
\begin{equation}
\begin{gathered}
\nIpsi_{ij}(f)(v) = \sum_{\ell \in \nM_M} \psi_{\ell} \; \nEb_{- \ell}(b_1) \, \nFb_{\ell}(v) \sum_{m \in \nM_M} f_m \, \widetilde{I}_{ij}^{(m - \ell)}(v) \,, \\
(i, j) \in \nM_1\,, \quad v \in \Omega\,,
\end{gathered}
\end{equation}
by summations along certain directions as well as fast inverse Fourier transforms fundamental integrals and accordingly for $\nJpsi_{ijk}(f)$ with $(i, j,k) \in \nM_2$, 
\item
on the other hand, the values of the fundamental quantities~$\nIphipsi_{ij}(f)$ for $(i, j) \in \nM_1$ and~$\nJphipsi_{ijk}(f)$ for $(i, j,k) \in \nM_2$ by fast inverse Fourier transforms,
\item
the products of their sums with the derivatives of the density function to obtain the components of the integral operator~$\nQc$,
\item   
and, in a final step, the divergence by multiplications with eigenvalues and fast inverse Fourier transforms. 
\end{enumerate}
%%%%%%%%%%%%%%%%%%%%%%%%%%%%%%%%%%%%%%%%%%%%%%%%%%%%%%%%%%%%%%%%%%%%%%%%%%%%%%%%%%%%%%%%%%%%%%%%%%%%%%%%%%%%%%%%%%%%%%%%%%%%%%%%%%%%%%%%%%%%%%%%%% 
%%%%%%%%%%%%%%%%%%%%%%%%%%%%%%%%%%%%%%%%%%%%%%%%%%%%%%%%%%%%%%%%%%%%%%%%%%%%%%%%%%%%%%%%%%%%%%%%%%%%%%%%%%%%%%%%%%%%%%%%%%%%%%%%%%%%%%%%%%%%%%%%%% 
\section{Numerical validations and comparisons}
\label{sec:Section3}
%%%%%%%%%%%%%%%%%%%%%%%%%%%%%%%%%%%%%%%%%%%%%%%%%%%%%%%%%%%%%%%%%%%%%%%%%%%%%%%%%%%%%%%%%%%%%%%%%%%%%%%%%%%%%%%%%%%%%%%%%%%%%%%%%%%%%%%%%%%%%%%%%% 
%%%%%%%%%%%%%%%%%%%%%%%%%%%%%%%%%%%%%%%%%%%%%%%%%%%%%%%%%%%%%%%%%%%%%%%%%%%%%%%%%%%%%%%%%%%%%%%%%%%%%%%%%%%%%%%%%%%%%%%%%%%%%%%%%%%%%%%%%%%%%%%%%% 
In this section, we are concerned with a thorough numerical validation and comparison of the approaches summarised in Table~\ref{tab:Table1}.
For this purpose, we next state test problems involving constant, regular, and singular integral kernels.
%%%%%%%%%%%%%%%%%%%%%%%%%%%%%%%%%%%%%%%%%%%%%%%%%%%%%%%%%%%%%%%%%%%%%%%%%%%%%%%%%%%%%%%%%%%%%%%%%%%%%%%%%%%%%%%%%%%%%%%%%%%%%%%%%%%%%%%%%%%%%%%%%% 
\subsection{Test problems}
\begin{figure}
\begin{center}
\includegraphics[width=6.6cm]{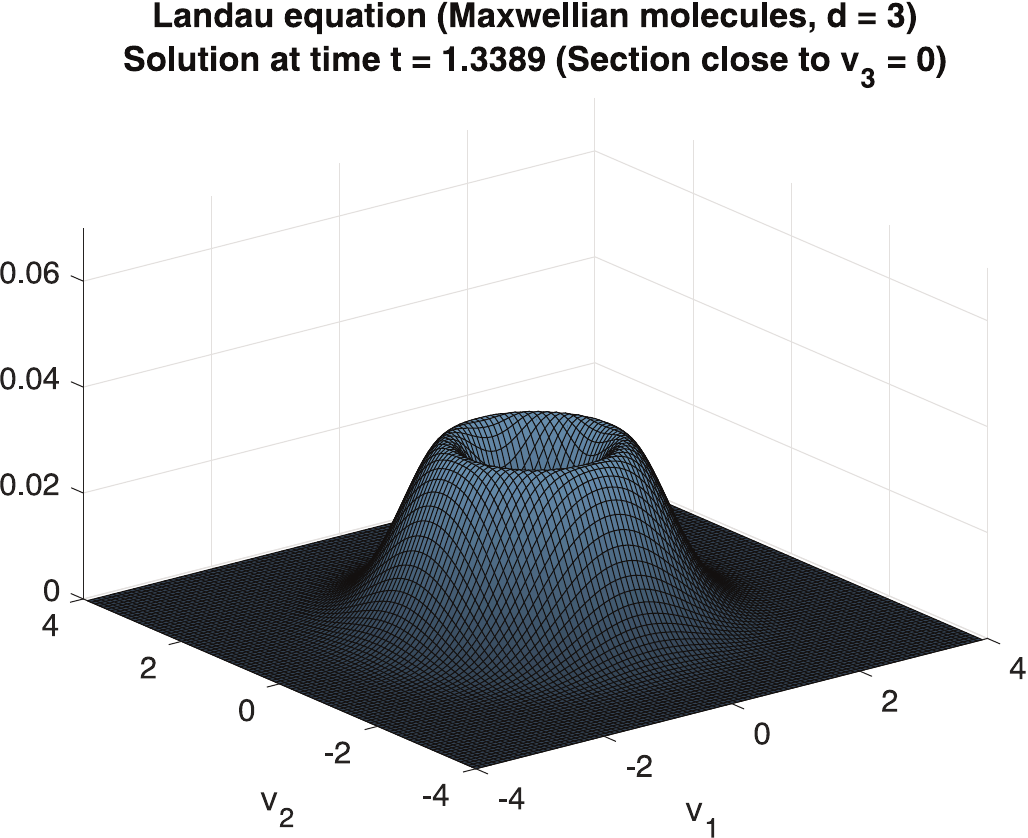}
\quad
\includegraphics[width=6.6cm]{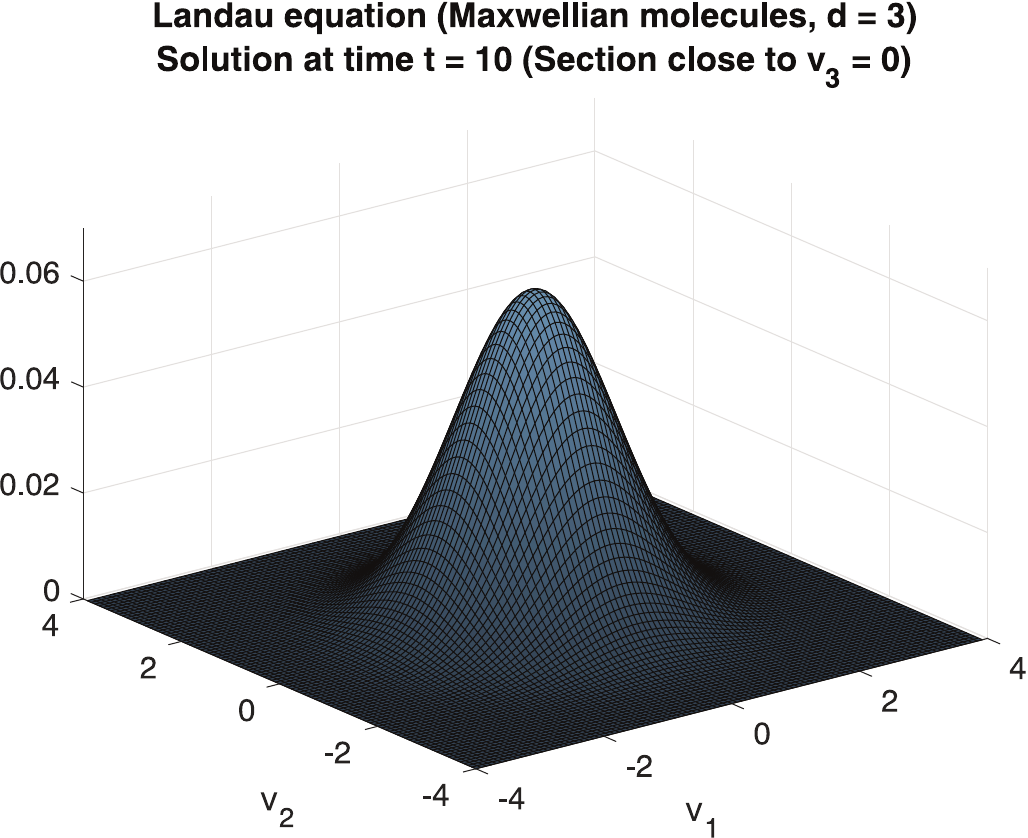}
\\[2mm]
\includegraphics[width=6.6cm]{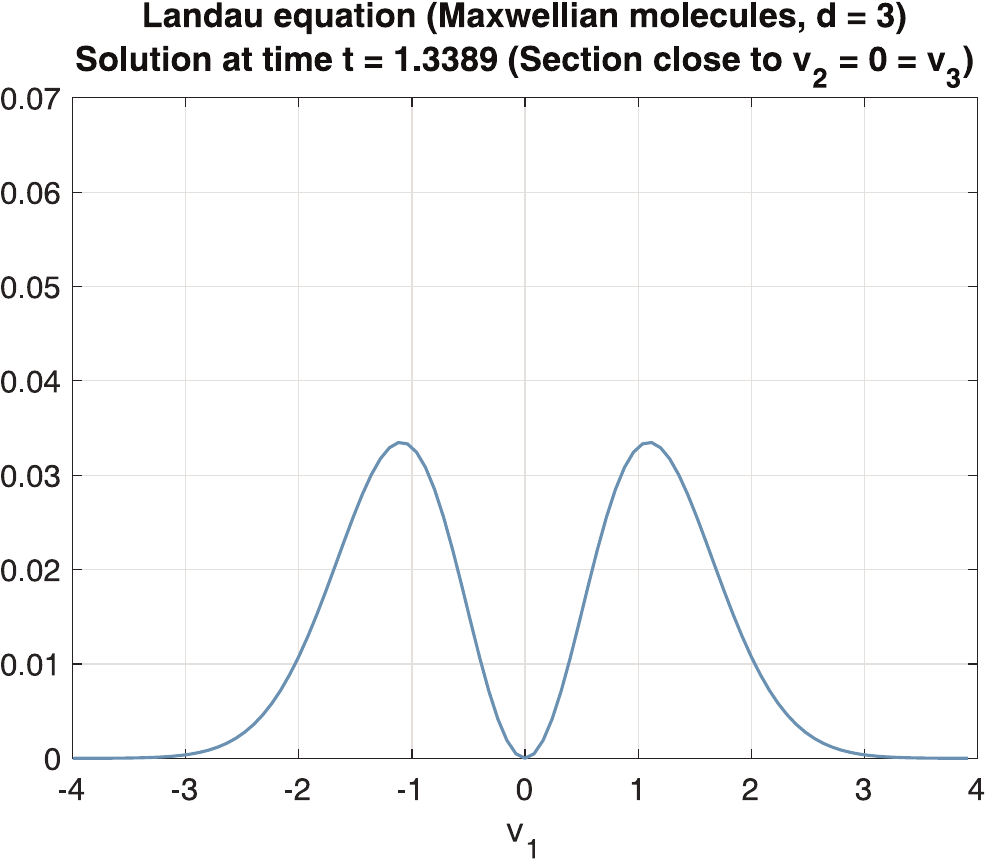}
\quad
\includegraphics[width=6.6cm]{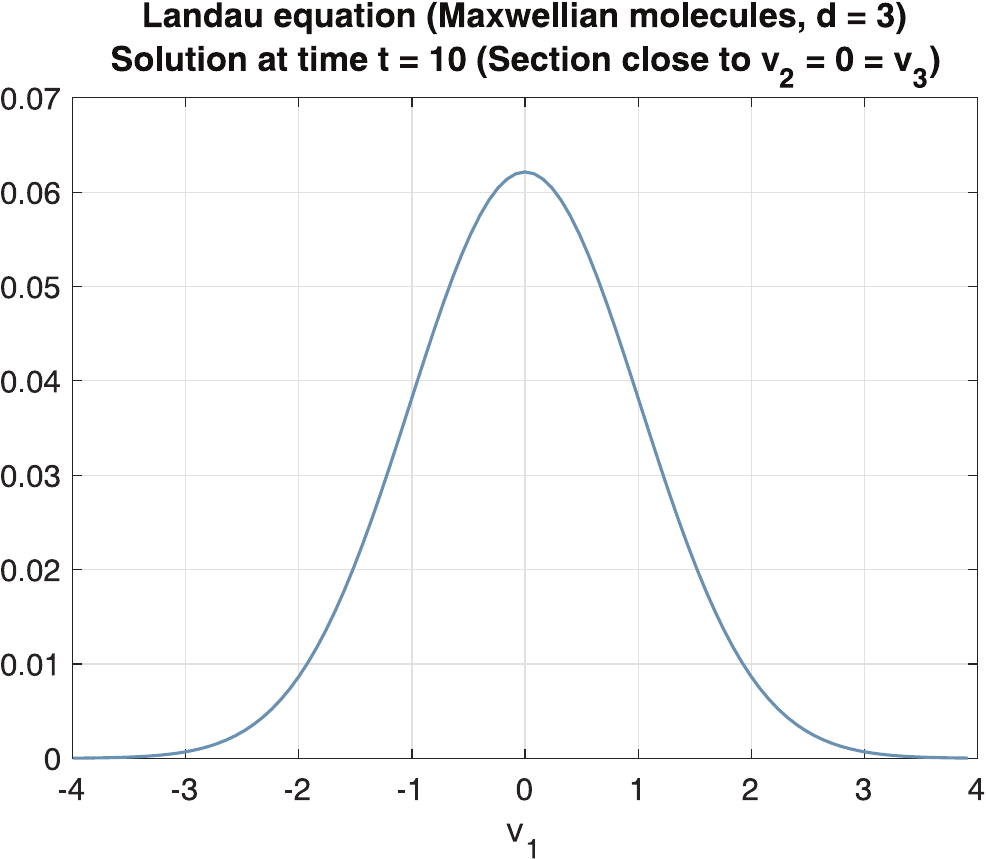}
\caption{Numerical illustration of the BKW solution to the Landau equation involving a constant kernel in three dimensions (Maxwellian molecules case), see~\eqref{eq:TestProblemA}.
The initial time $t_0 = 6 \, \ln(\frac{5}{4})$ is chosen in such a way that the non-negativity of the solution is ensured.}
\label{fig:SolutionAd3}
\end{center}
\end{figure}

%%%%%%%%%%%%%%%%%%%%%%%%%%%%%%%%%%%%%%%%%%%%%%%%%%%%%%%%%%%%%%%%%%%%%%%%%%%%%%%%%%%%%%%%%%%%%%%%%%%%%%%%%%%%%%%%%%%%%%%%%%%%%%%%%%%%%%%%%%%%%%%%%% 
\MyParagraph{Test problem~A (Constant integral kernel, Unbounded domain)}
In the special case of Maxwellian molecules in two and three dimensions, the integral kernels are defined by certain constants
\begin{subequations}
\label{eq:TestProblemA}
\begin{equation}
\begin{gathered}
d = 2\,, \quad \varphi = C = \tfrac{1}{16}\,, \\
d = 3\,, \quad \varphi = C = \tfrac{1}{24}\,, 
\end{gathered}
\end{equation}
and the underlying velocity domains coincide with the Euclidian spaces. 
In our numerical experiments, we make use of the fact that particular choices of the density functions permit to determine the associated integral and Landau operators by straightforward calculations 
\begin{equation}
\begin{gathered}
d = 2\,, \\
f(v) = \tfrac{1}{\pi} \, \ee^{- \abs{v}^2} \, \abs{v}^2\,, \\
\nQc(f, f)(v) = - \, \tfrac{1}{16 \, \pi} \, \ee^{- \abs{v}^2} \, \big(\abs{v}^2 - 2\big) \, v\,, \\
\nQ(f, f)(v) = \tfrac{1}{8 \, \pi} \, \ee^{- \abs{v}^2} \, \big(\abs{v}^4 - 4 \, \abs{v}^2 + 2\big)\,, \\
v \in \Omega = \RR^2\,,
\end{gathered}
\end{equation}
\begin{equation}
\begin{gathered}
d = 3\,, \\
f(v) = \tfrac{1}{2 \, \pi^{3/2}} \, \ee^{- \abs{v}^2} \, \big(2 \, \abs{v}^2 - 1\big)\,, \\
\nQc(f, f)(v) = - \, \tfrac{1}{12 \, \pi^{3/2}} \, \ee^{- \abs{v}^2} \, \big(\abs{v}^2 - \tfrac{5}{2}\big) \, v\,, \\
\nQ(f, f)(v) = \tfrac{1}{6 \, \pi^{3/2}} \, \ee^{- \abs{v}^2} \, \big(\abs{v}^4 - 5 \, \abs{v}^2 + \tfrac{15}{4}\big)\,, \\
v \in \Omega = \RR^3\,.
\end{gathered}
\end{equation}
Further numerical tests for the Landau equation rely on the knowledge of the BKW solutions 
\begin{equation}
\begin{gathered}
d = 2\,, \quad (\alpha_1, \alpha_2, \alpha_3) = \big(2, 1, \tfrac{1}{8}\big)\,, \\
d = 3\,, \quad (\alpha_1, \alpha_2, \alpha_3) = \big(\tfrac{5}{2}, \tfrac{3}{2}, \tfrac{1}{6}\big)\,, \\
K(t) = 1 - \tfrac{1}{2} \, \ee^{- \, \alpha_3 \, t}\,, \\
f(v, t) =  \tfrac{1}{(2 \, \pi \, K(t))^{d/2}} \, \ee^{- \frac{1}{2} \, \frac{1}{K(t)} \, \abs{v}^2} \,
\big(\alpha_1 - \alpha_2 \, \tfrac{1}{K(t)} + \tfrac{1}{2} \, \tfrac{1 - K(t)}{(K(t))^2} \, \abs{v}^2\big)\,, \\
v \in \Omega = \RR^d\,, \quad t \in [t_0, T]\,, 
\end{gathered}
\end{equation}
\end{subequations}
see~\cite[Ex.~1 and~2]{CarrilloEtAl2020}.
For three dimensions, the profile of the BKW solution is illustrated in Figure~\ref{fig:SolutionAd3}.

\MyParagraph{Test problem~B (Regular integral kernel, Bounded domain)}
First artificial test problems in two and three dimensions involve regular integral kernels and bounded velocity domains 
\begin{subequations}
\label{eq:TestProblemB}
\begin{equation}
\varphi(v) = \cos(v_1) \cdots \cos(v_d)\,, \quad v \in \Omega = [- \pi, \pi]^d\,.
\end{equation}
The prescribed density functions
\begin{equation}
f(v) = \sin(v_1) \cdots \sin(v_d) \,, \quad v \in \Omega = [- \pi, \pi]^d\,,
\end{equation}
are chosen such that the integral operators result from straightforward calculations
\begin{equation}
\begin{gathered}
d = 2\,, \\
q_1(v) = \tfrac{\pi^2}{8} \cos(2 \, v_2) \, \big(2 \, v_1 \cos(2 \, v_1) - \sin(2 \, v_1) - 2 \, v_1\big)\,, \\ 
q_2(v) = - \, \tfrac{\pi^2}{4} \cos(2 \, v_1) \, \sin(v_2) \, \big(2 \, v_2 \, \sin(v_2) + \cos(v_2)\big)\,, \\ 
\nQc(f, f)(v) = \big(q_1(v), q_2(v)\big)^T\,, \\
v \in \Omega = [- \pi, \pi]^2\,,
\end{gathered}
\end{equation}
\begin{equation}
\begin{gathered}
d = 3\,, \\    
q_{11}(v) = \big(\cos(2 \, v_3) - \tfrac{1}{2}\big) \, \big(v_1 \cos(2 \, v_1) - \tfrac{1}{2} \, \sin(2 \, v_1) - v_1\big) \cos(2 \, v_2)\,, \\
q_{12}(v) = \big(\tfrac{1}{4} \, \sin(2 \, v_1) - \tfrac{1}{2} \, v_1 \cos(2 \, v_1) \big) \, \cos(2 \, v_3) + v_1 \big(\cos^2(v_3) - \tfrac{1}{2}\big)\,, \\
q_1(v) = - \, \tfrac{\pi^3}{4} \, \big(q_{11}(v) + q_{12}(v)\big)\,, \\
q_{21}(v) = \big(\cos(2 \, v_3) - \tfrac{1}{2}\big) \, \big(v_2 \cos(2 \, v_2) - \tfrac{1}{2} \, \sin(2 \, v_2) - v_2\big) \cos(2 \, v_1)\,, \\
q_{22}(v) = \big(\tfrac{1}{4} \, \sin(2 \, v_2) - \tfrac{1}{2} \, v_2 \cos(2 \, v_2) \big) \cos(2 \, v_3) + v_2 \big(\cos^2(v_3) - \tfrac{1}{2}\big)\,, \\
q_2(v) = - \, \tfrac{\pi^3}{4} \, \big(q_{21}(v) + q_{22}(v)\big)\,, \\
q_{31}(v) = \tfrac{1}{4} \, \big(2 \, v_3 \cos(2 \, v_3) - \sin(2 \, v_3) - 2 \, v_3\big)\,, \\
q_{32}(v) = 2 \cos(2 \, v_1) \cos(2 \, v_2) - \cos(2 \, v_1) - \cos(2 \, v_2)\,, \\
q_3(v) = - \, \tfrac{\pi^3}{4} \, q_{31}(v) \, q_{32}(v)\,, \\
\nQc(f, f)(v) = \big(q_1(v), q_2(v), q_3(v)\big)^T\,, \\
v \in \Omega = [- \pi, \pi]^3\,.
\end{gathered}
\end{equation}
The associated Landau operators are obtained by differentiation
\begin{equation}
\begin{gathered}
d = 2\,, \\
q_1(v) = \big(2 \, v_2 \sin(2 \, v_2) + 1\big) \cos(2 \, v_1)\,, \\
q_2(v) = \big(2 \, v_1 \sin(2 \, v_1) + 1\big) \cos(2 \, v_2)\,, \\
\nQ(f, f)(v) = - \, \tfrac{\pi^2}{4} \, \big(q_1(v) + q_2(v)\big)\,, \\
v \in \Omega = [- \pi, \pi]^2\,,
\end{gathered}
\end{equation}
\begin{equation}
\begin{gathered}
d = 3\,, \\    
q_{11}(v) = v_2 \sin(2 \, v_2) \cos(2 \, v_3) + v_3 \cos(2 v_2) \sin(2 \, v_3)\,, \\
q_{12}(v) = \tfrac{1}{2} \, \big(- v_2 \sin(2 \, v_2) - v_3 \sin(2 \, v_3) + \cos(2 \, v_2) + \cos(2 \, v_3) - 1\big)\,, \\
q_1(v) = \big(q_{11}(v) + q_{12}(v)\big) \cos(2 \, v_1)\,, \\
q_{21}(v) = \big(v_1 \sin(2 \, v_1) + \tfrac{1}{2}\big) \cos(2 \, v_3)\,, \\
q_{22}(v) = - \, \tfrac{1}{2} \, \big(v_1 \sin(2 \, v_1) + v_3 \sin(2 \, v_3) + 1\big)\,, \\
q_2(v) = \big(q_{21}(v) + q_{22}(v)\big) \, \cos(2 \, v_2)\,, \\ 
q_3(v) = - \, \tfrac{1}{2} \, \big(v_1 \sin(2 \, v_1) + v_2 \sin(2 \, v_2)\big) \cos(2 \, v_3) - \cos^2(v_3) + \tfrac{1}{2}\,, \\
\nQ(f, f)(v) = \tfrac{\pi^3}{2} \, \big(q_1(v) + q_2(v) + q_3(v)\big)\,, \\ 
v \in \Omega = [- \pi, \pi]^3\,.
\end{gathered}
\end{equation}
\end{subequations}

\MyParagraph{Test problem~C (Regular integral kernel, Unbounded domain)}
Further artificial test problems in two and three dimensions are defined by Gaussian-like integral kernels and density functions 
\begin{subequations}
\label{eq:TestProblemC}
\begin{equation}
\begin{gathered}
d = 2\,, \quad \varphi(v) = \ee^{- v_1^2 - 2 \, v_2^2}\,, \quad
f(v) = \ee^{- \frac{1}{2} v_1^2 - \frac{1}{4} v_2^2}\,, \\
d = 3\,, \quad 
\varphi(v) = \ee^{- v_1^2 - 2 \, v_2^2 - 3 \, v_3^2}\,, \quad 
f(v) = \ee^{- \frac{1}{2} v_1^2 - \frac{1}{4} v_2^2 - \frac{1}{8} v_3^2}\,, \\
v \in \Omega = \RR^d\,.
\end{gathered}
\end{equation}
The associated integral and Landau operators are given by 
\begin{equation}
\begin{gathered}
d = 2\,, \\
q_1(v) = \sqrt{6} \, \pi \, \ee^{- \frac{5}{6} v_1^2 - \frac{17}{36} v_2^2}\,, \\   
\nQc(f, f)(v) = q_1(v) \, \Big(- \, \tfrac{1}{2187} \, v_1 \, \big(v_2^2 + 18\big), \tfrac{1}{729} \, v_2 \, \big(v_1^2 + 3\big)\Big)^T\,, \\
\nQ(f, f)(v) = - \, \tfrac{1}{13122} \, q_1(v) \, \big(7 \, v_1^2 \, v_2^2 - 198 \, v_1^2 + 57 \, v_2^2 + 54)\,, \\
v \in \Omega = \RR^2\,,
\end{gathered}
\end{equation}
\begin{equation}
\begin{gathered}
d = 3\,, \\
q_1(v) = \sqrt{3} \, \pi^{3/2} \, \ee^{- \frac{5}{6} v_1^2 - \frac{17}{36} v_2^2 - \frac{49}{200} v_3^2}\,, \\
q_{21}(v) = - \, \tfrac{1}{6834375} \, v_1 \, \big(1250 \, v_2^2 + 243 \, v_3^2 + 46800\big)\,, \\
q_{22}(v) = \tfrac{1}{2278125} \, v_2 \, \big(1250 \, v_1^2 - 9 \, v_3^2 + 2850\big)\,, \\
q_{23}(v) = \tfrac{1}{91125} \, v_3 \, \big(27 \, v_1^2 + v_2^2 + 99\big)\,, \\
\nQc(f, f)(v) = 2 \, q_1(v) \, \big(q_{21}(v), q_{22}(v), q_{23}(v)\big)^T\,, \\
q_{31}(v) = 17500 \, v_1^2 \, v_2^2 + 7047 \, v_1^2 \, v_3^2 + 135 \, v_2^2 \, v_3^2\,, \\ 
q_{31}(v) = - \, 1005300 \, v_1^2 + 111000 \, v_2^2 + 46899 \, v_3^2 + 369900\,, \\ 
Q(f, f)(v) = - \, \tfrac{1}{41006250} \, q_1(v) \, \big(q_{31}(v) + q_{32}(v)\big) \,, \\
v \in \Omega = \RR^3\,.
\end{gathered}
\end{equation}
\end{subequations}

\begin{figure}
\begin{center}
\includegraphics[width=6.6cm]{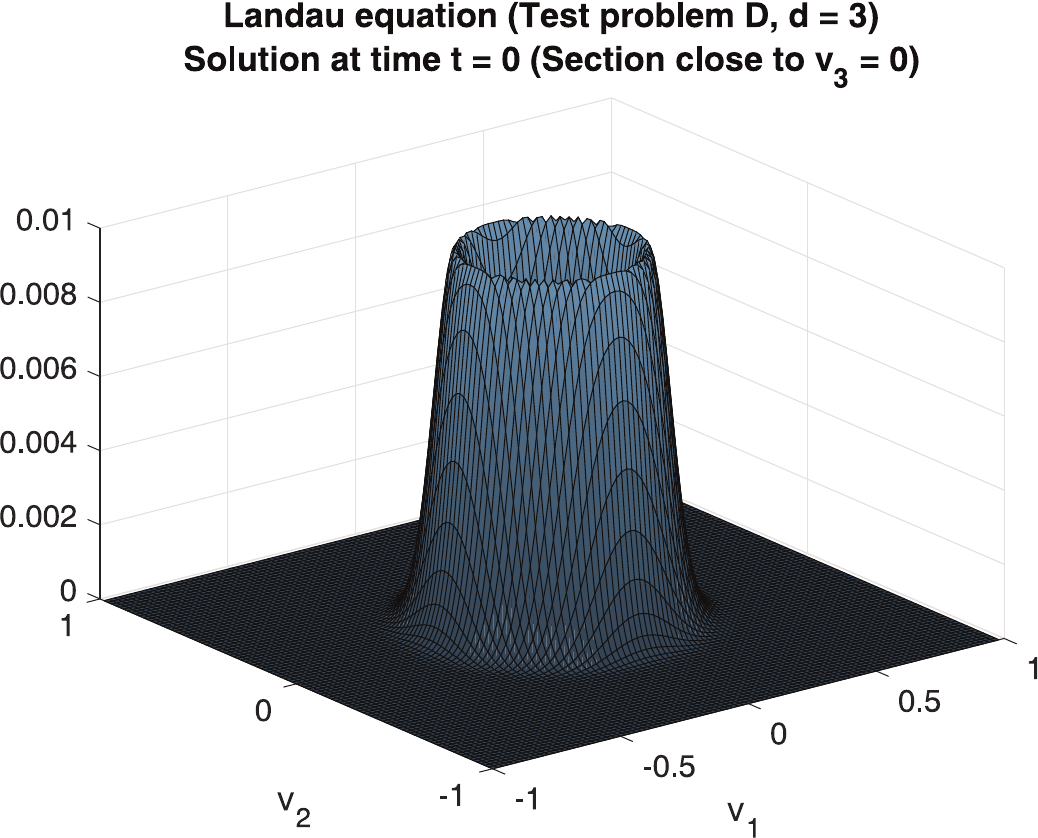}
\quad
\includegraphics[width=6.6cm]{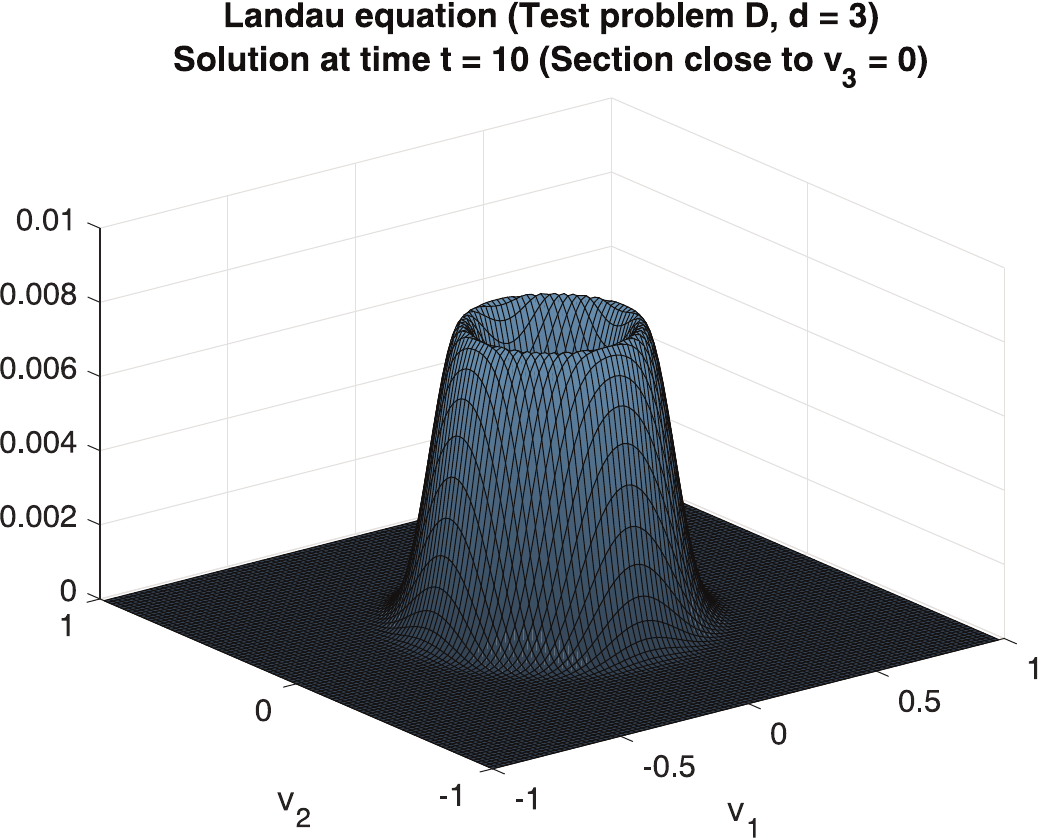}
\\[2mm]
\includegraphics[width=6.6cm]{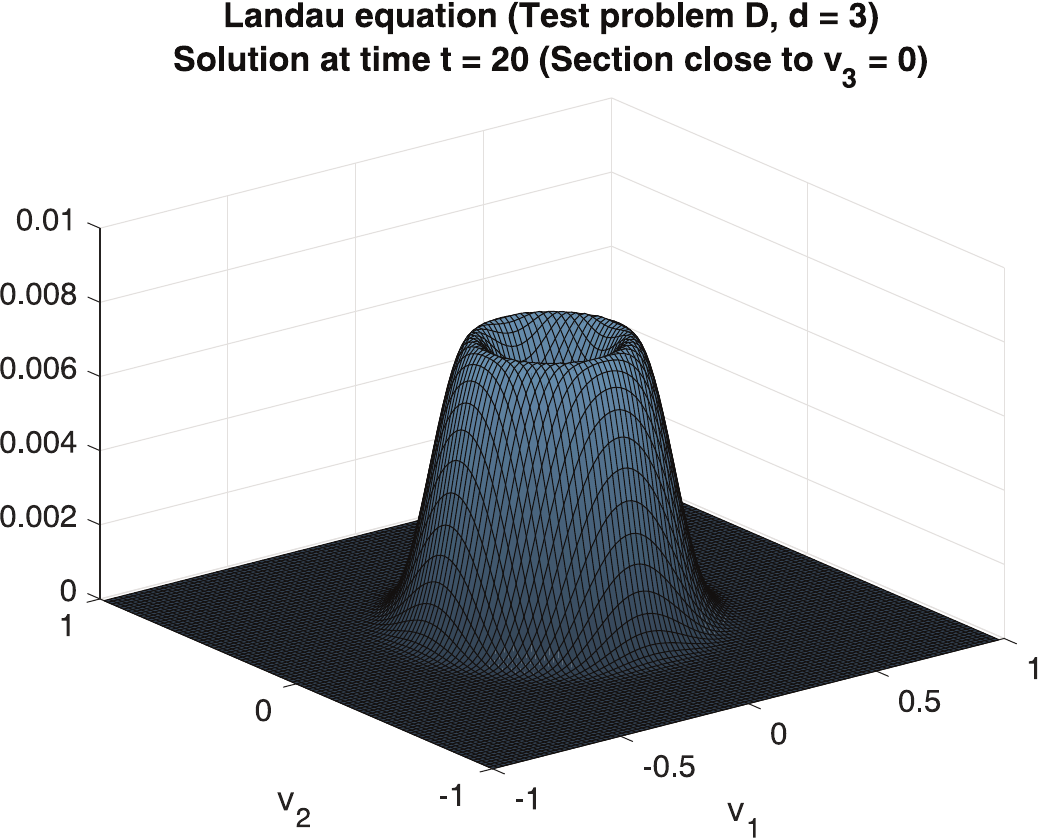}
\quad
\includegraphics[width=6.6cm]{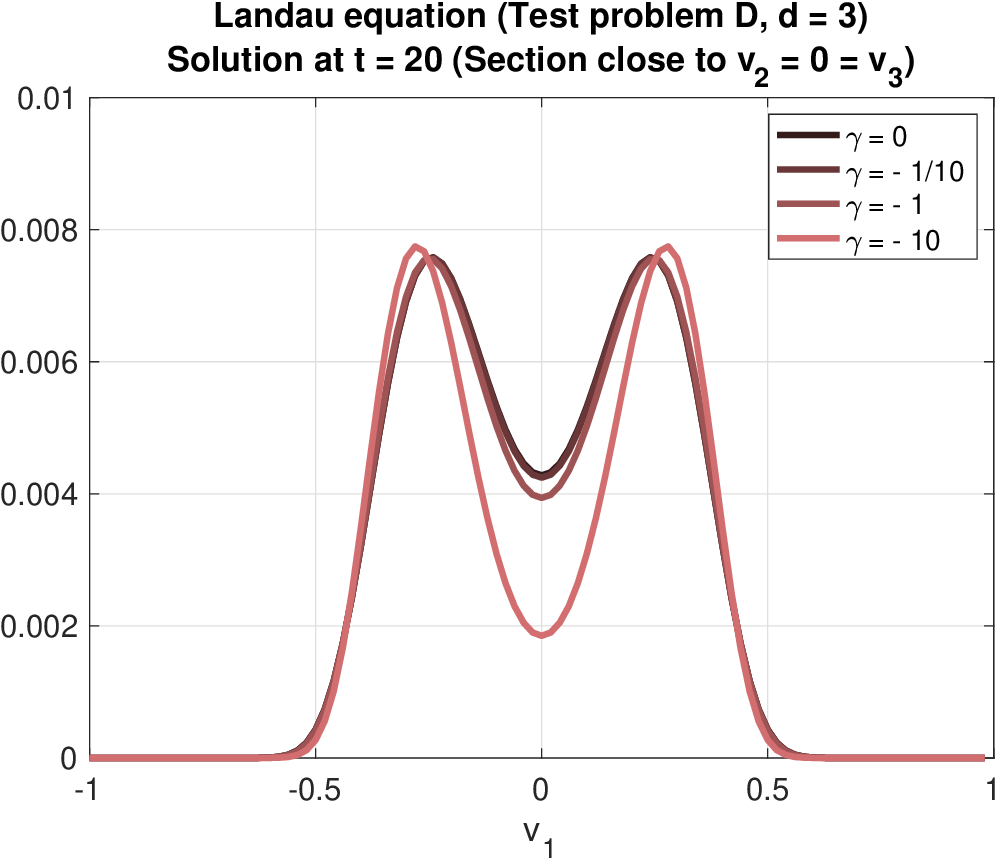}
\caption{Numerical illustration of the solution to the Landau equation involving a singular kernel in three dimensions. 
Comparison of the solution profiles for different exponents $\gamma \in \{0, - \, \tfrac{1}{10}, - \, 1, - \, 10\}$, see~\eqref{eq:TestProblemDCoulombGeneral}.}
\label{fig:SolutionDd3}
\end{center}
\end{figure}

\MyParagraph{Test problem~D (Singular integral kernel, Unbounded domain)}
As final test problem, we study the numerical evaluation of Landau operators that involve singular integral kernels.
Accordingly to our general setting and in the lines of~\cite[Ex.~3 and~4]{CarrilloEtAl2020} illustrating a two-dimensional anisotropic solution and the three-dimensional Rosenbluth problem with Coulomb potential, we in particular set 
\begin{equation}
\label{eq:TestProblemDCoulombGeneral}
\begin{gathered}
\varphi(z) = C \, \abs{z}^{\beta} \, \ee^{\gamma \, \abs{z}}\,, \quad \beta = - \, 3\,, \quad \gamma \in \{0, - \, \tfrac{1}{10}, - \, 1, - \, 10\}\,, \quad z \in \RR^d\,, \\
d = 2\,, \\
C = \tfrac{1}{16}\,, \quad f(v) = \tfrac{1}{4 \, \pi} \, \big(\ee^{- \frac{1}{2} \, ((v_1 + 2)^2 + (v_2 - 1)^2)} + \ee^{- \frac{1}{2} \, (v_1^2 + (v_2 + 1)^2)}\big)\,, \\
d = 3\,, \\
C = \tfrac{1}{4 \, \pi}\,, \quad f(v) = \tfrac{1}{c_1^2} \, \ee^{- c_1/c_2^2 \, (\abs{v} - c_2)^2}\,, \quad c_1 = 10\,, \quad c_2 = \tfrac{3}{10}\,, \\
v \in \Omega = \RR^d\,,
\end{gathered}
\end{equation}
see~\eqref{eq:KernelCoulombGeneral}.
For the relevant three-dimensional case, the profile of the solution to the associated Landau equation is shown in Figure~\ref{fig:SolutionDd3}.
%%%%%%%%%%%%%%%%%%%%%%%%%%%%%%%%%%%%%%%%%%%%%%%%%%%%%%%%%%%%%%%%%%%%%%%%%%%%%%%%%%%%%%%%%%%%%%%%%%%%%%%%%%%%%%%%%%%%%%%%%%%%%%%%%%%%%%%%%%%%%%%%%%
\subsection{Numerical results}
\label{subsec:NumericalResults}
%%%%%%%%%%%%%%%%%%%%%%%%%%%%%%%%%%%%%%%%%%%%%%%%%%%%%%%%%%%%%%%%%%%%%%%%%%%%%%%%%%%%%%%%%%%%%%%%%%%%%%%%%%%%%%%%%%%%%%%%%%%%%%%%%%%%%%%%%%%%%%%%%% 
\MyParagraph{Reliability and efficiency}
In our numerical tests for the different approaches outlined in Table~\ref{tab:Table1}, we first resort to the test problems A--C with known solutions, see~\eqref{eq:TestProblemA} to~\eqref{eq:TestProblemC}.
We in particular monitor the achieved accuracy when evaluating the Landau operator and integrating the associated Landau equation in time.
In addition, we measure the reliability of our approaches in a long-term integration through the conservation of mass, momentum, and energy as well as the decay of entropy, see~\eqref{eq:MassMomentumEnergyEntropy}.
Afterwards, we extend our studies to the most relevant cases of singular integral kernels~\eqref{eq:TestProblemDCoulombGeneral}. 
It is notable that our approach is flexible with regard to modifications of Coulomb interactions involving a coupling of all velocity directions.
In these general cases, quadrature approximations on relatively small subdomains are needed.
In special cases, where a decomposition in certain velocity directions is possible, the corresponding component of the code is adapted and the computational effort reduces accordingly. 
We recall that our strategy for the regularisation of the singular integral kernel is illustrated in Figure~\ref{fig:Kernel}.
With the conclusion from Table~\ref{tab:Table2} in mind, we use quadrature approximations based on few grid points. 
We point out that the observed computation times in \textsc{Matlab} are erratic and thus serve as rough indicators for the overall cost.
We expect improvements for implementations based on efficient software packages.

\begin{figure}
\begin{center}
\includegraphics[width=6.6cm]{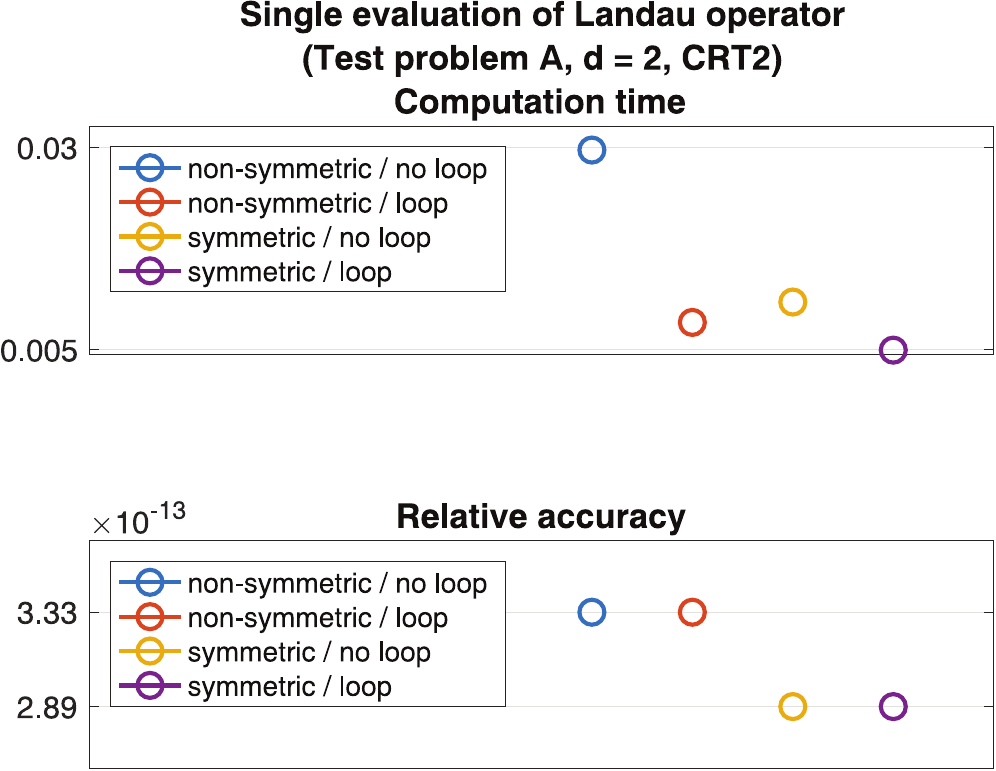}
\quad
\includegraphics[width=6.6cm]{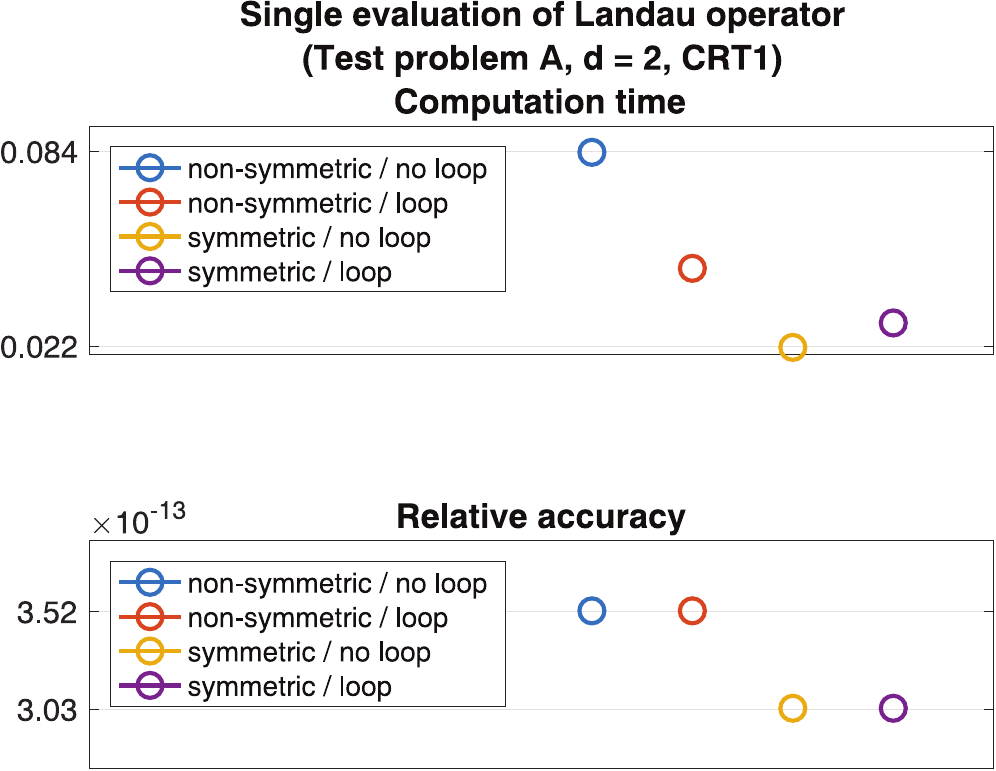} \\[8mm]
\includegraphics[width=6.6cm]{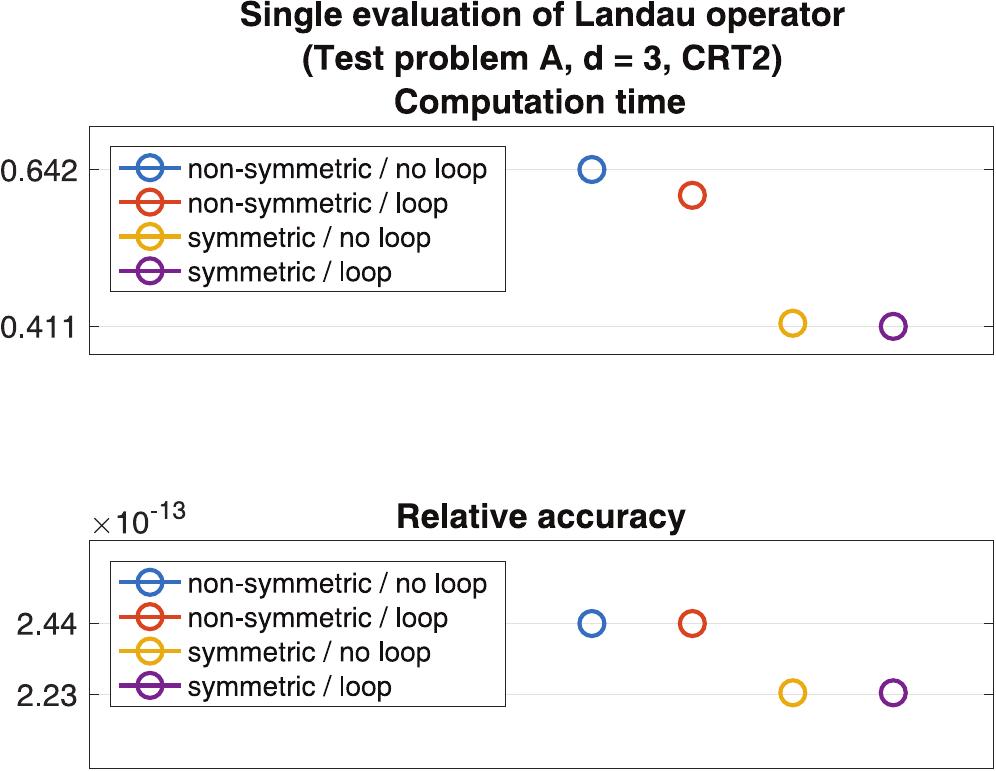}
\quad
\includegraphics[width=6.6cm]{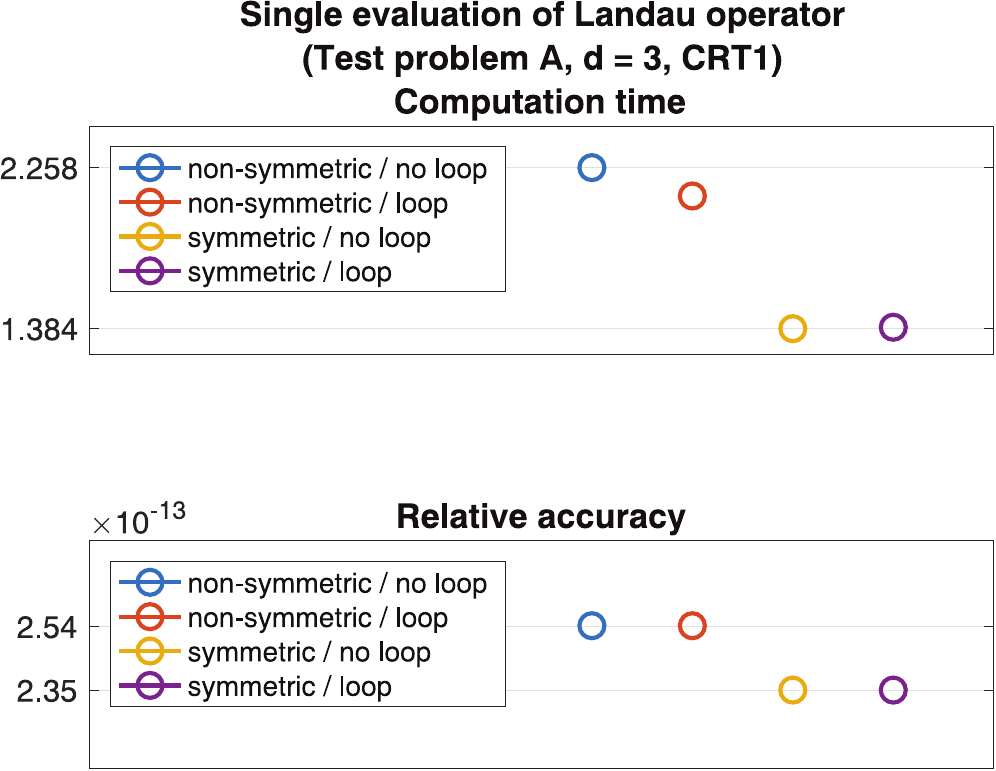} 
\caption{Test problem~A (Maxwellian molecules case) in two ($d = 2$) and three ($d = 3$) dimensions.
The evaluation of the Landau operator is based on the approaches CRT2 (left) and CRT1 (right) described in Section~\ref{sec:SectionCRT1-2}. 
For different implementations (velocity domains defined by non-symmetric versus symmetric intervals, computation of fundamental integrals without or with for-loops), consistent results concerning accuracy and computation time are observed.}
\label{fig:Comparison1}
\end{center}
\end{figure}

\MyParagraph{First validation}
In sight of the complexity of the problem, it is reasonable to validate our theoretical considerations in two and three dimensions step-by-step.
First numerical tests rely on the knowledge of the Landau operators in the Maxwellian molecules case (Test problem~A) for special choices of the density functions, see~\eqref{eq:TestProblemA}. 
We study the approaches CRT2 and CRT1, described in Section~\ref{sec:SectionCRT1-2}, which are suitable for regular kernels. 
On the one hand, we vary the truncated domain and the Fourier spectral discretisation in each direction 
\begin{equation*}
\begin{gathered}
d = 2\,, \\  
\Omb = [- \, 9, 10] \times [- \, 10, 11]\,, \quad M = (100, 110)\,, \\
d = 3\,, \\  
\Omb = [- \, 9, 10] \times [- \, 10, 11] \times [- \, 11,12]\,, \quad M = (100, 110, 120)\,, 
\end{gathered}
\end{equation*}
and on the other hand, we set 
\begin{equation}
\label{eq:DomainSymmetric}  
\Omb = [-10, 10]^d\,, \quad M_i = 100\,, \quad i \in \{1, \dots, d\}\,, 
\end{equation}
see also~\eqref{eq:Omb} and~\eqref{eq:nM}.
Besides, we contrast implementations using for-loops over index sets such as~\eqref{eq:nM123} or not, respectively.
Due to the fact that consistent results regarding accuracy and computation time are observed in all cases, see Figure~\ref{fig:Comparison1}, we proceed with thorough numerical comparisons of the different approaches.

\begin{figure}
\begin{center}
\includegraphics[width=4.2cm]{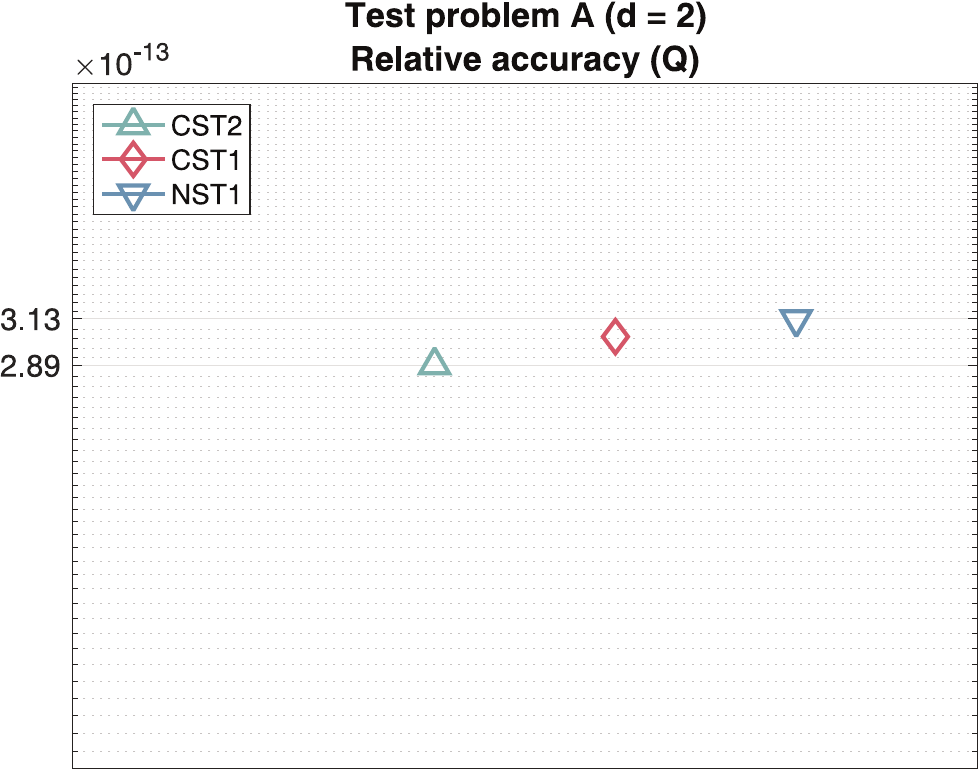} \quad
\includegraphics[width=4.2cm]{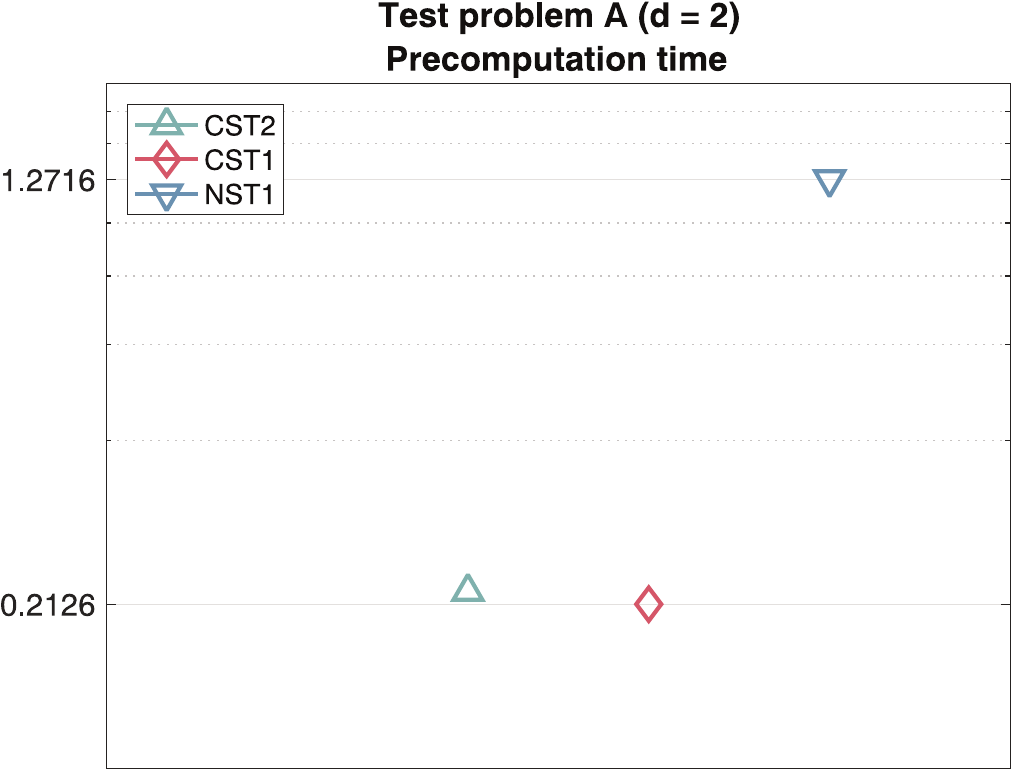} \quad
\includegraphics[width=4.2cm]{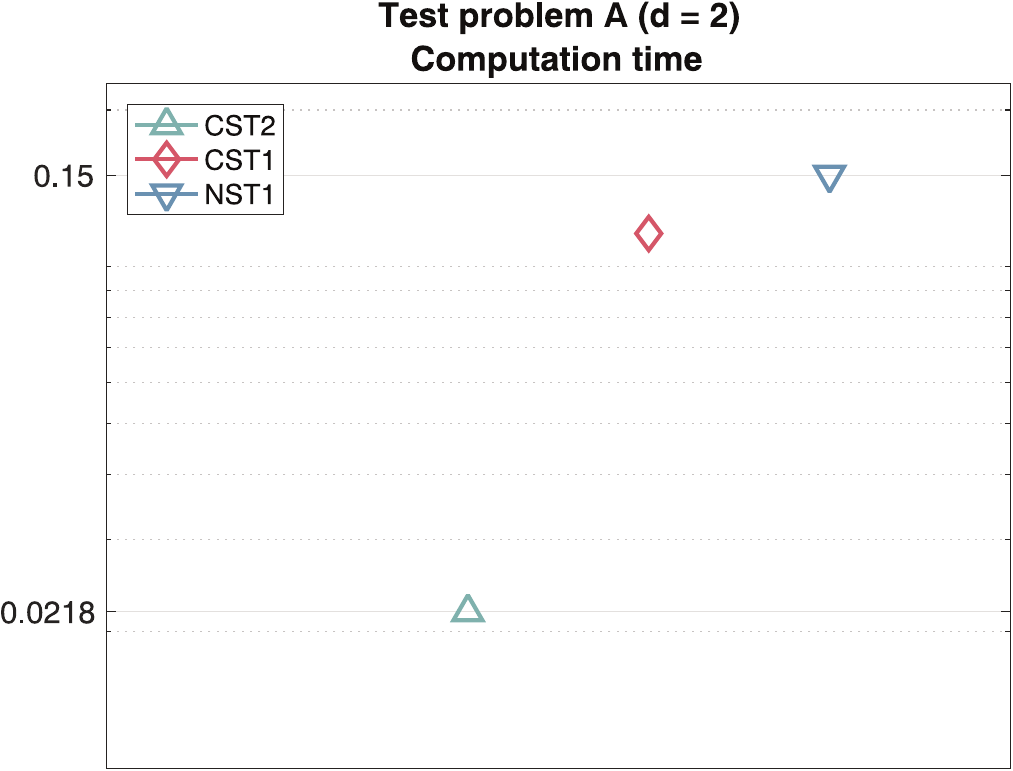} \\[2mm]
\includegraphics[width=4.2cm]{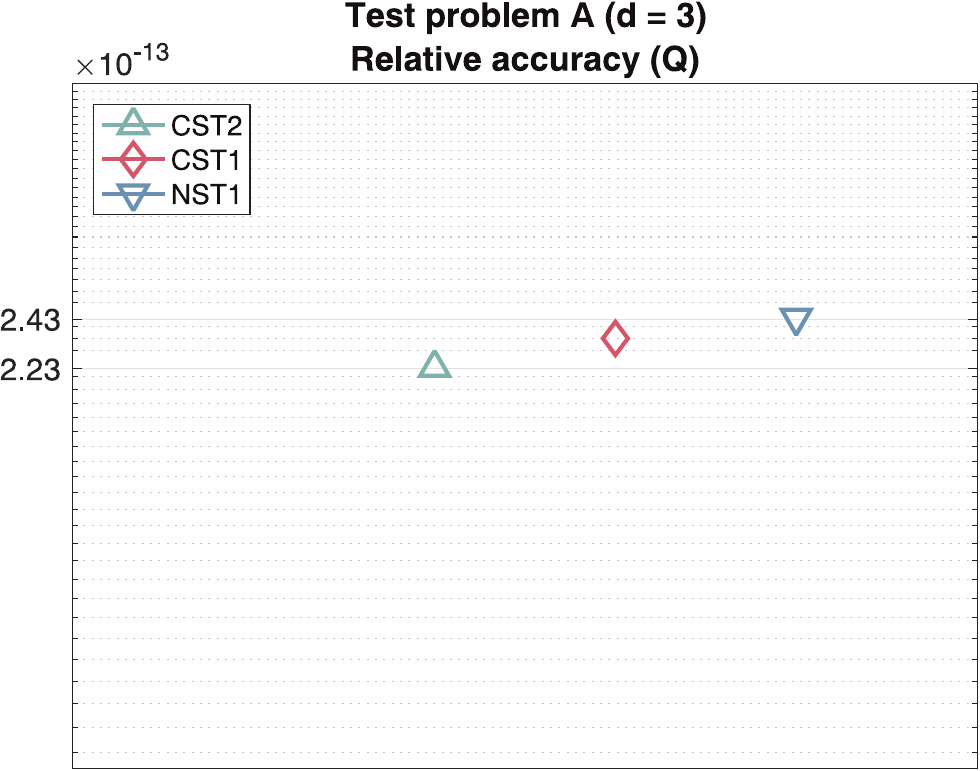} \quad
\includegraphics[width=4.2cm]{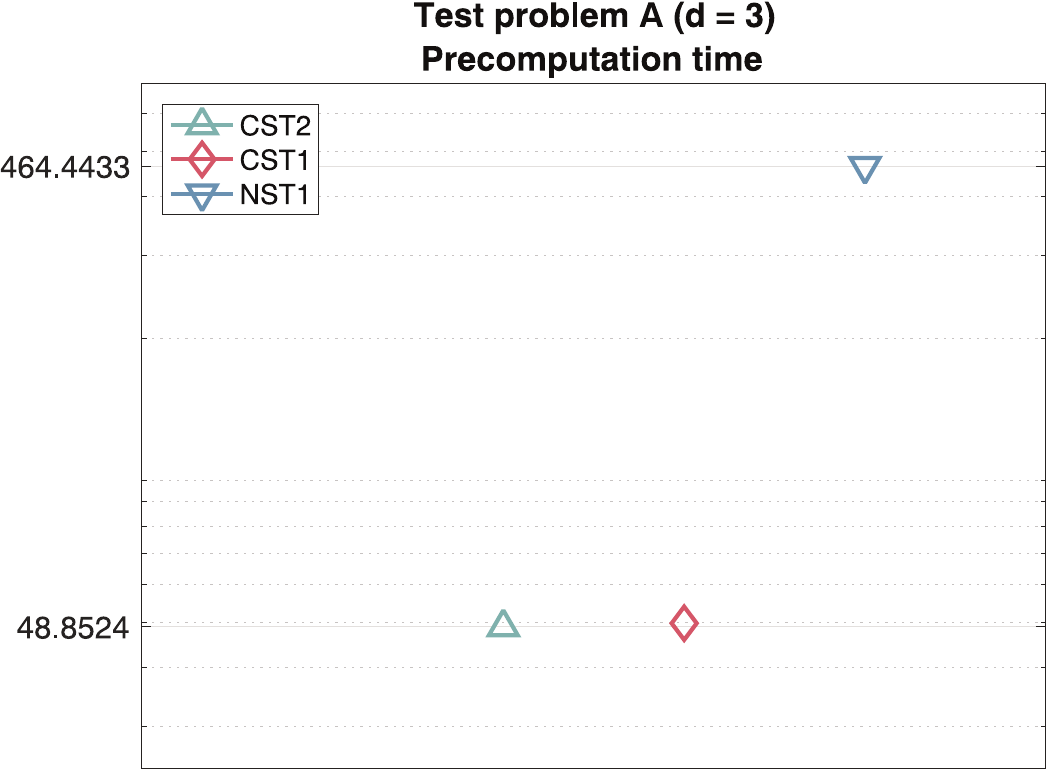} \quad 
\includegraphics[width=4.2cm]{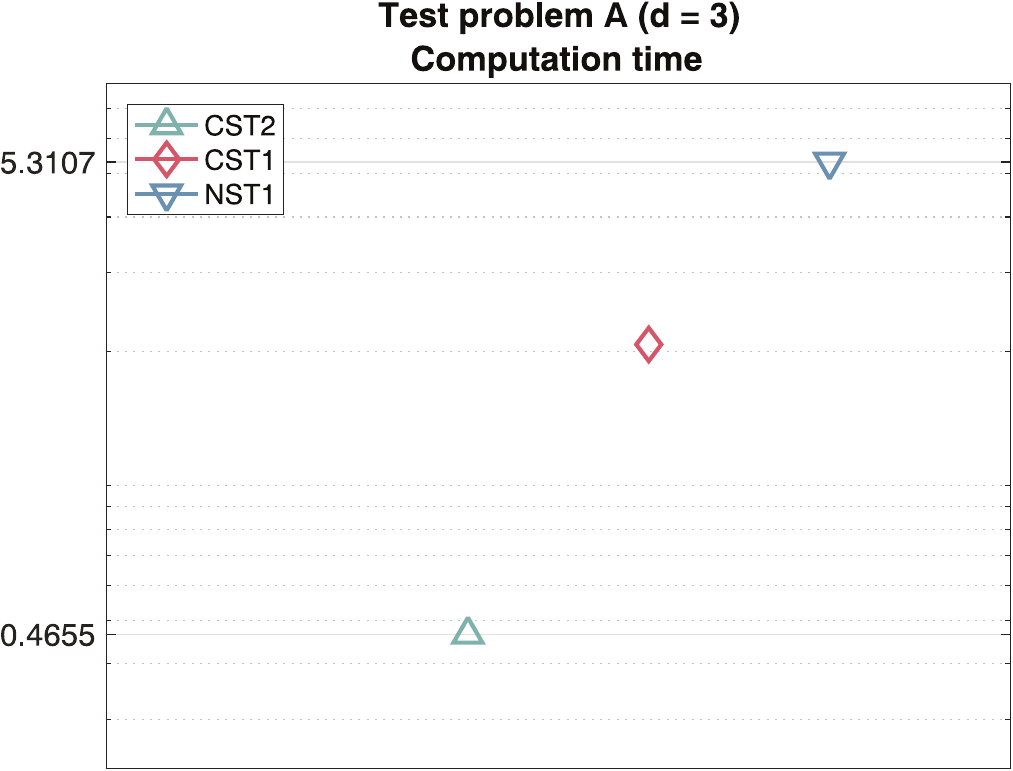} \\[2mm]
\includegraphics[width=4.2cm]{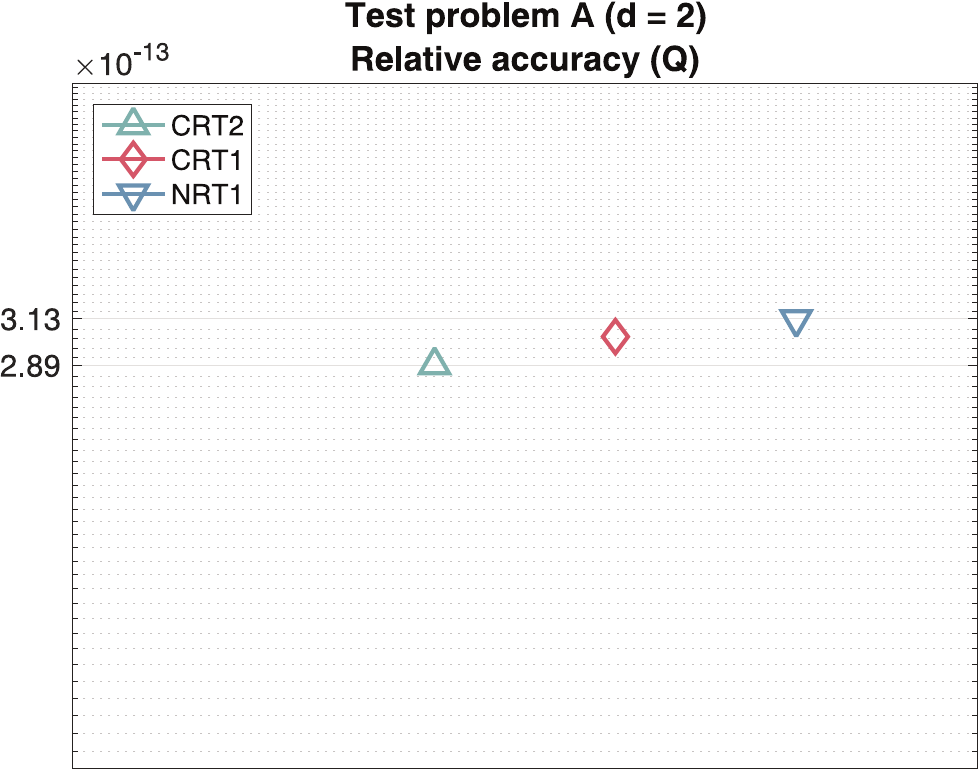} \quad
\includegraphics[width=4.2cm]{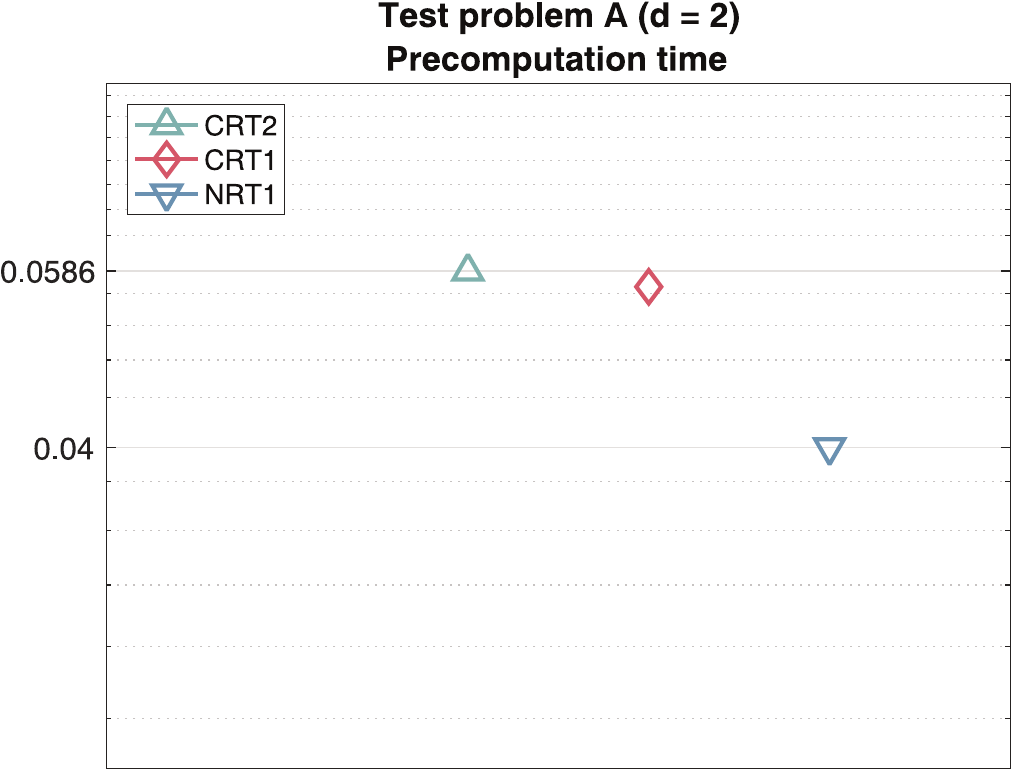} \quad
\includegraphics[width=4.2cm]{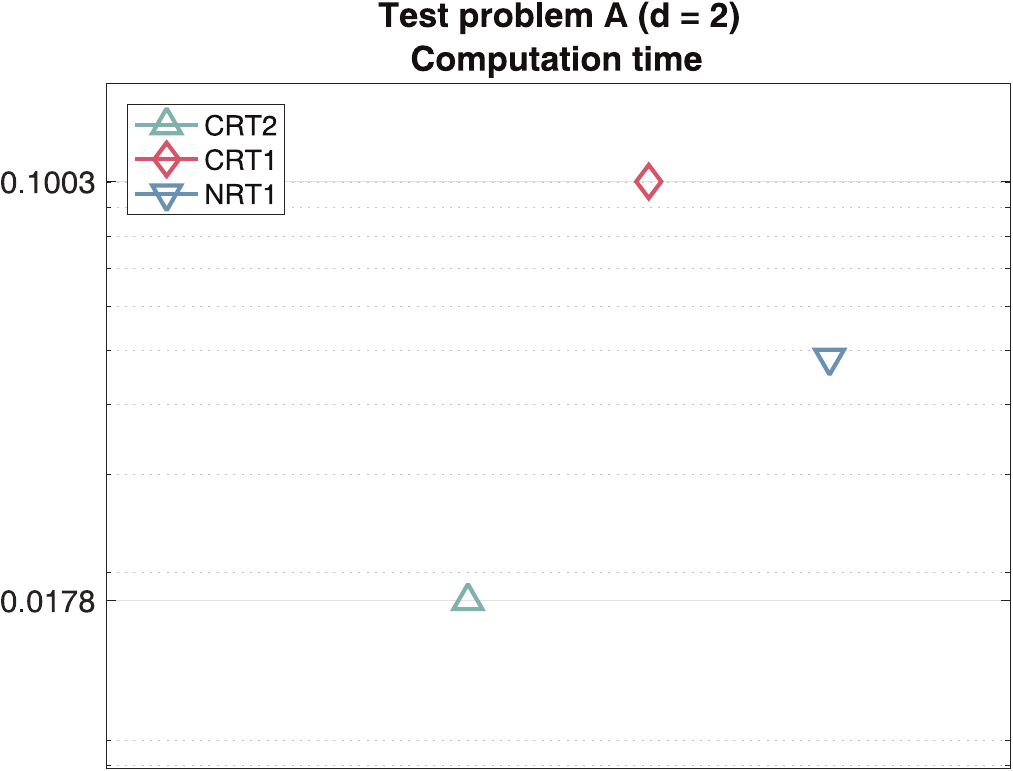} \\[2mm]
\includegraphics[width=4.2cm]{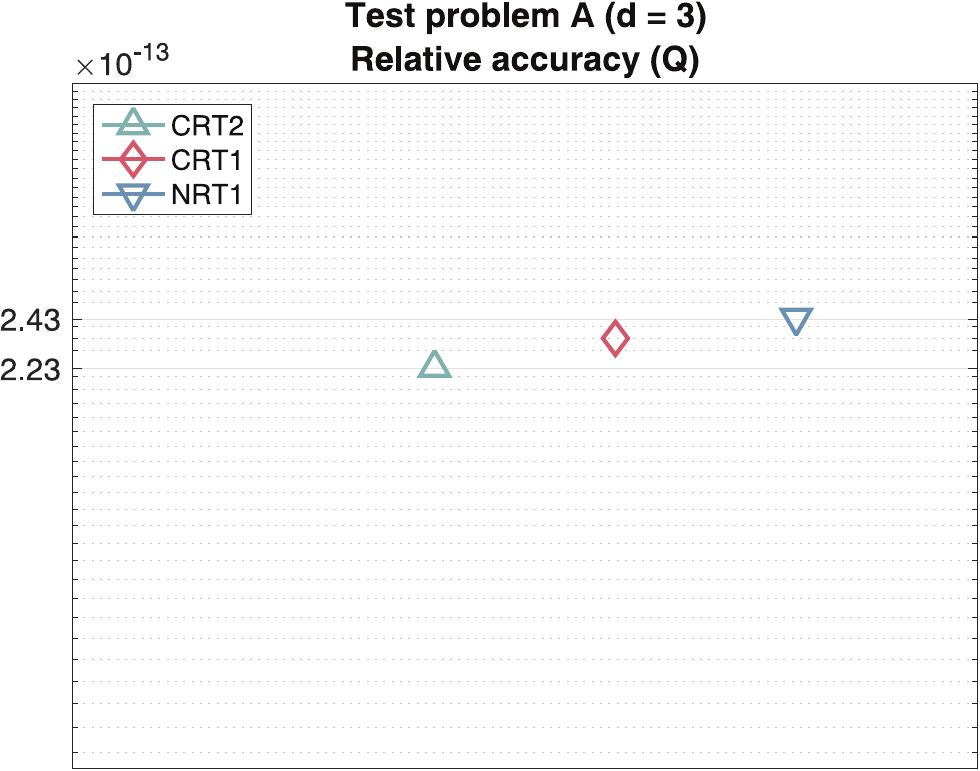} \quad
\includegraphics[width=4.2cm]{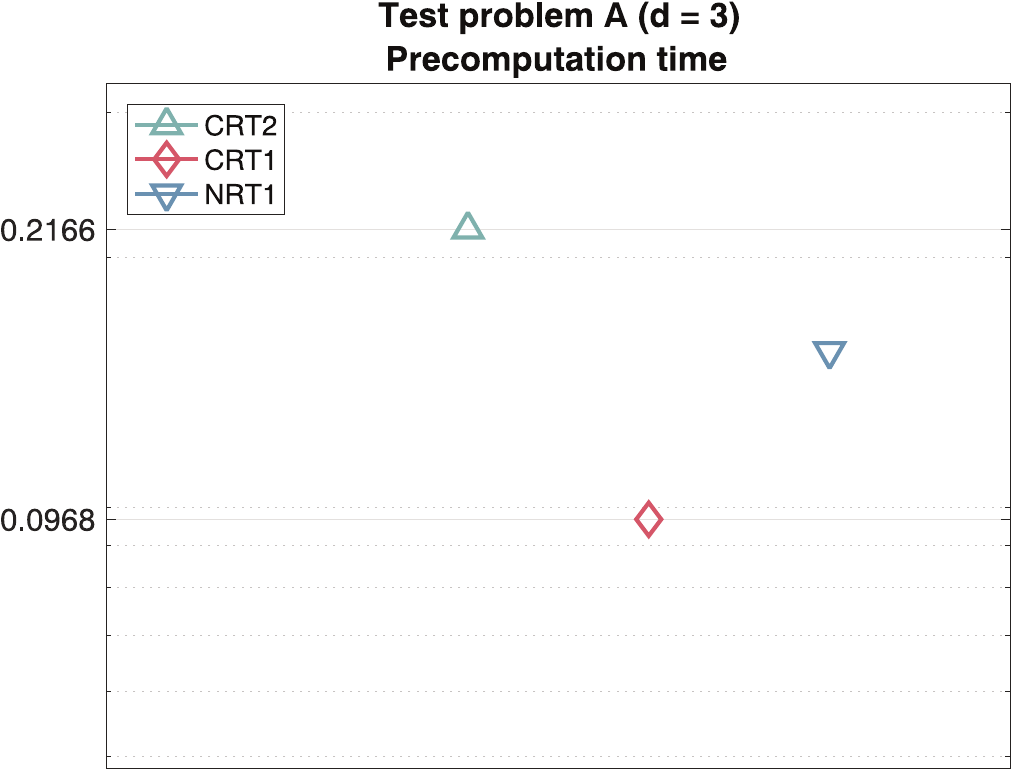} \quad
\includegraphics[width=4.2cm]{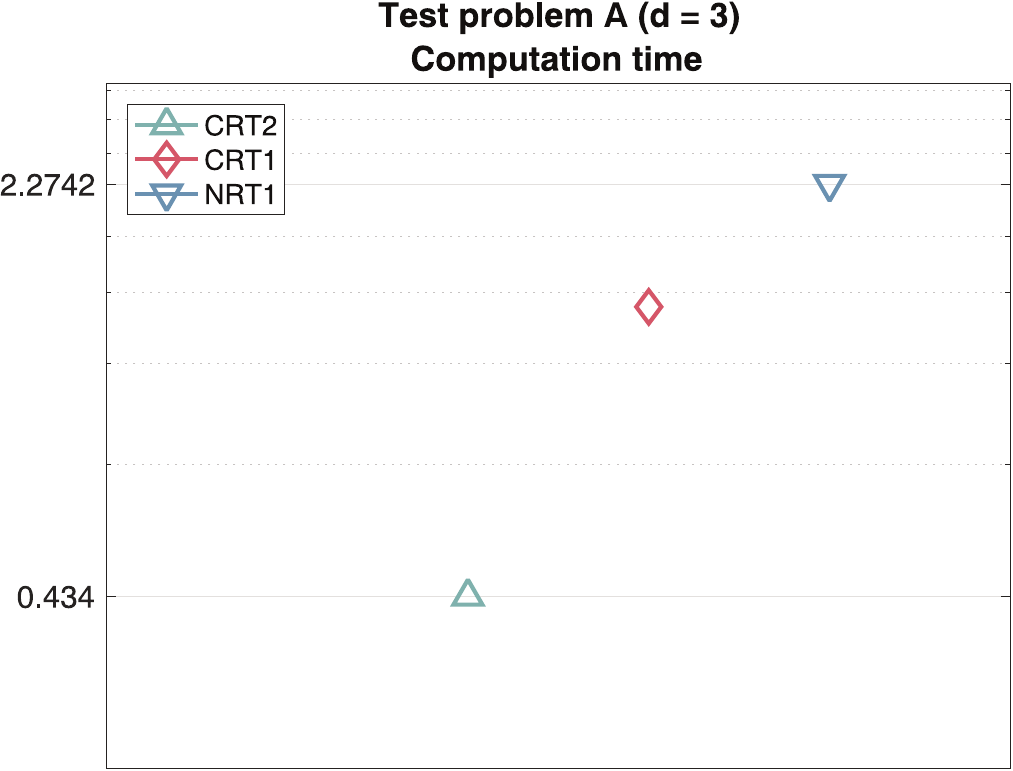} 
\caption{Test problem~A (Maxwellian molecules case) in two and three dimensions. Numerical comparisons of relative accuracies as well as precomputation and computation times.
First and second rows: General approaches for the evaluation of the associated Landau operator based on the conservative form (CST2, CST1) and the non-conservative form (NST1).
Third and fourth rows: Simplifications to regular kernels (CRT2, CRT1, NRT1).
In all cases, highly accurate results are obtained.}
\label{fig:ComparisonA}
\end{center}
\end{figure}

\begin{figure}
\begin{center}
\includegraphics[width=6.6cm]{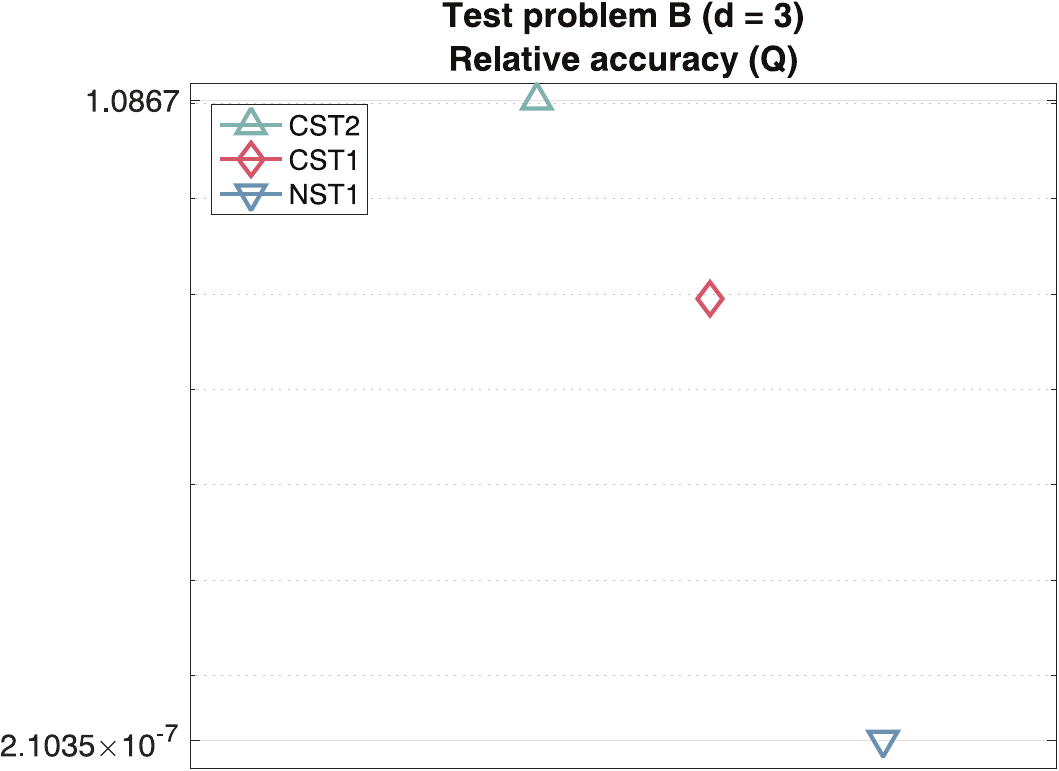} 
\quad
\includegraphics[width=6.6cm]{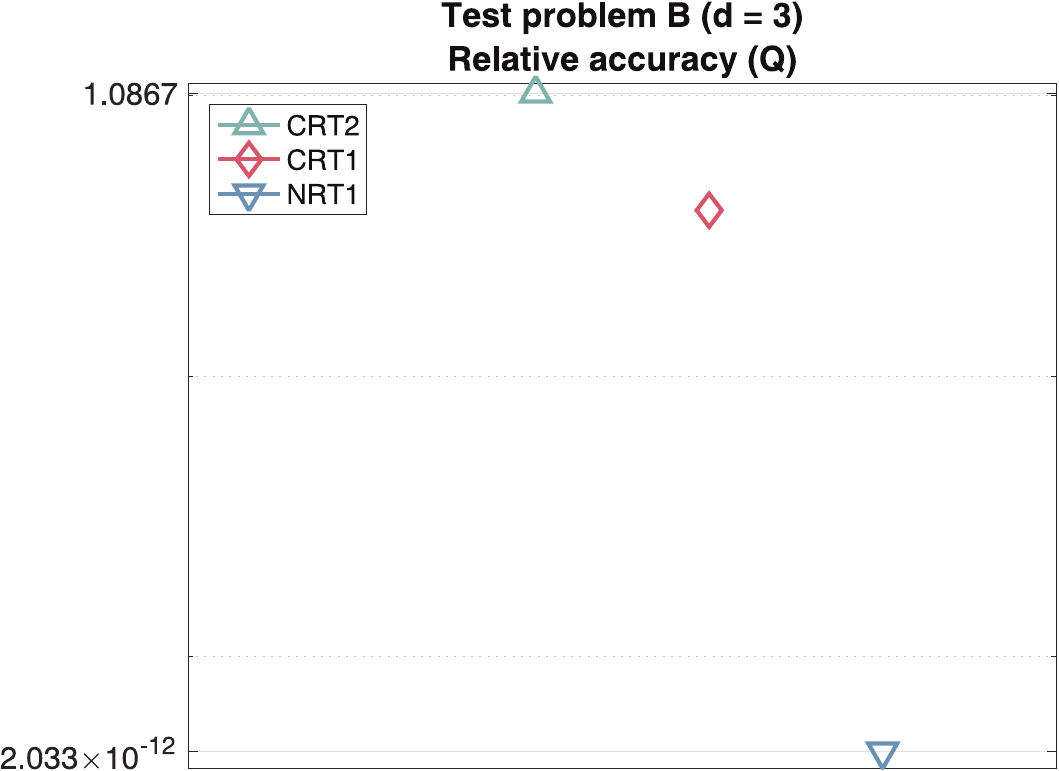}
\caption{Test problem~B (regular integral kernel, bounded domain, known solution). The approaches based on the non-conservative formulation yield a satisfactory result (NST1, NRT1), whereas the approaches based on the conservative formulation and numerical differentiation are not suitable (CST2, CST1, CRT2, CRT1).}
\label{fig:ComparisonB}
\end{center}
\end{figure}

\begin{figure}
\begin{center}
\includegraphics[width=4.2cm]{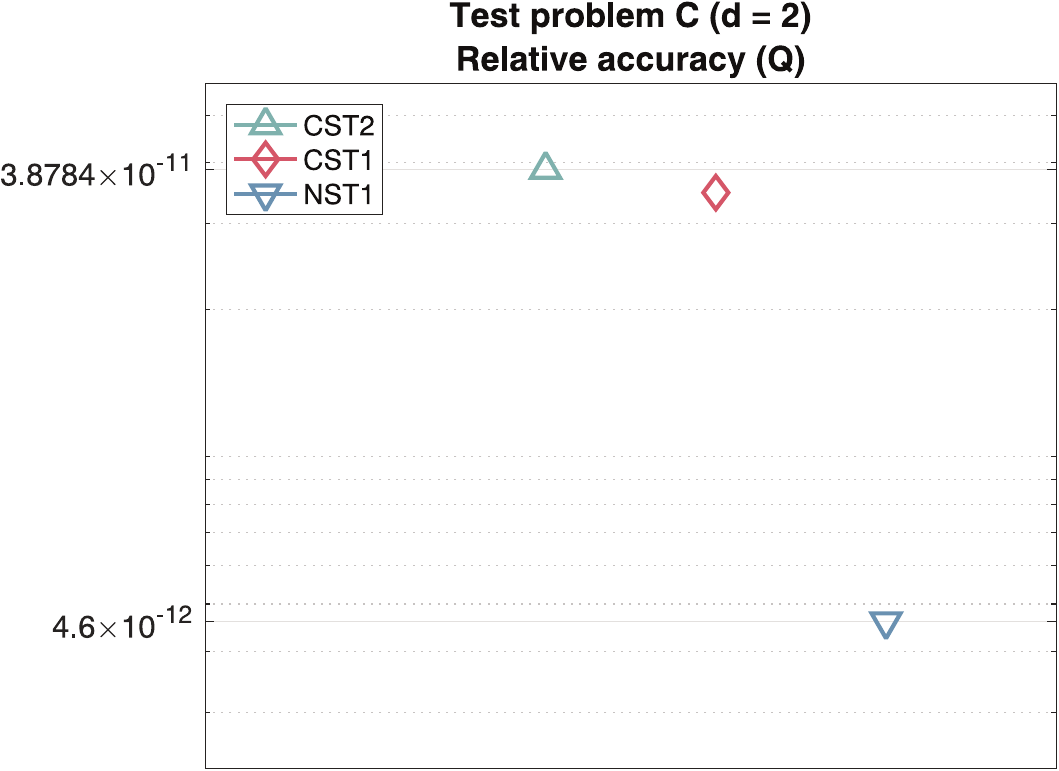} \quad
\includegraphics[width=4.2cm]{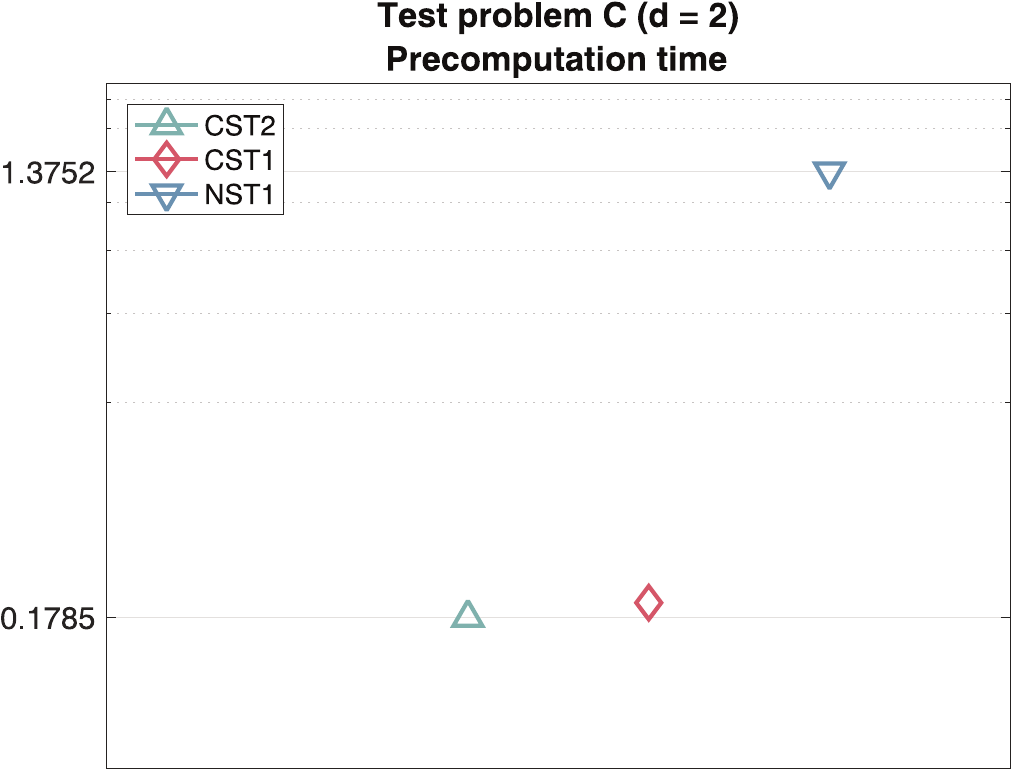} \quad
\includegraphics[width=4.2cm]{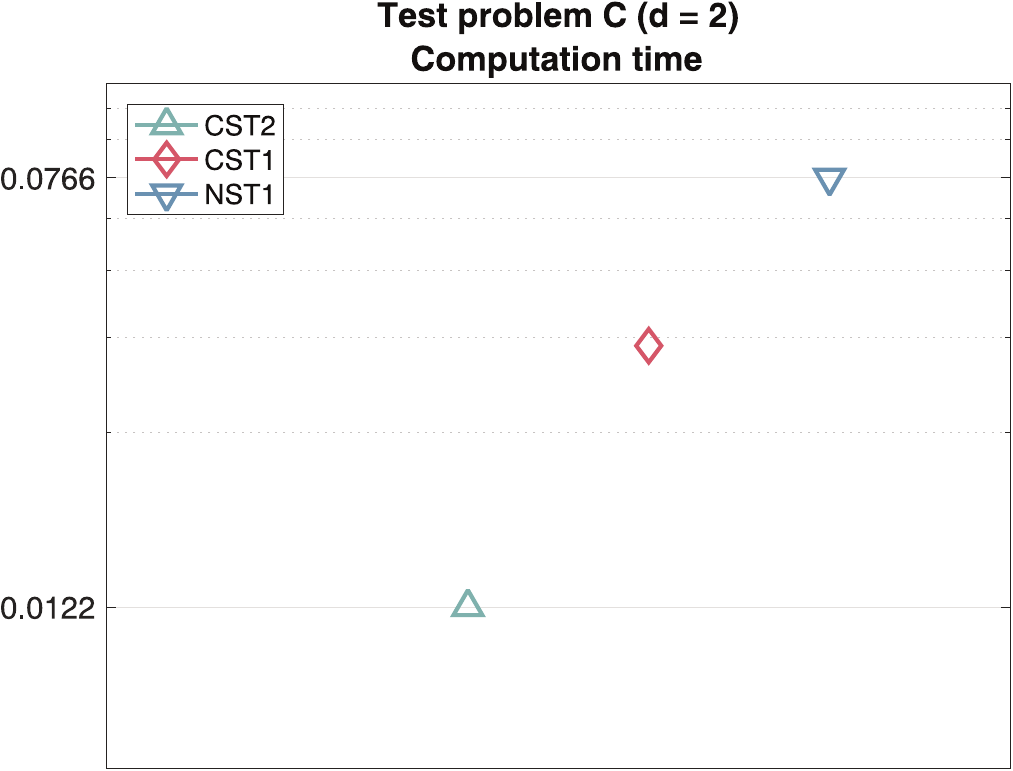} \\[2mm]
\includegraphics[width=4.2cm]{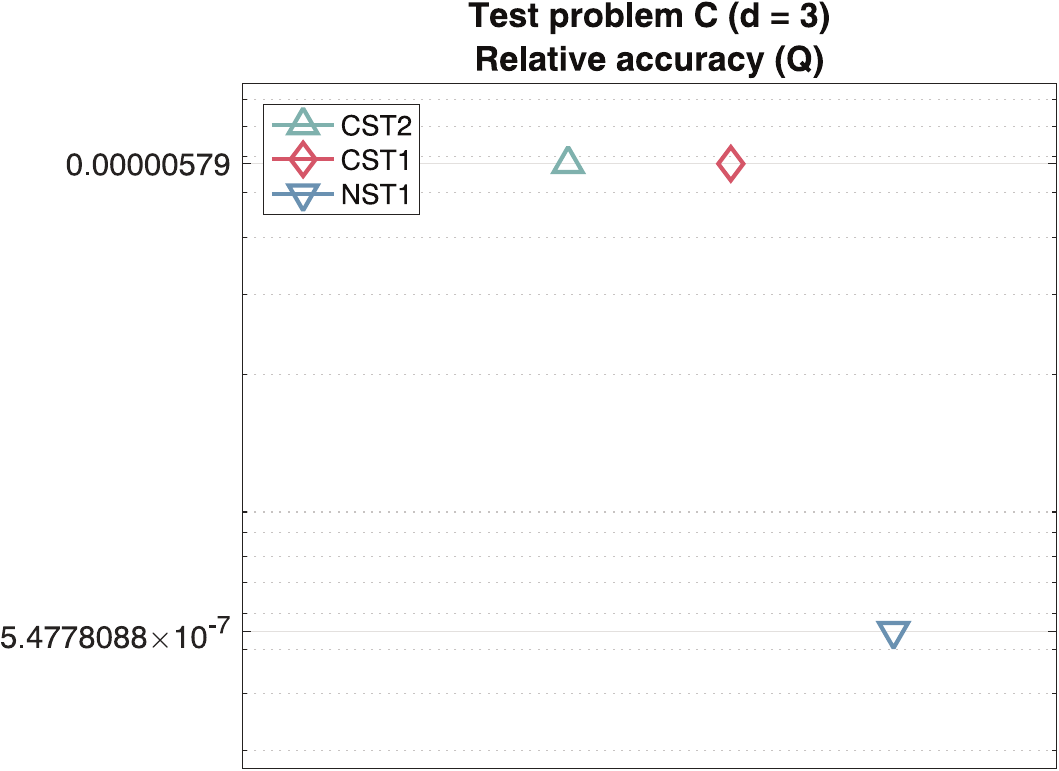} \quad
\includegraphics[width=4.2cm]{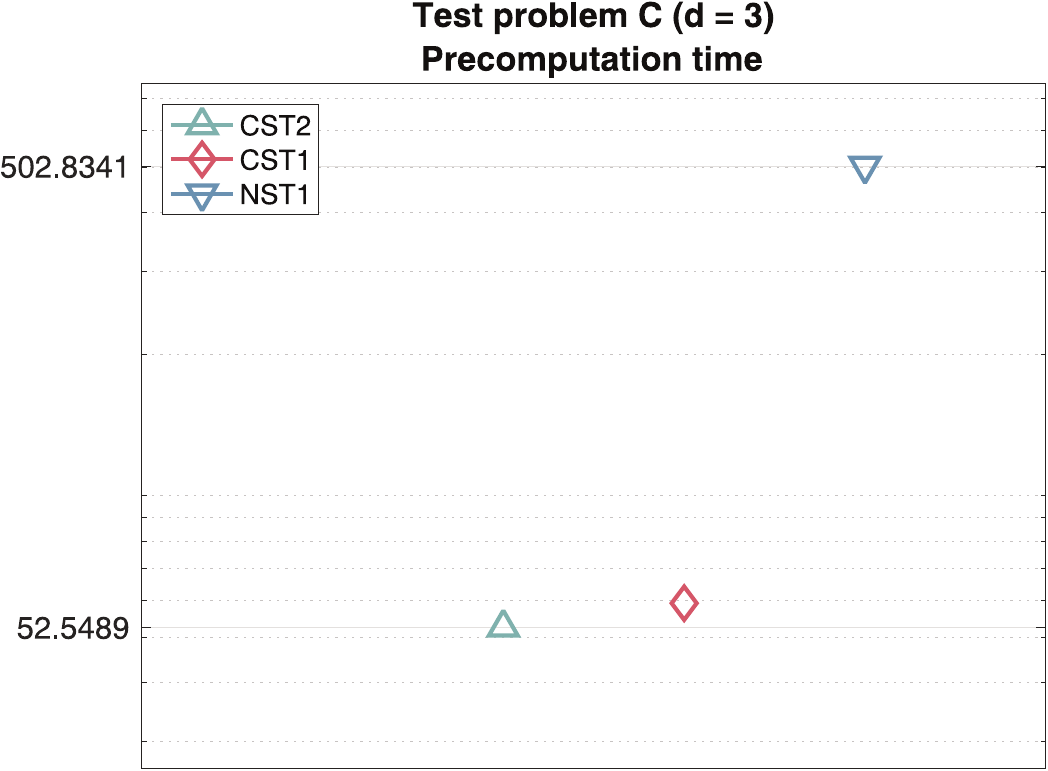} \quad
\includegraphics[width=4.2cm]{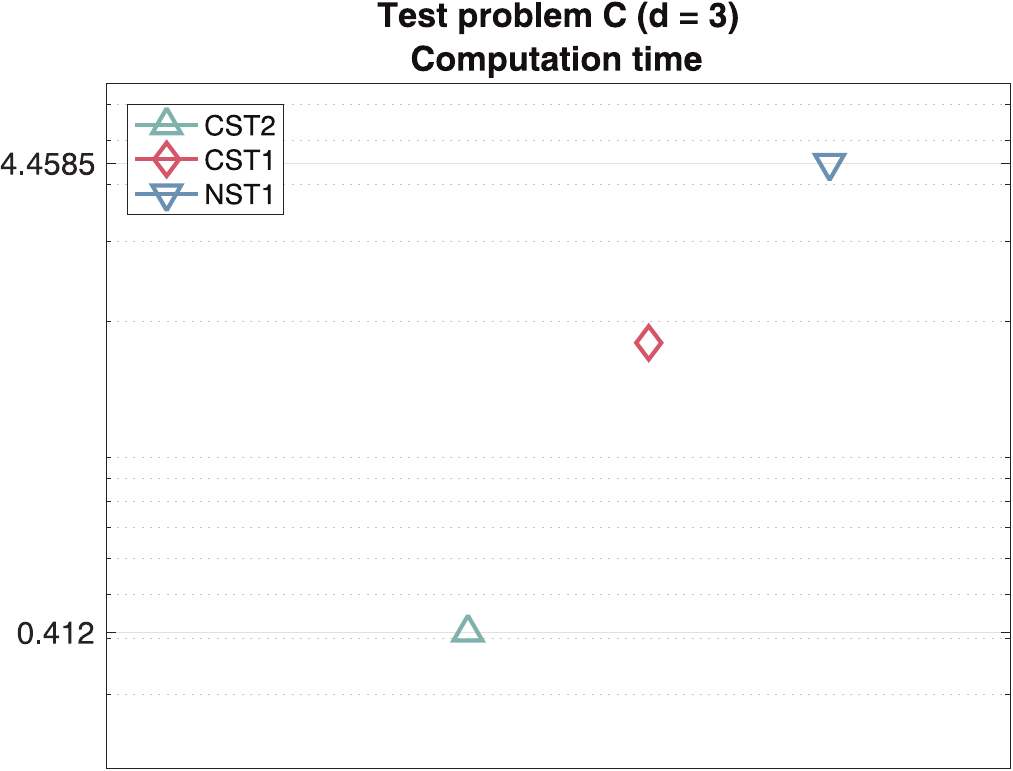} \\[2mm]
\includegraphics[width=4.2cm]{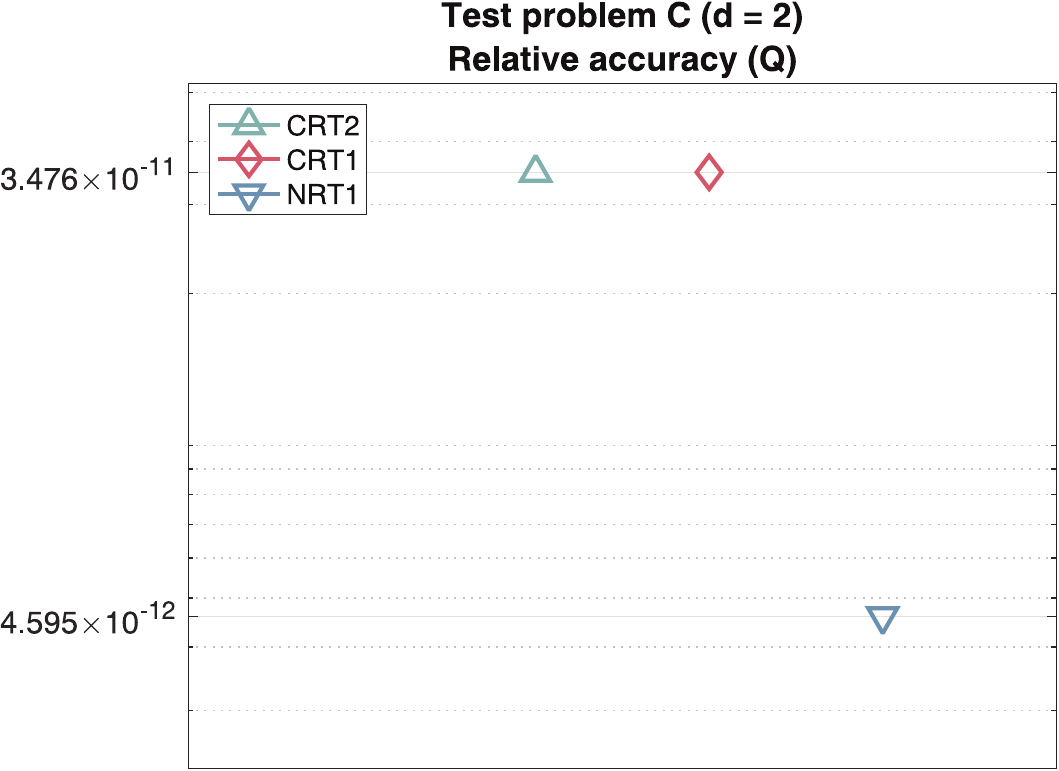} \quad
\includegraphics[width=4.2cm]{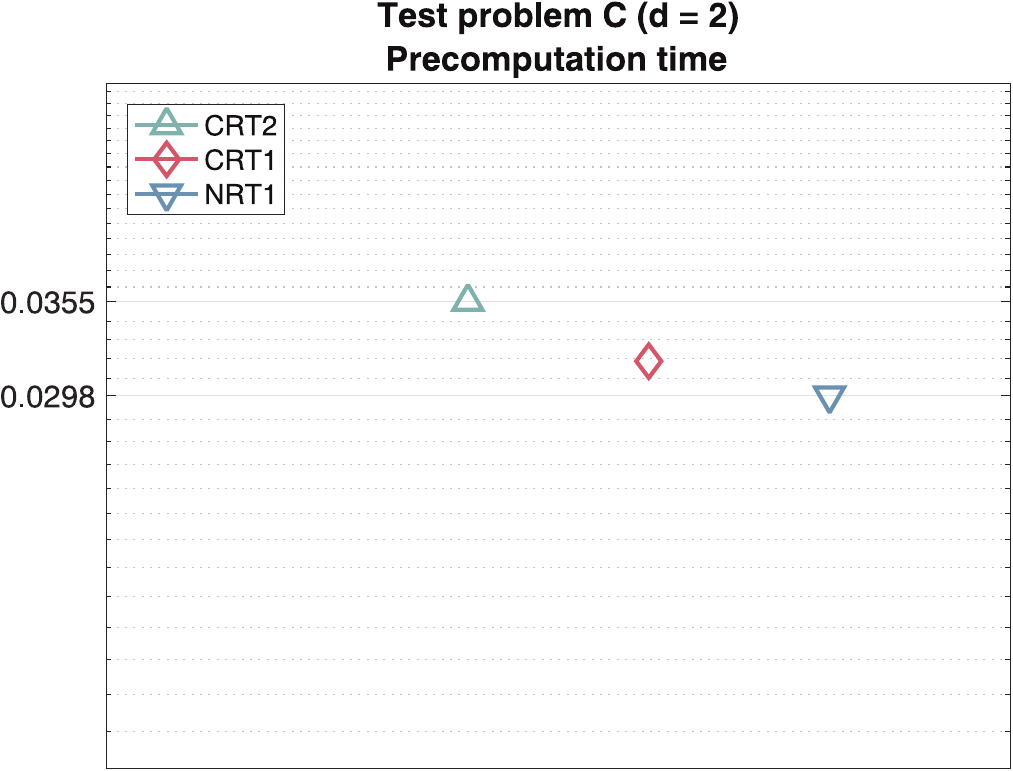} \quad
\includegraphics[width=4.2cm]{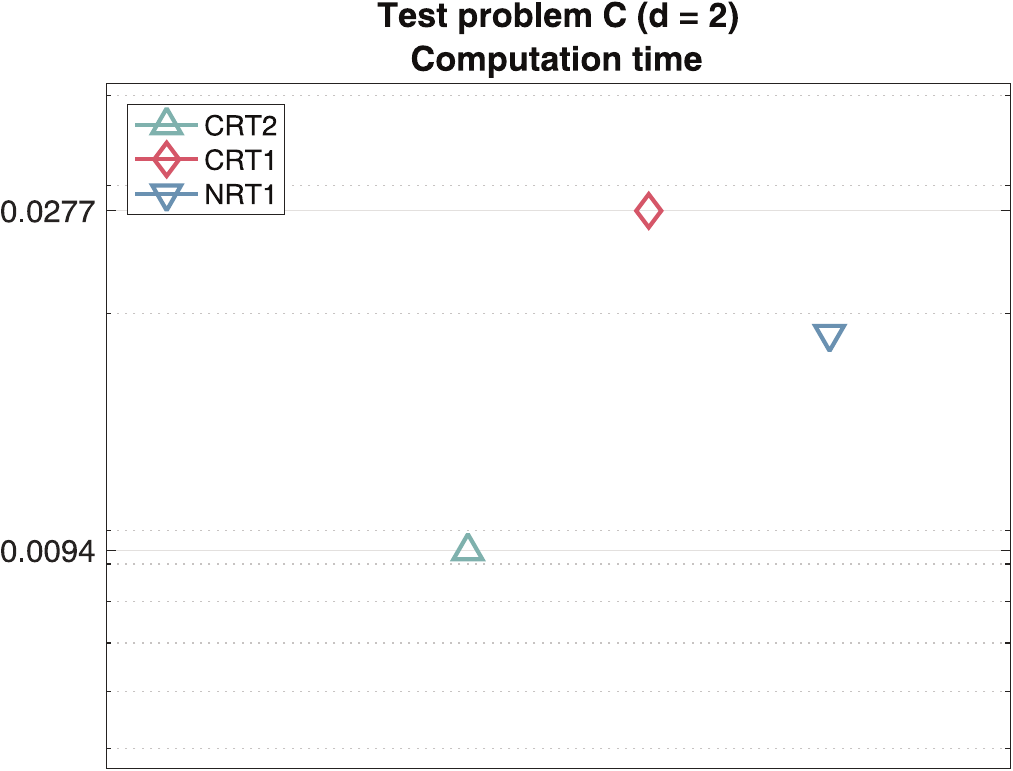} \\[2mm]
\includegraphics[width=4.2cm]{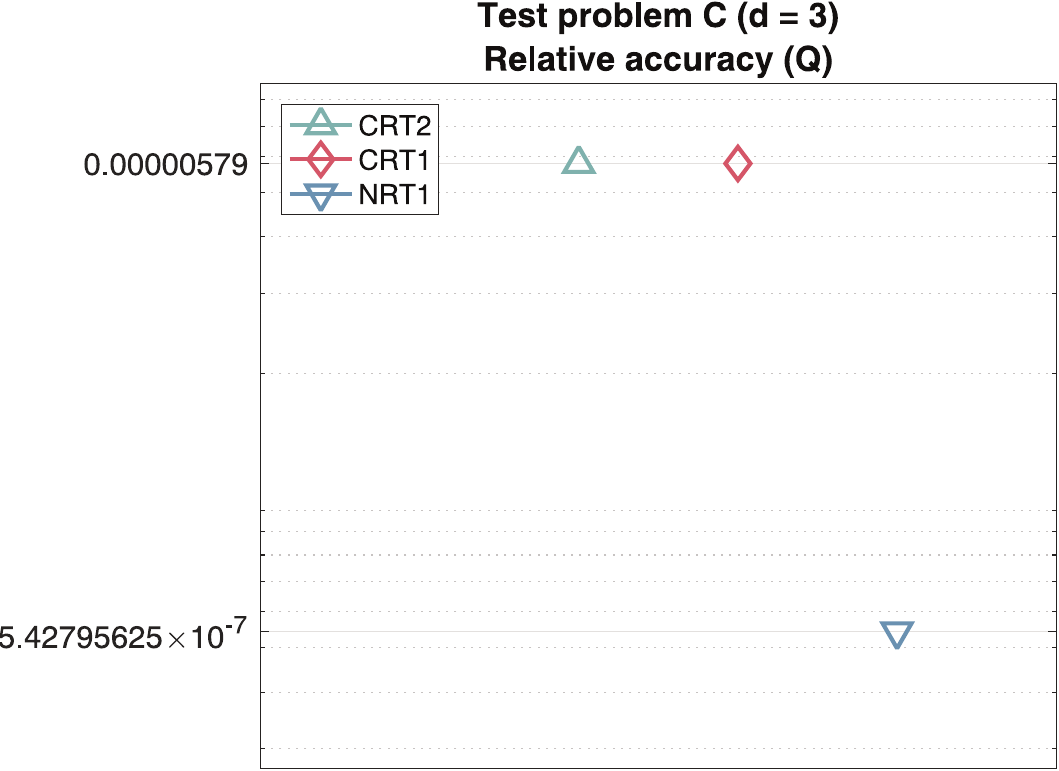} \quad
\includegraphics[width=4.2cm]{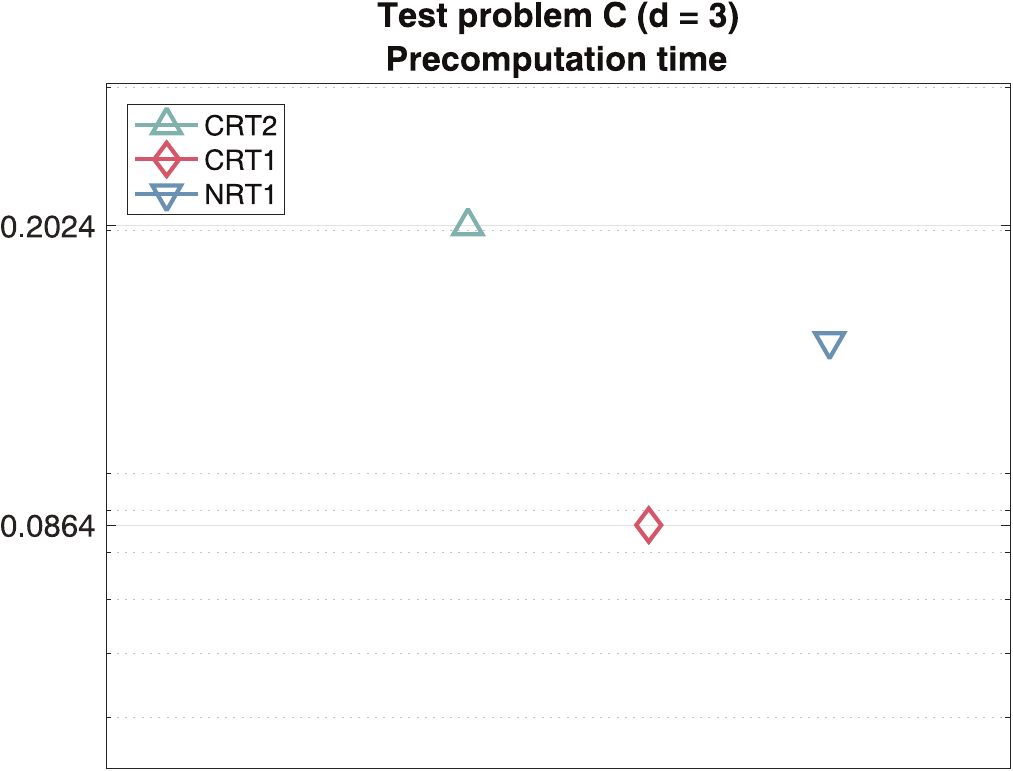} \quad
\includegraphics[width=4.2cm]{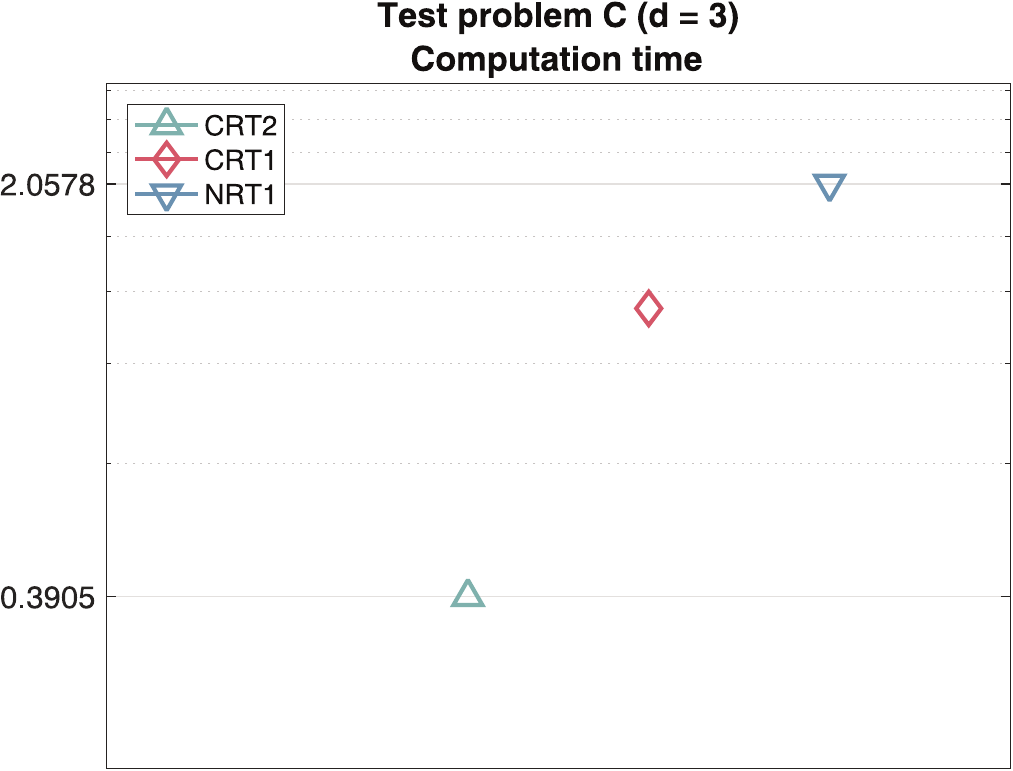}
\caption{Test problem~C (regular integral kernel, unbounded domain, known solution) in two and three dimensions. Numerical comparisons of relative accuracies as well as precomputation and computation times.
General approaches for the evaluation of the associated Landau operator based on the conservative form (CST2, CST1) and the non-conservative form (NST1).
Simplifications to regular kernels (CRT2, CRT1, NRT1).}
\label{fig:ComparisonC}
\end{center}
\end{figure}

\begin{figure}
\begin{center}
\includegraphics[width=4.2cm]{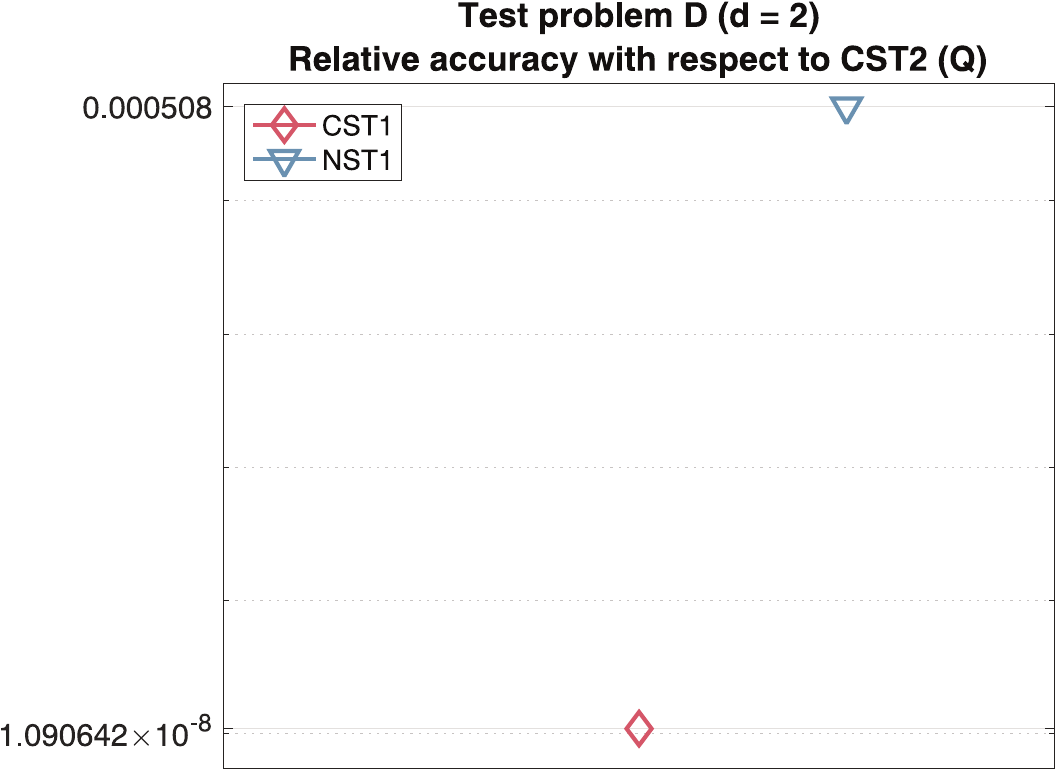} \quad
\includegraphics[width=4.2cm]{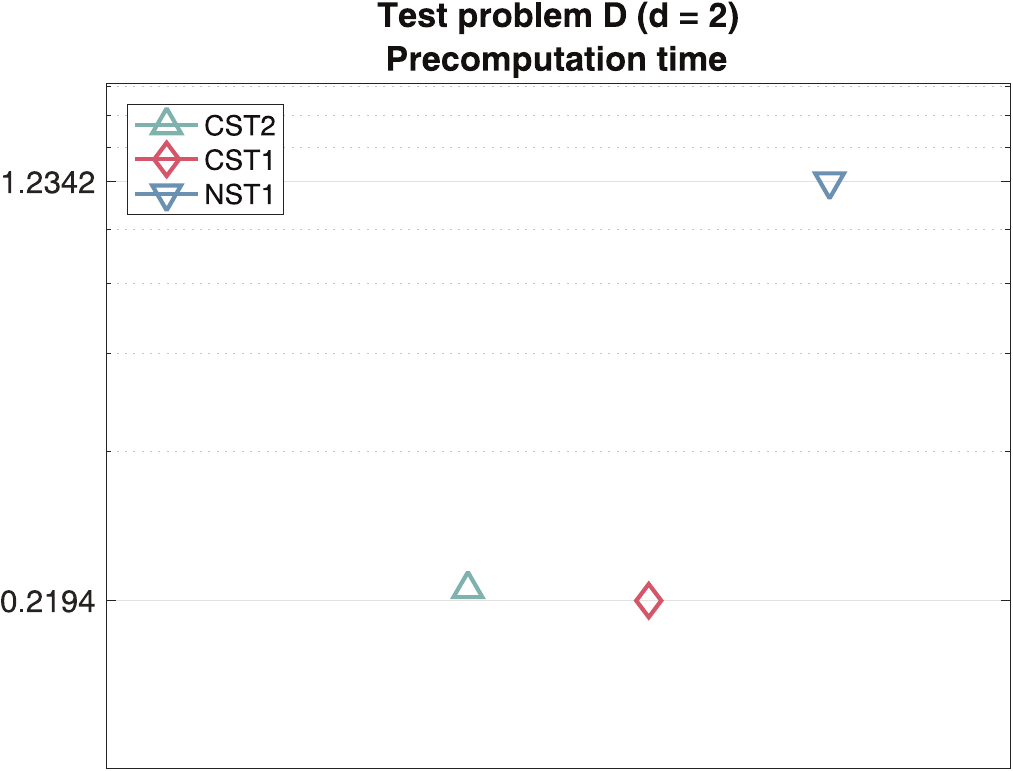} \quad
\includegraphics[width=4.2cm]{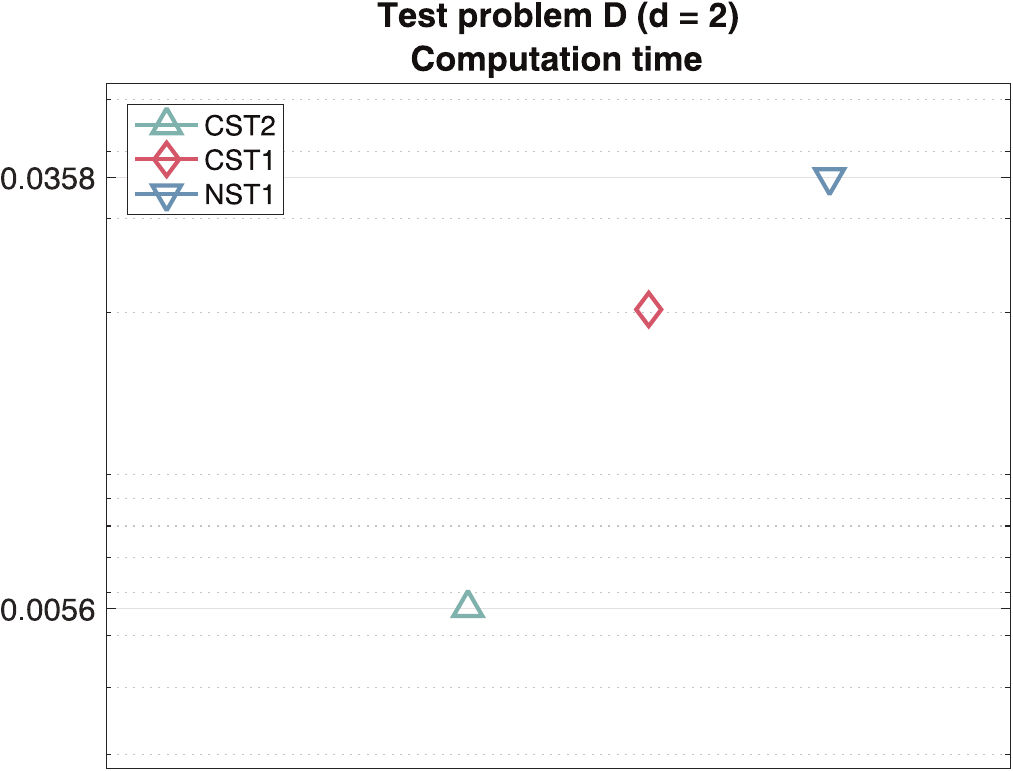} \\[2mm]
\includegraphics[width=4.2cm]{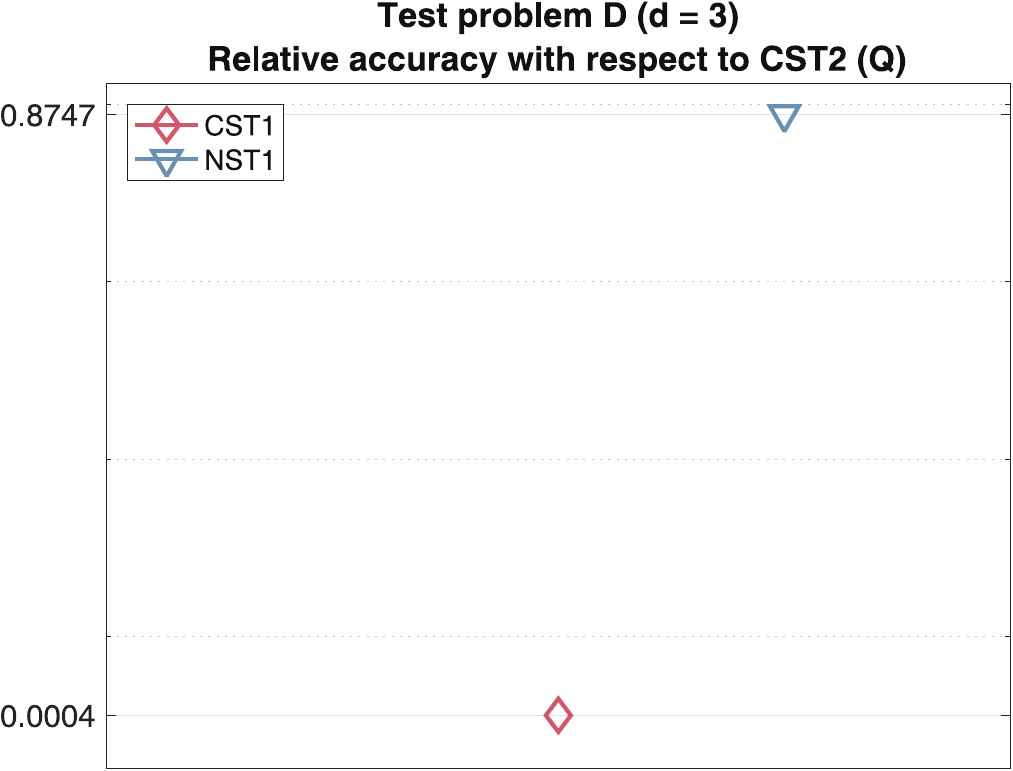} \quad
\includegraphics[width=4.2cm]{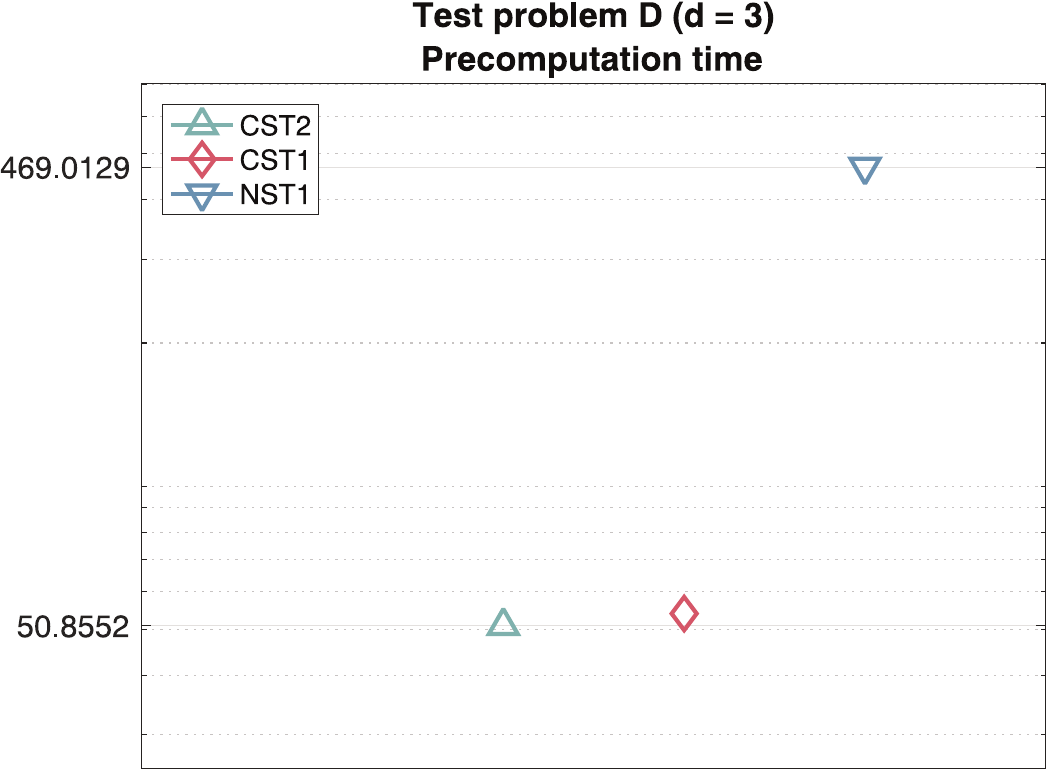} \quad
\includegraphics[width=4.2cm]{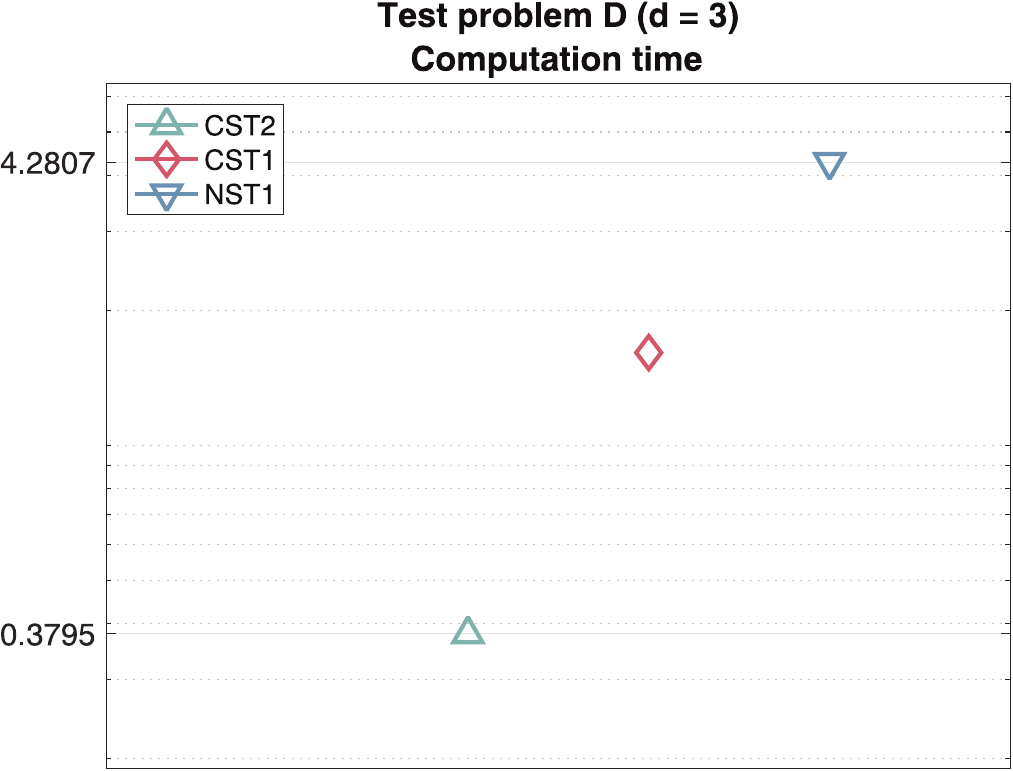}
\caption{Test problem~D in two and three dimensions.
Numerical comparisons of relative accuracies with respect to the first approximation obtained by approach CST2 as well as precomputation and computation times.}
\label{fig:ComparisonD}
\end{center}
\end{figure}

\MyParagraph{Numerical comparisons}
We next contrast the results obtained for the general approaches CST2, CST1, NST1 and their simplifications CRT2, CRT1, NRT1. 
Setting again~\eqref{eq:DomainSymmetric}, we consider a Landau operator involving a constant integral kernel (Maxwellian molecules case, Test problem~A) and a regular kernel (Test problem~C), respectively, see Figures~\ref{fig:ComparisonA} and~\ref{fig:ComparisonC}.
Moreover, we demonstrate the issues of a bounded domain and numerical differentiation (Test problem~B), see Figure~\ref{fig:ComparisonB}. Specifically, the integral operator does not fulfill periodicity requirements and hence is not well represented by Fourier series.
In addition, we compare computation times and relative accuracies with respect to approach CST2 when evaluating the Landau operator involving a singular kernel (Test problem~D), see Figure~\ref{fig:ComparisonD}.
Presumably, in the latter case, the main source of approximation errors is linked to numerical quadrature.   
Due to the fact that the non-conservative approach NST1 does not preserve mass and requires quadrature approximations of the singular integral kernel and its first-order derivatives, we select the conservative approach CST2 with a reduced computation time for further numerical tests.

\begin{figure}
\begin{center}
\includegraphics[width=4.2cm]{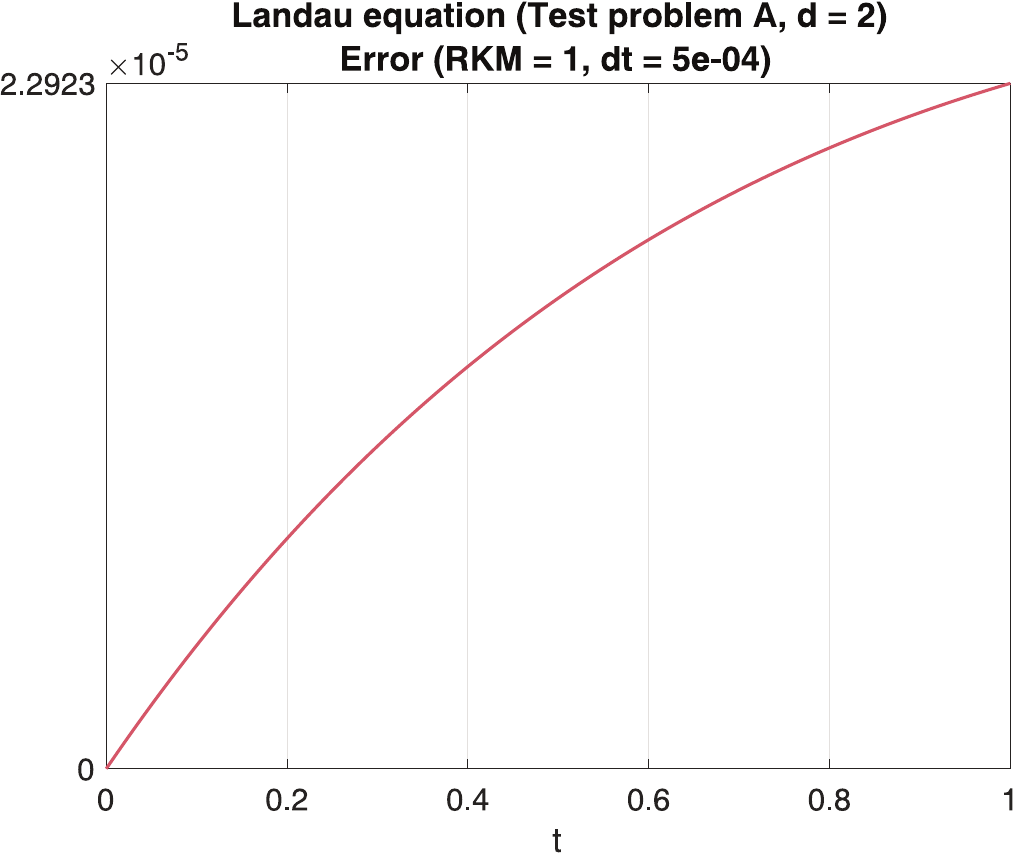} \quad
\includegraphics[width=4.2cm]{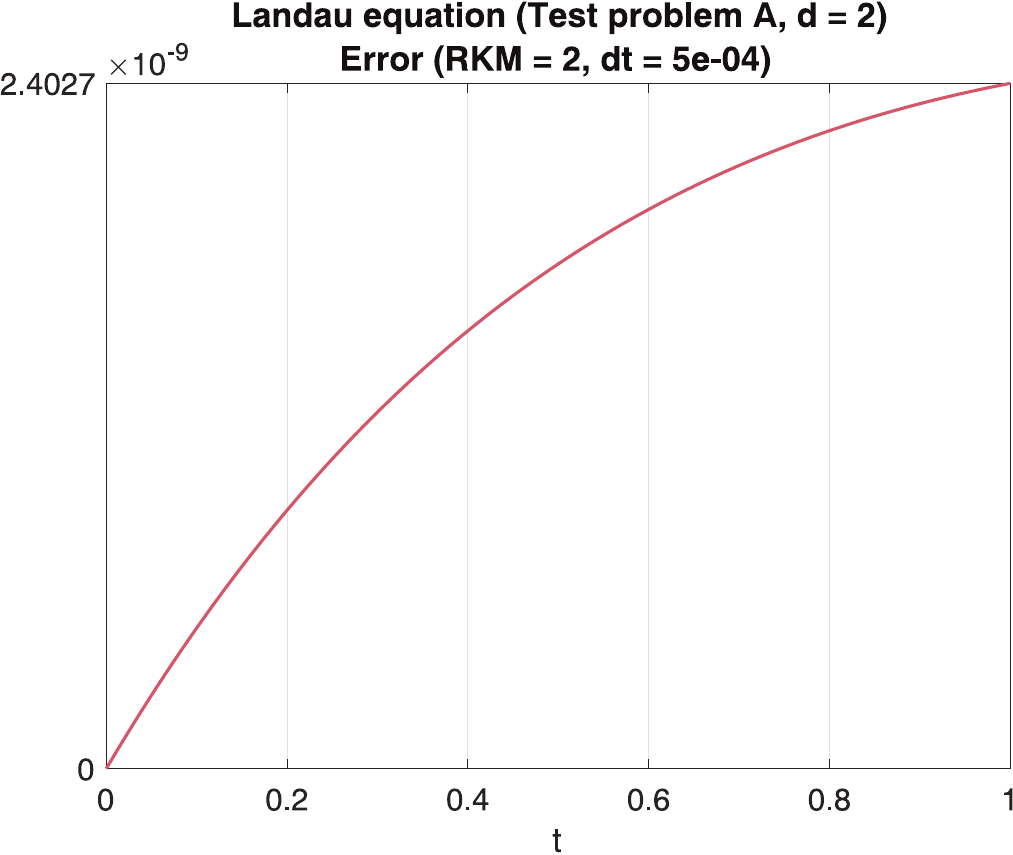} \quad
\includegraphics[width=4.2cm]{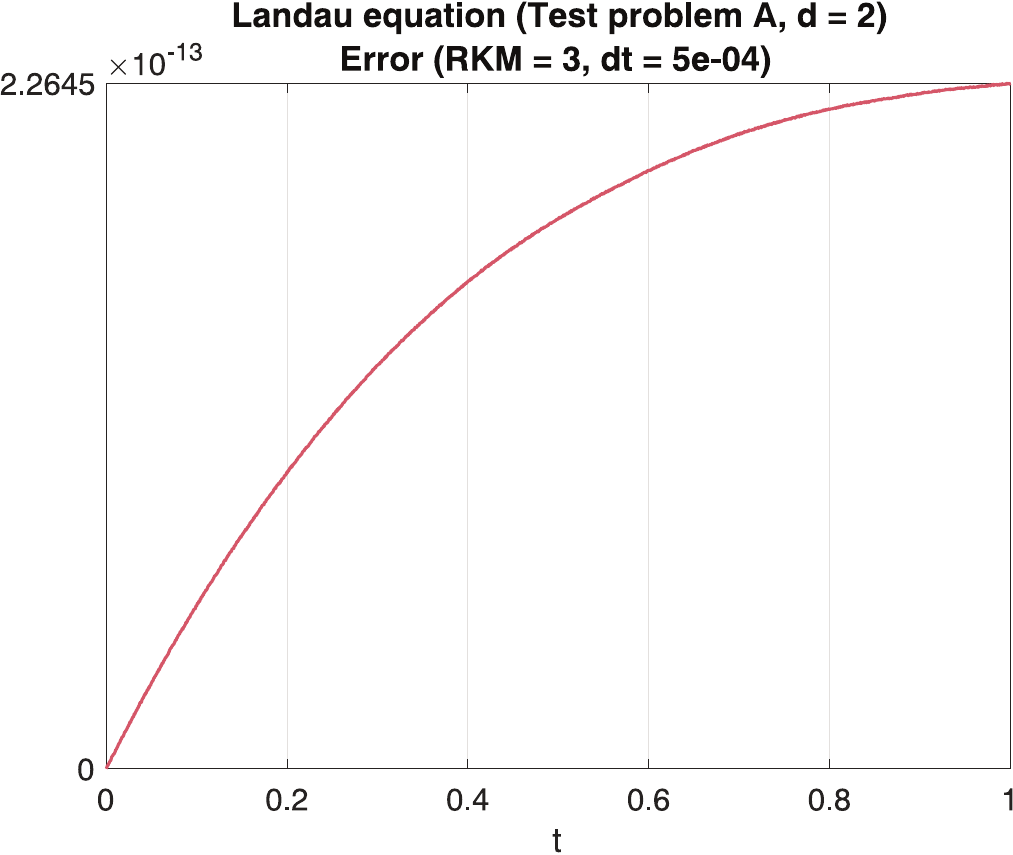} \\[2mm]
\includegraphics[width=4.2cm]{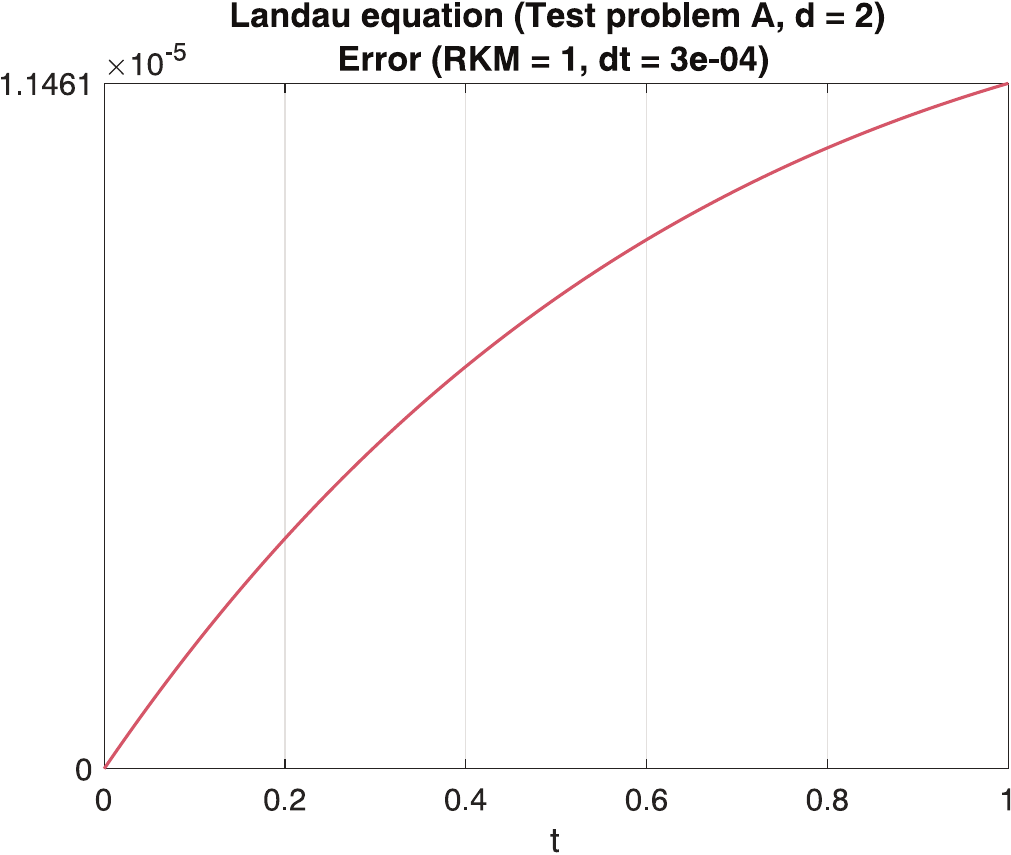} \quad
\includegraphics[width=4.2cm]{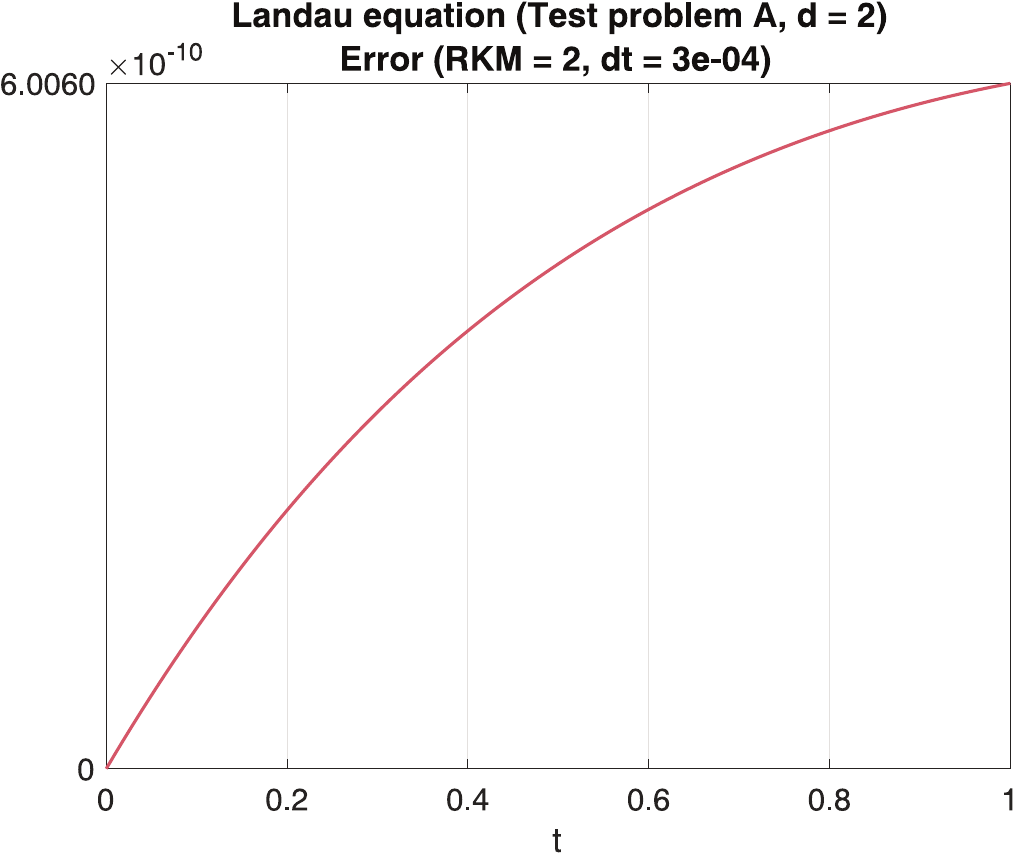} \quad
\includegraphics[width=4.2cm]{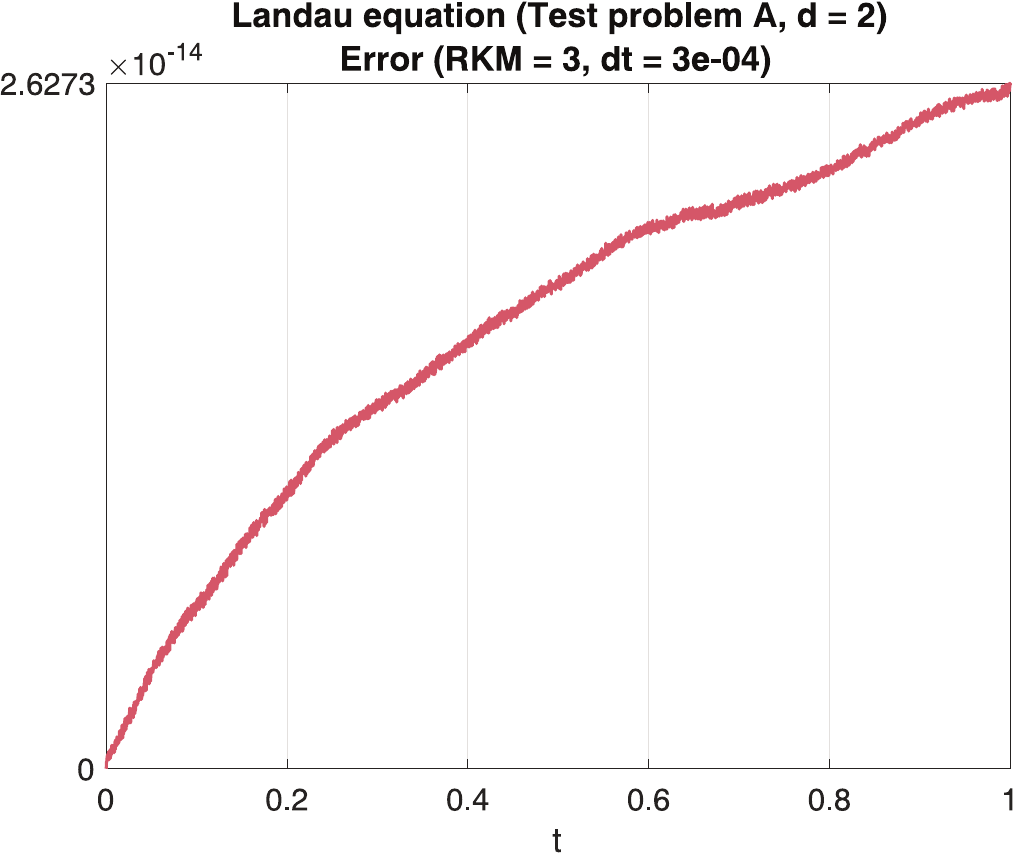} 
\caption{Test problem A (Maxwellian molecules case) in two dimensions.
Time integration of the Landau equation based on $100 \times 100$ Fourier functions and explicit Runge--Kutta methods of orders $p \in \{1, 2, 3\}$.
The absolute errors over time with respect to the known solution values, obtained for a certain time increment and a reduced increment, confirm the orders of convergence.}
\label{fig:ErrorA2d}
\end{center}
\end{figure}

\begin{figure}
\begin{center}
\includegraphics[width=4.2cm]{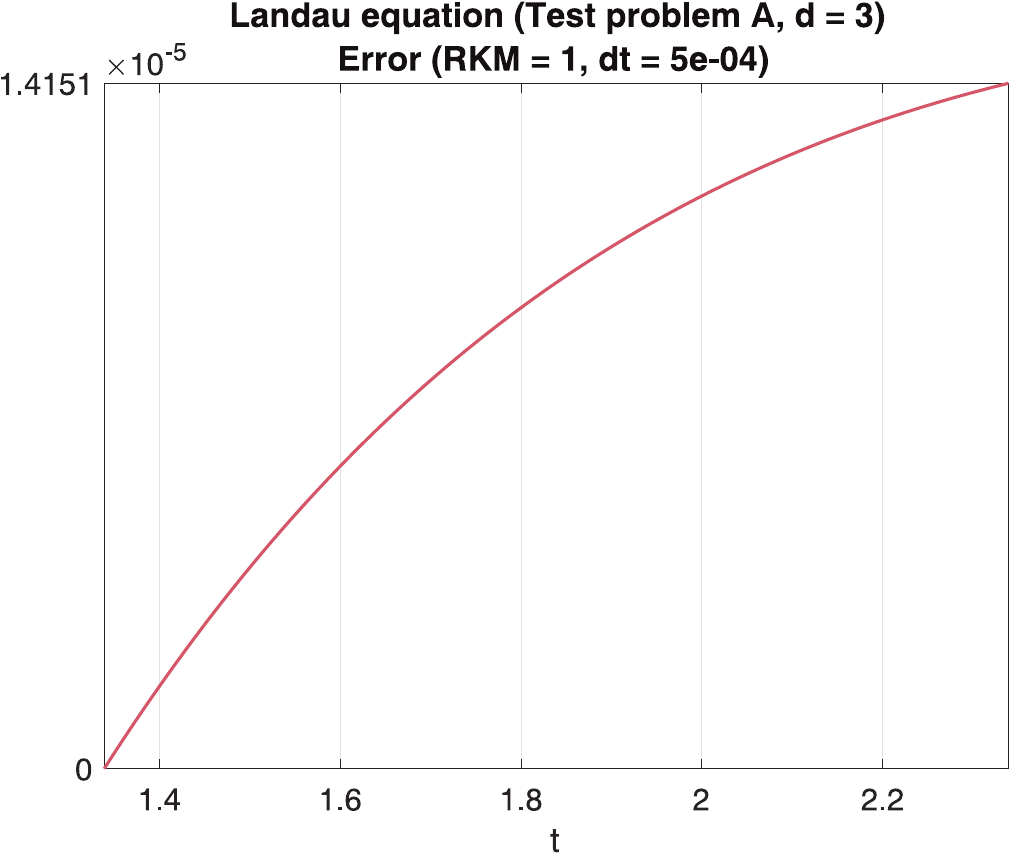} 
\quad
\includegraphics[width=4.2cm]{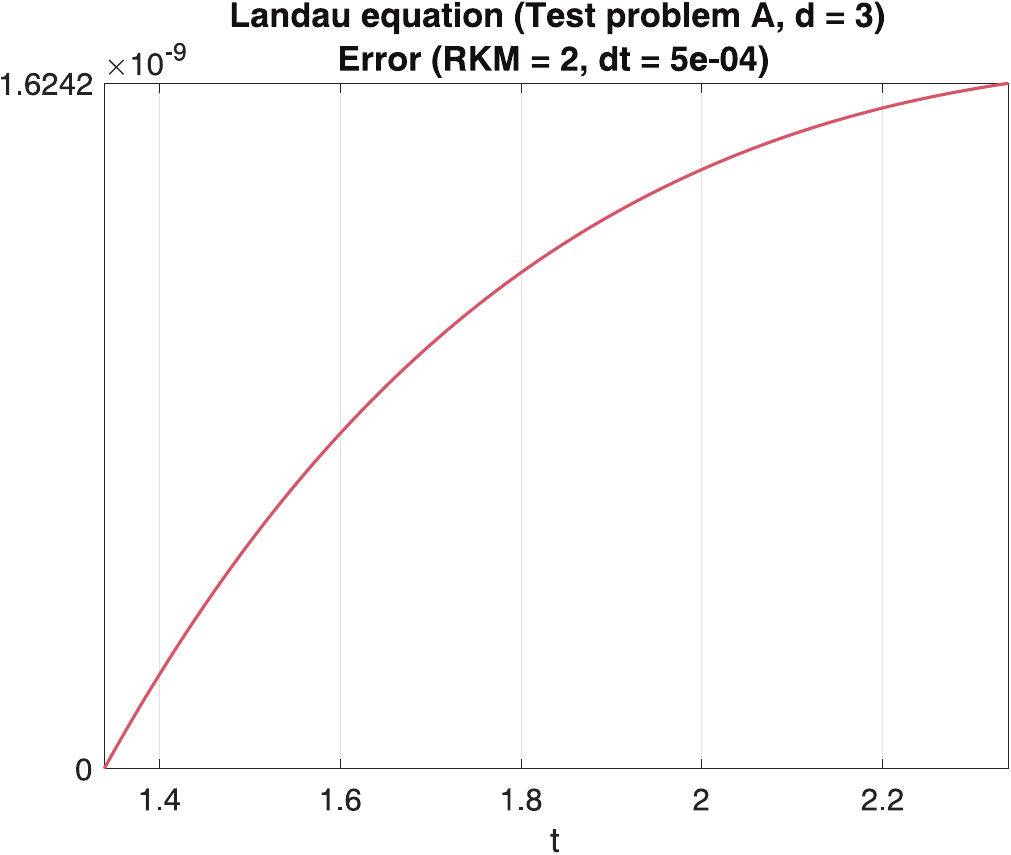} 
\quad
\includegraphics[width=4.2cm]{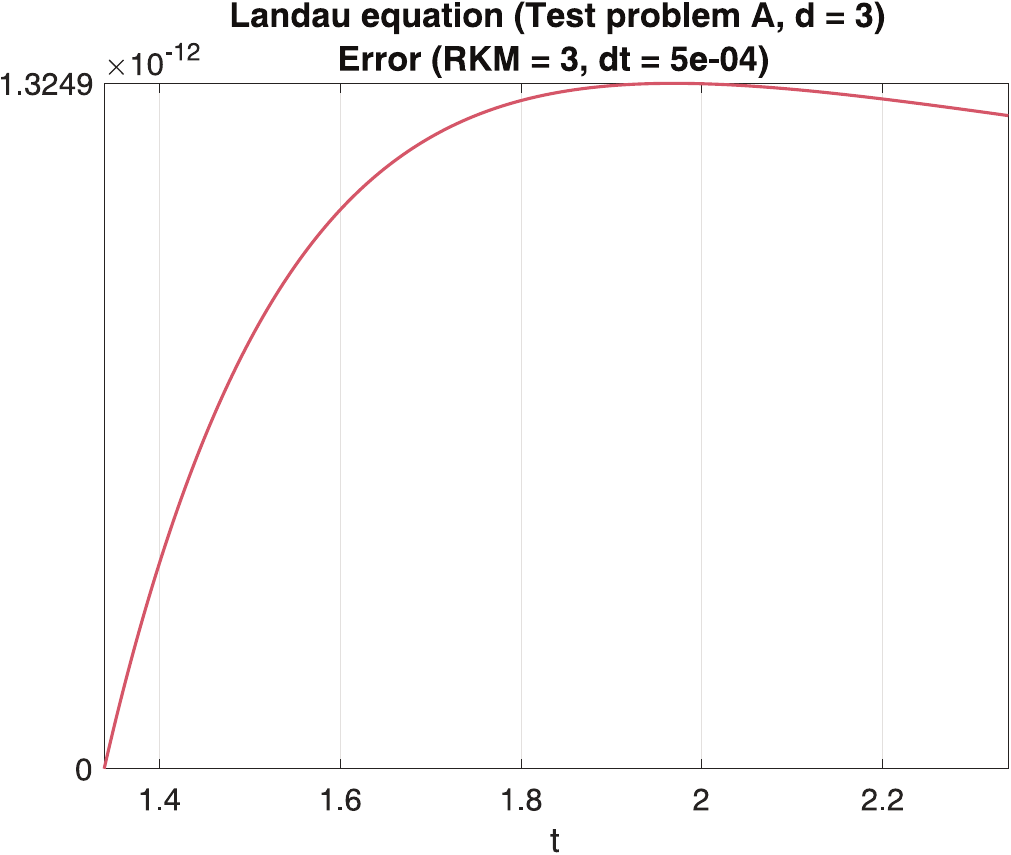} 
\caption{Test problem A (Maxwellian molecules case) in three dimensions.
Time integration of the Landau equation based on $64 \times 64 \times 64$ Fourier functions and explicit Runge--Kutta methods of orders $p \in \{1, 2, 3\}$.
Absolute errors over time with respect to the known solution values, obtained for a certain time increment.}
\label{fig:ErrorA3d}
\end{center}
\end{figure}

\begin{figure}
\begin{center}
\includegraphics[width=4.2cm]{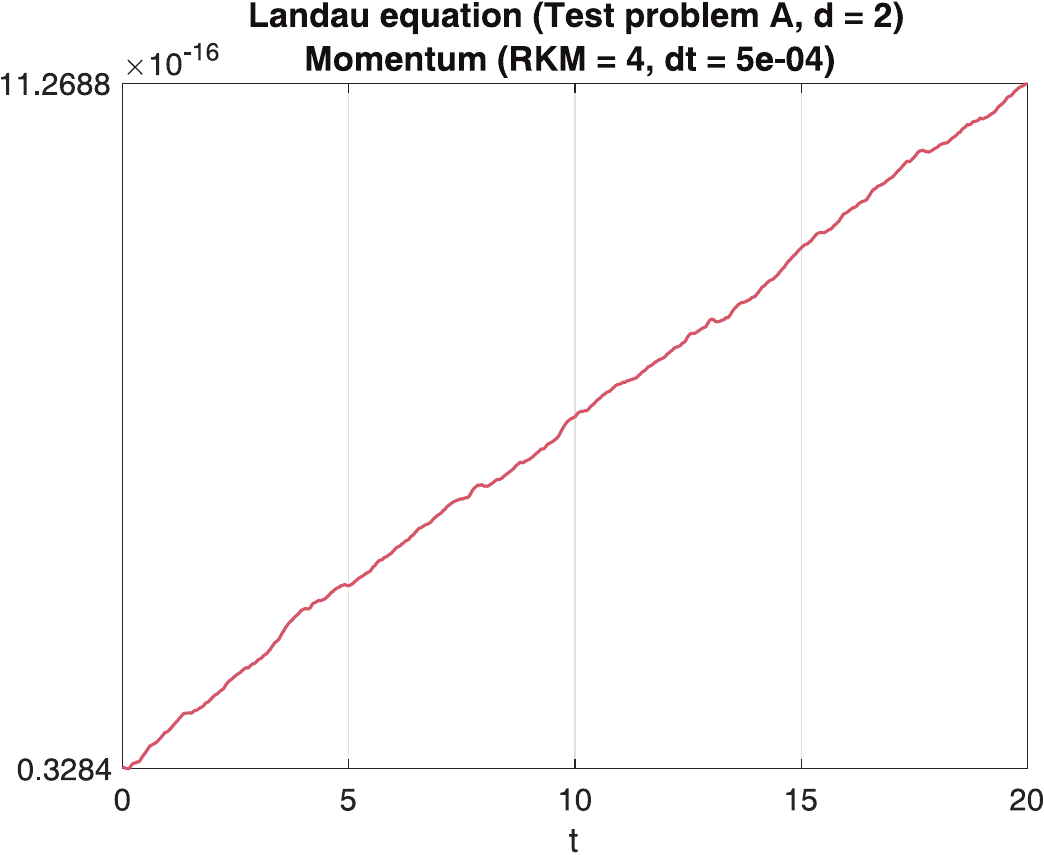} 
\quad
\includegraphics[width=4.2cm]{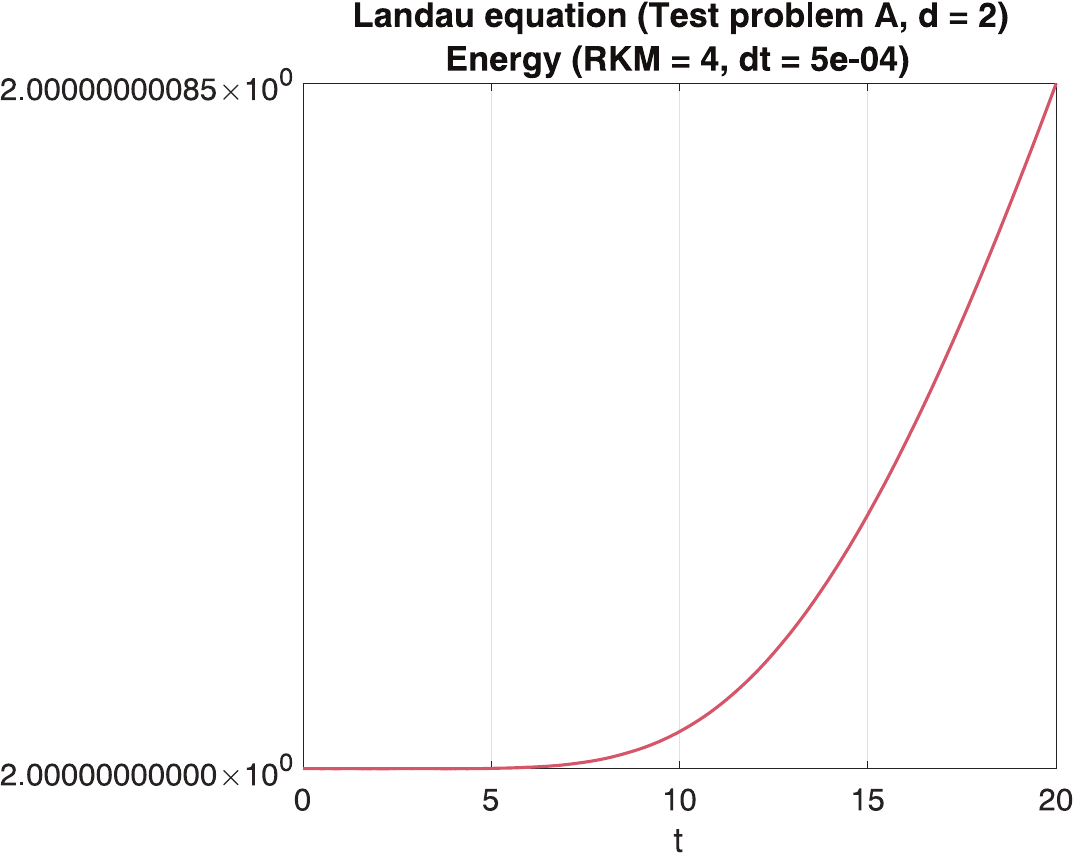} 
\quad
\includegraphics[width=4.2cm]{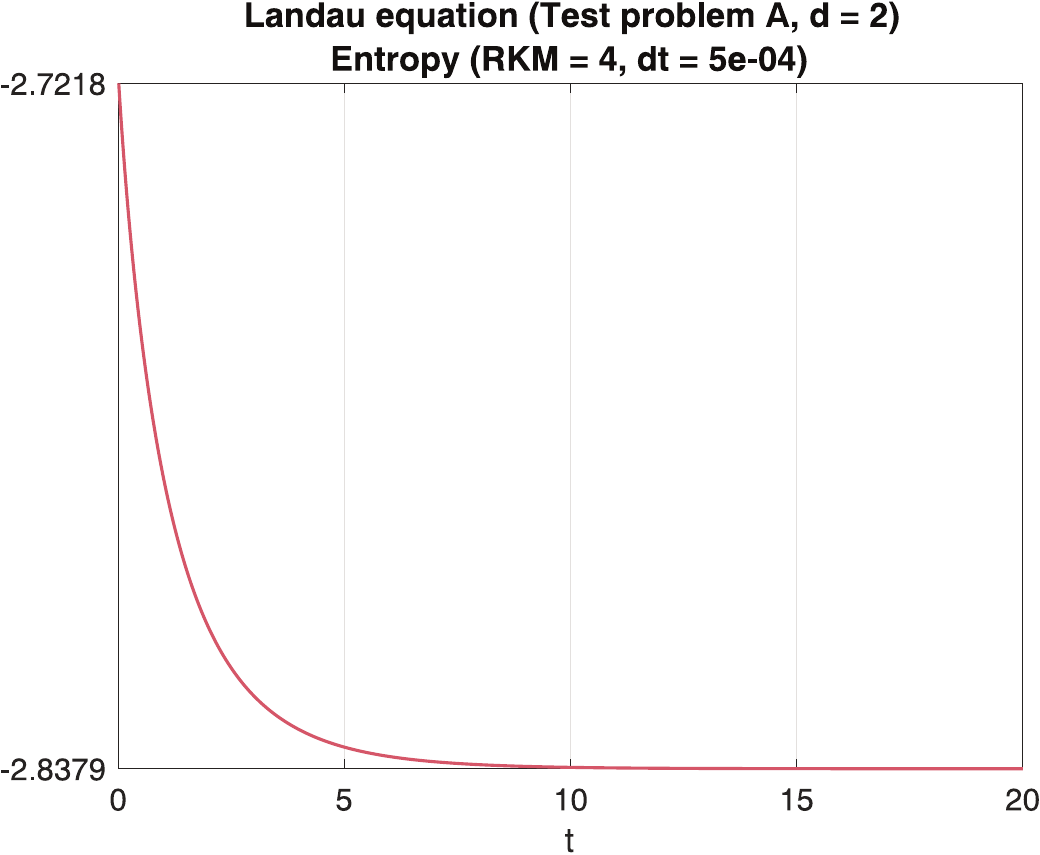} \\[2mm]
\includegraphics[width=4.2cm]{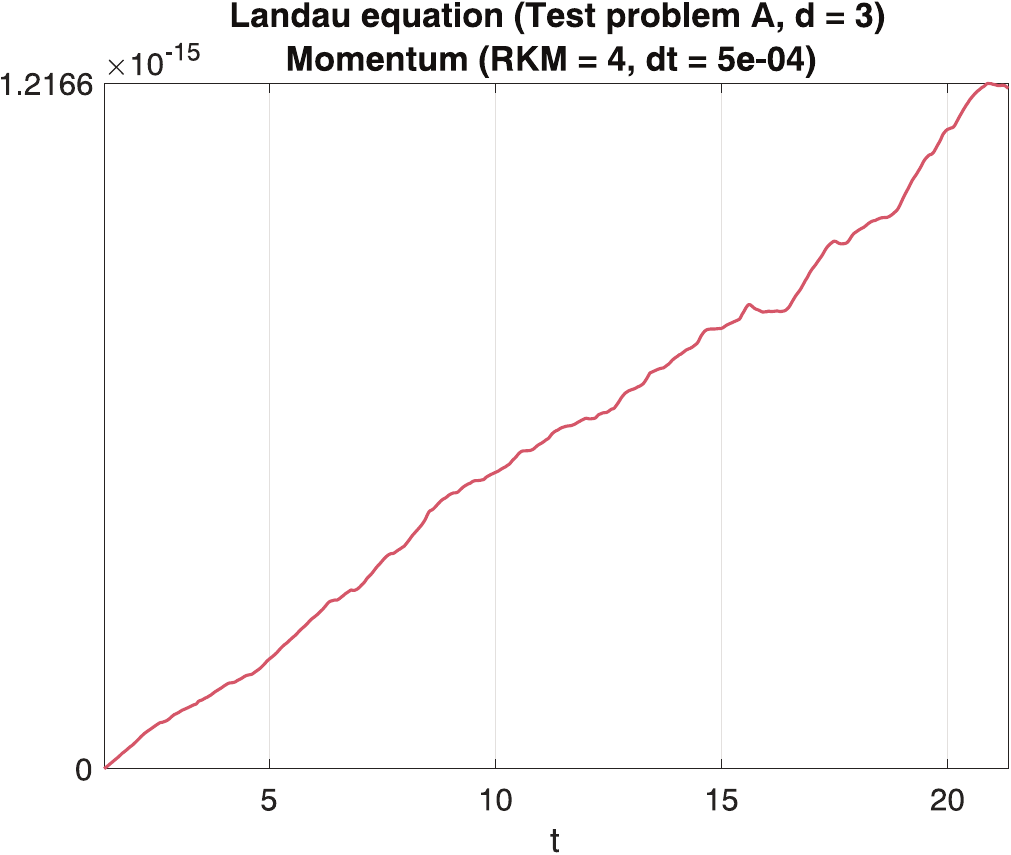} 
\quad
\includegraphics[width=4.2cm]{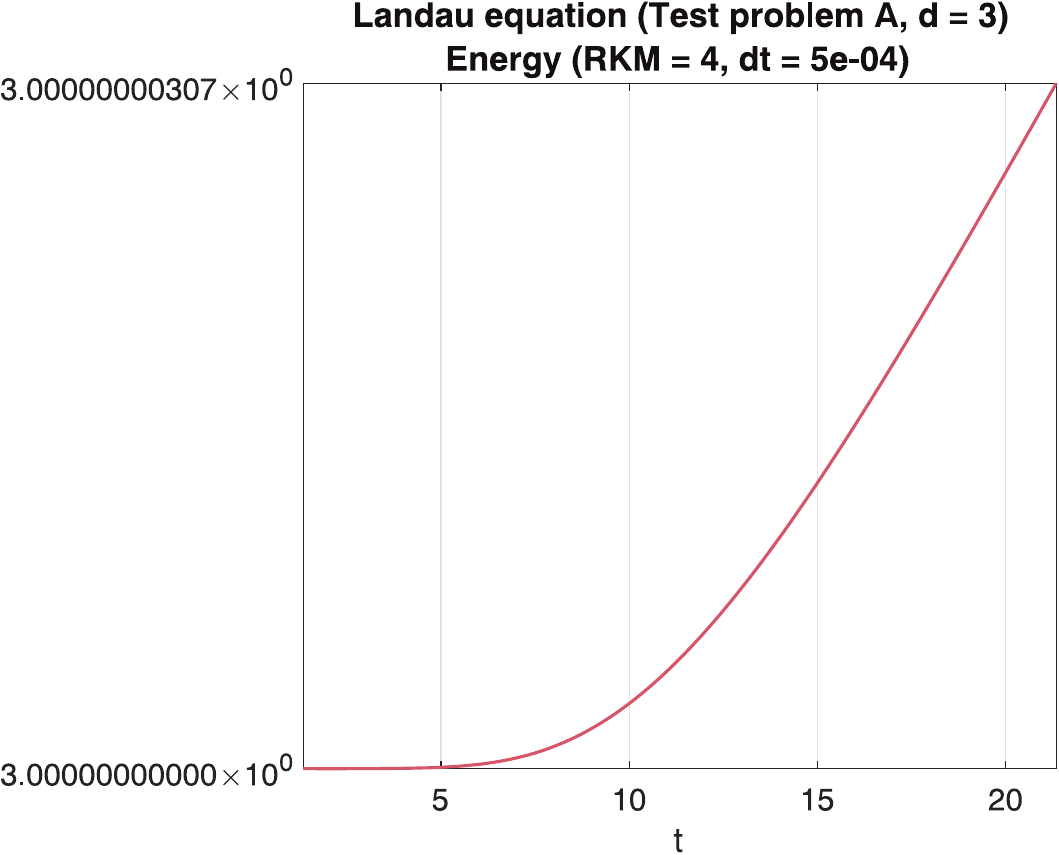} 
\quad
\includegraphics[width=4.2cm]{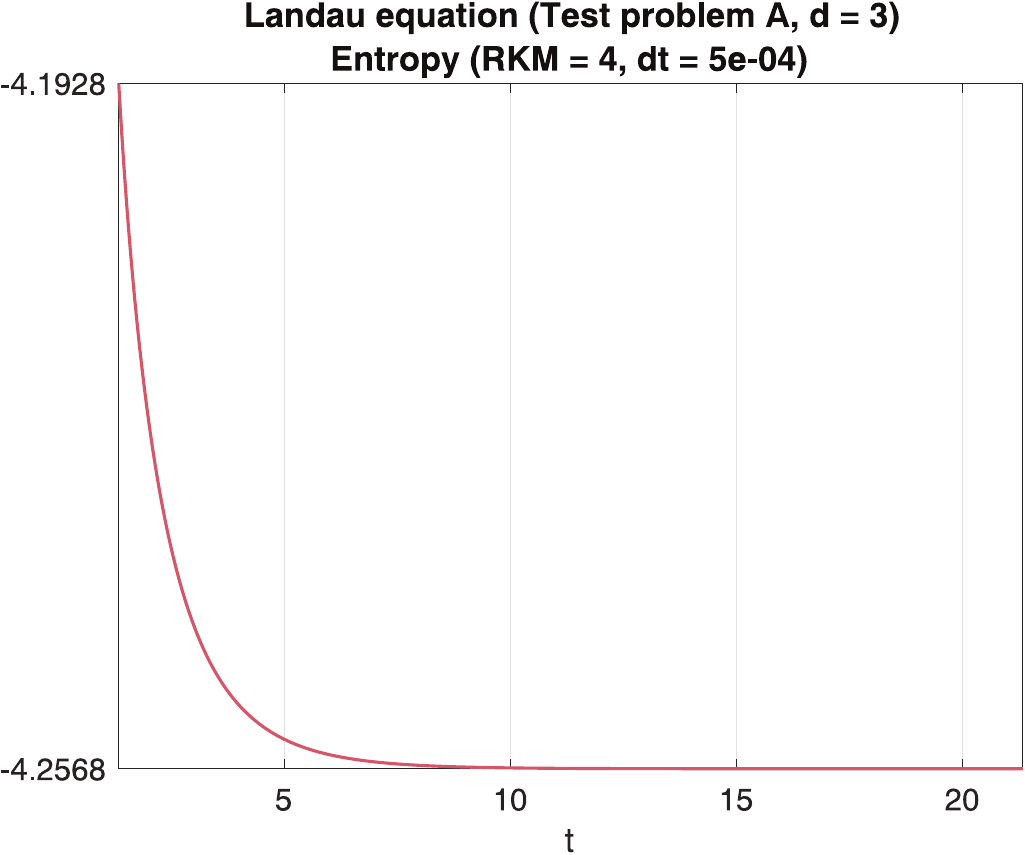} 
\caption{Test problem A (Maxwellian molecules case) in two and three dimensions.
Time integration of the Landau equation based on an explicit Runge--Kutta method of order $p = 4$.
Mass is conserved. First momentum, energy, and entropy over time.}
\label{fig:MassEtcA}
\end{center}
\end{figure}

\begin{figure}
\begin{center}
\includegraphics[width=4.2cm]{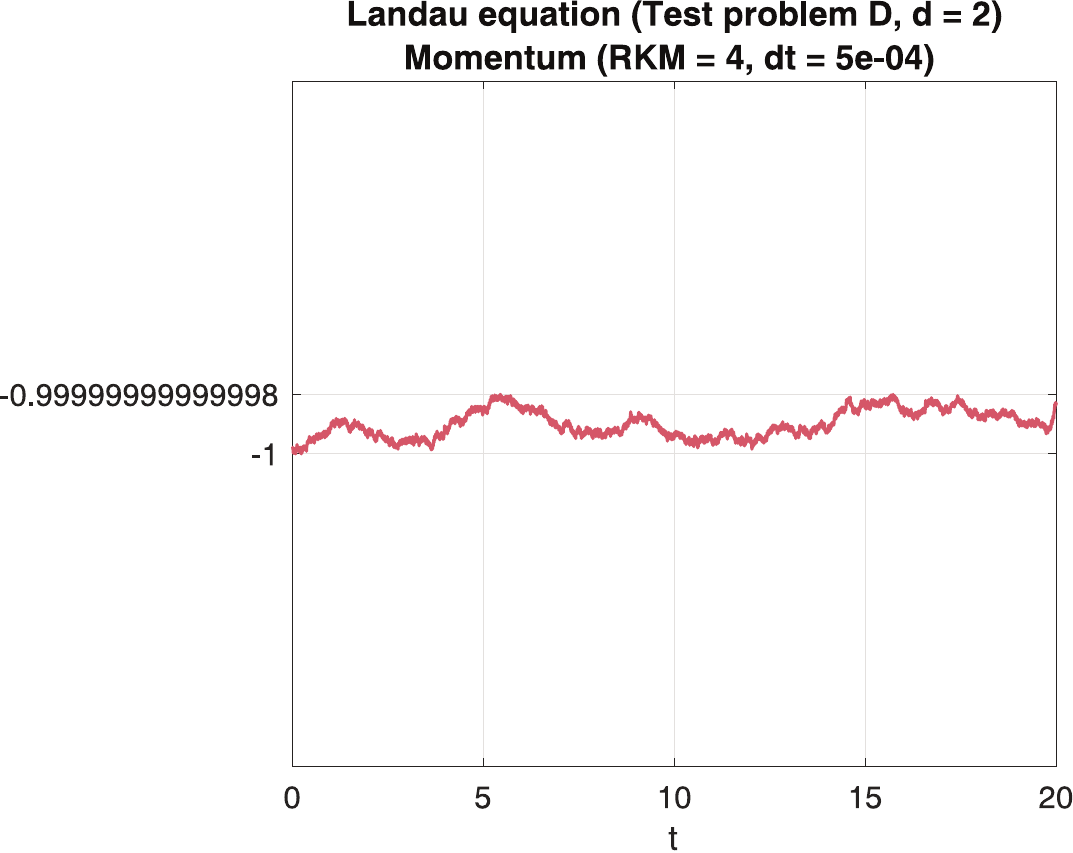} 
\quad
\includegraphics[width=4.2cm]{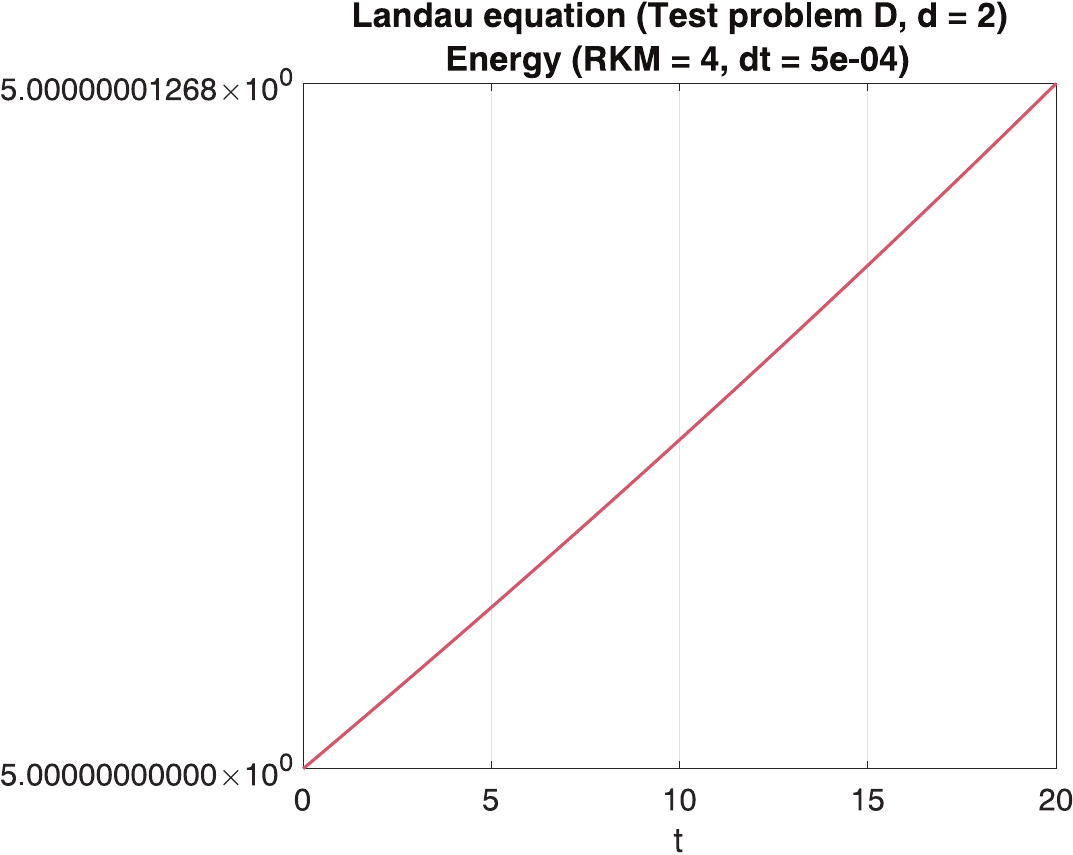} 
\quad
\includegraphics[width=4.2cm]{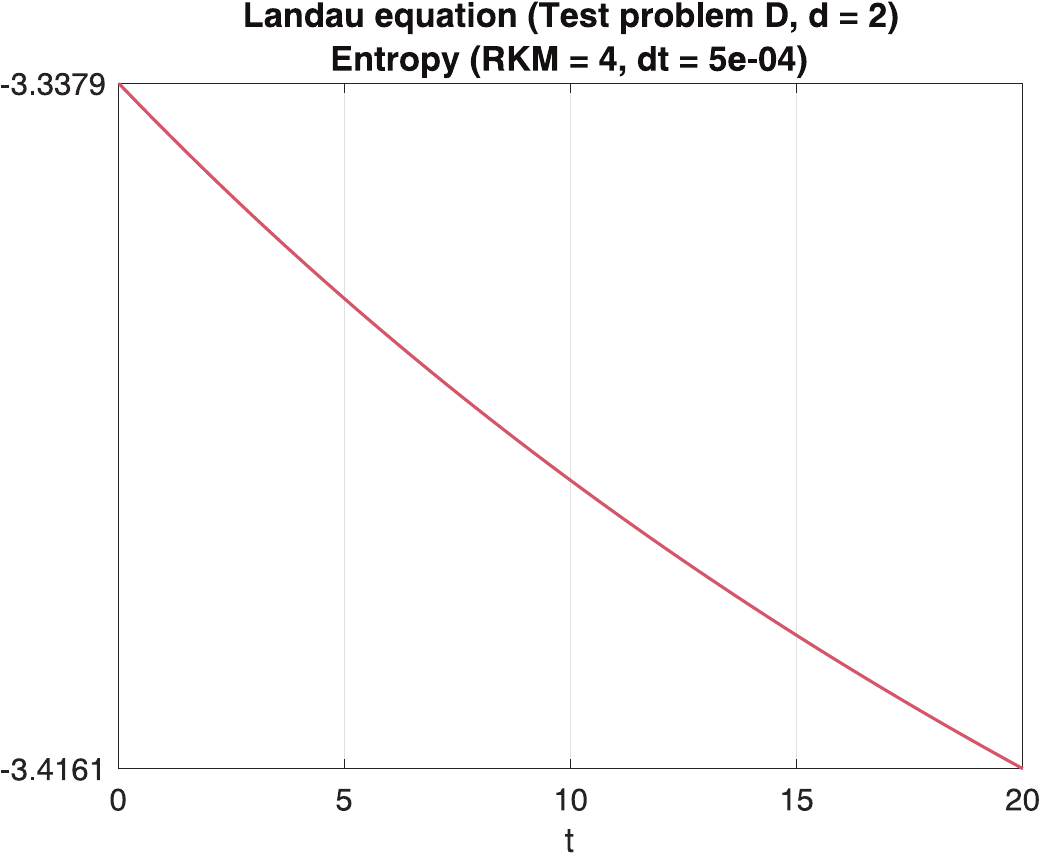} \\[2mm]
\includegraphics[width=4.2cm]{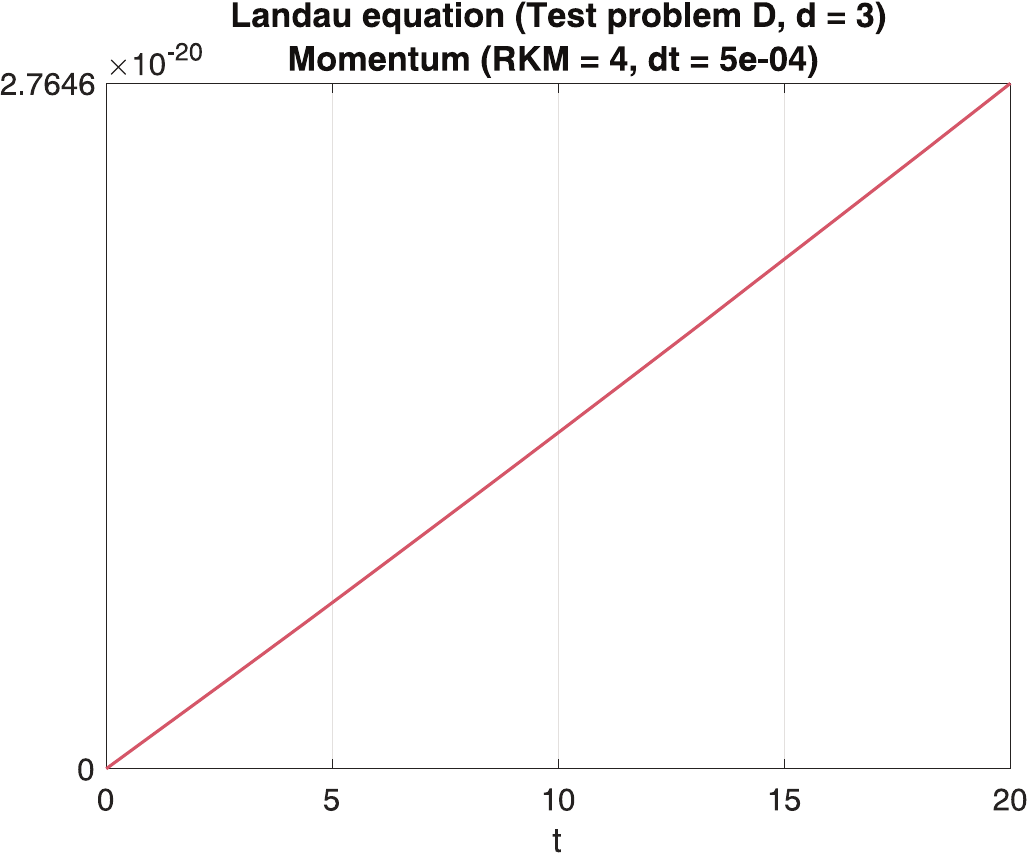} 
\quad
\includegraphics[width=4.2cm]{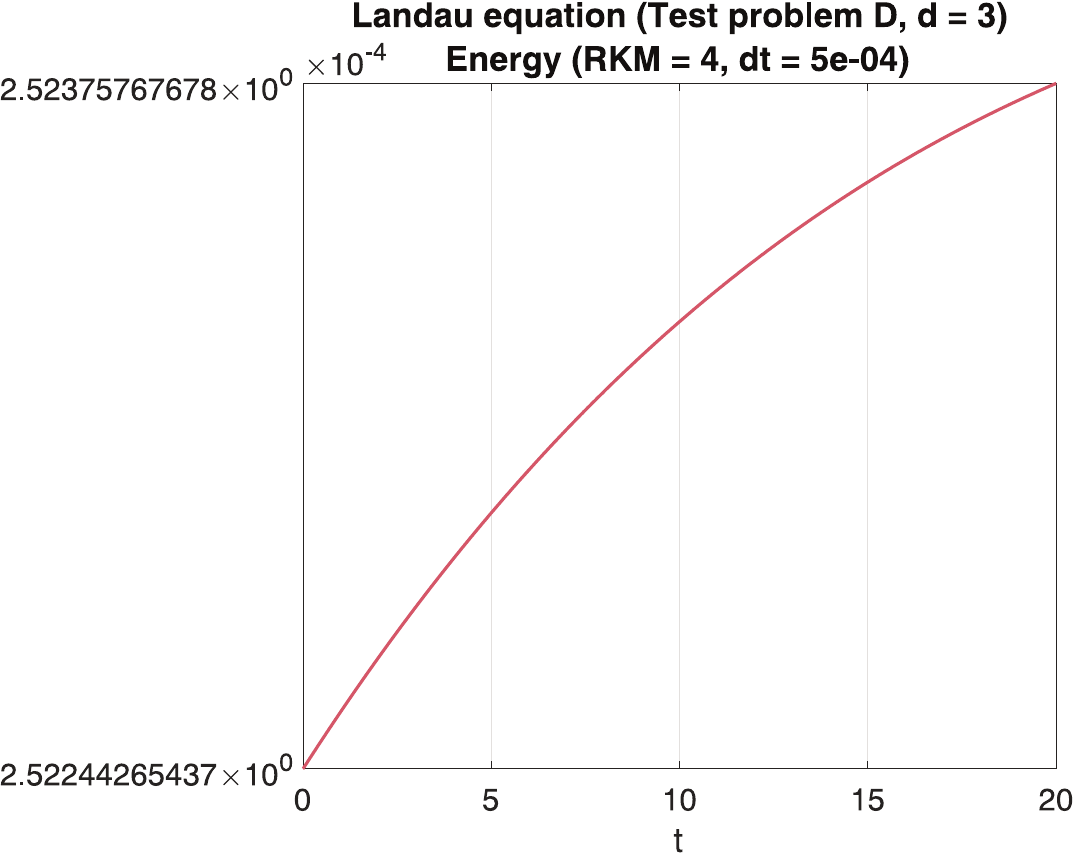} 
\quad
\includegraphics[width=4.2cm]{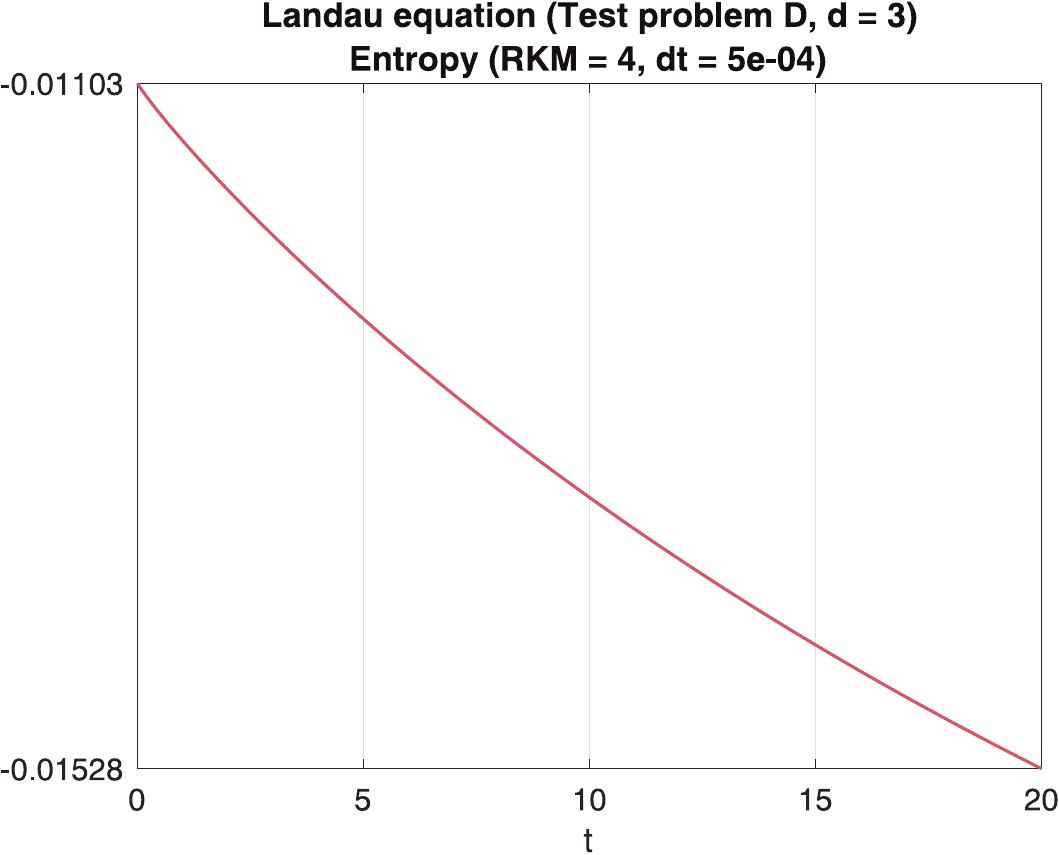} 
\caption{Test problem D in two and three dimensions.
Time integration of the Landau equation based on an explicit Runge--Kutta method of order $p = 4$.
Mass is conserved. First momentum, energy, and entropy over time.}
\label{fig:MassEtcD}
\end{center}
\end{figure}

\begin{figure}
\begin{center}
\includegraphics[width=4.2cm]{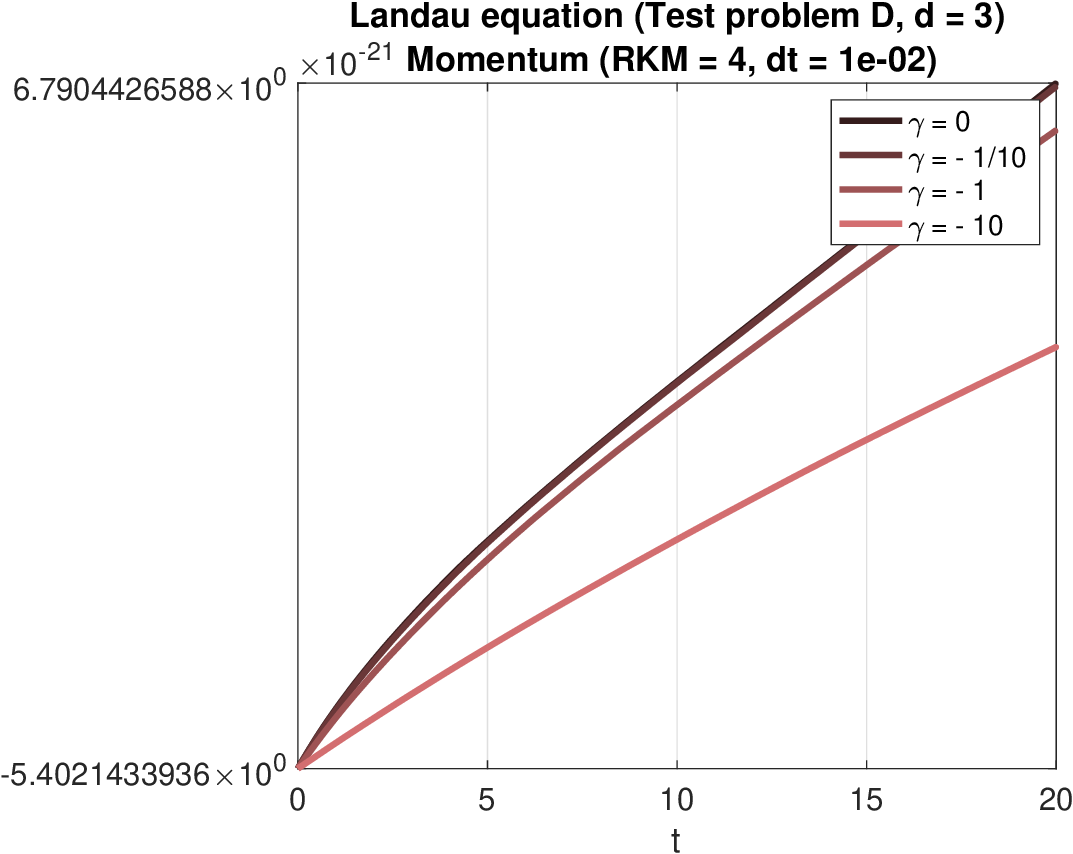} 
\quad
\includegraphics[width=4.2cm]{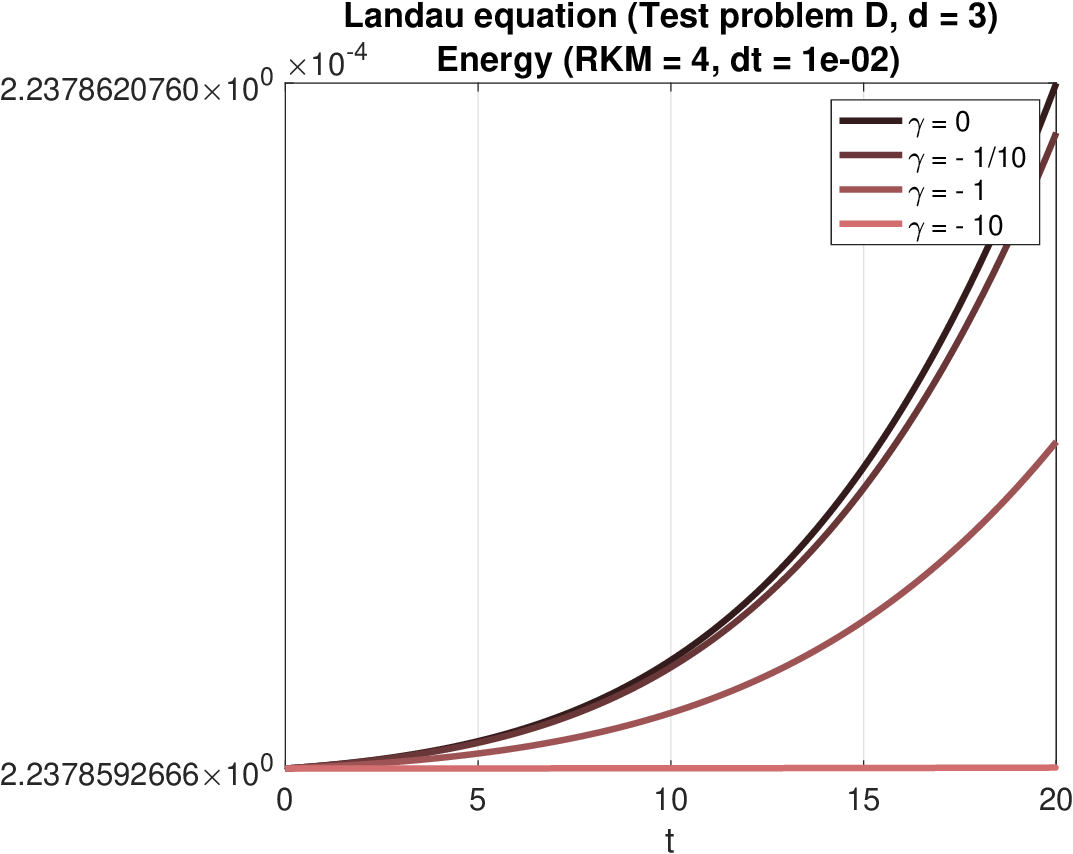} 
\quad
\includegraphics[width=4.2cm]{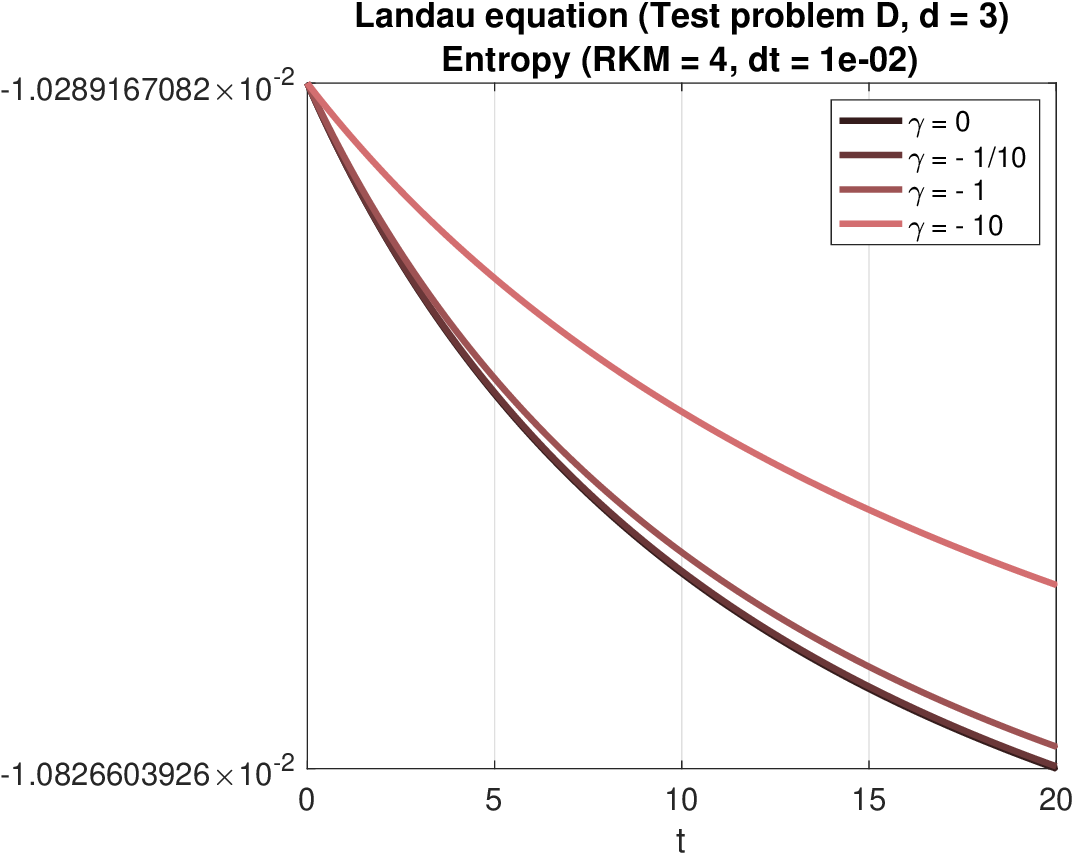} 
\caption{Corresponding results for test problem~D in three dimensions with a significantly smaller truncated domain and a reduced number of time steps.
Different choices of the exponent $\gamma \in \{0, - \, \tfrac{1}{10}, - \, 1, - \, 10\}$.}
\label{fig:MassEtcD3dReduced}
\end{center}
\end{figure}

\MyParagraph{Time integration}
For the time integration of the Landau equation, as indicated in Section~\ref{sec:Introduction}, we combine our approach CST2 with explicit Runge--Kutta methods of non-stiff orders $p \in \{1, 2, 3, 4\}$. 
We once again set $\Omb = [-10, 10]^d$ and choose $100 \times 100$ Fourier functions in two dimensions or $64 \times 64 \times 64$ Fourier functions in three dimensions, respectively.
The results displayed in Figures~\ref{fig:ErrorA2d} and~\ref{fig:ErrorA3d} confirm highly accurate numerical outcomes for the Maxwellian molecules case. Accordingly, we observe conservation of mass as well as nearby conservation of momentum and energy over time, see Figure~\ref{fig:MassEtcA}.
Regarding the decay of entropy, as negative contributions in the range of machine precision may appear, we employ a projection on the strictly positive density function values.
The corresponding results for the significantly more demanding cases of singular integral kernels (Test problem~D) are shown in Figure~\ref{fig:MassEtcD}.
As a final test, we prescribe a significantly smaller truncated domain and larger equidistant time steps.
We observe expected quantitative effects, that is, a certain loss of accuracy, but still retain qualitatively reliable results, see Figure~\ref{fig:MassEtcD3dReduced}. 

\MyParagraph{Solution profiles}
For the relevant cases of Maxwellian molecules (Test problem A) and singular integral kernels (Test problem D), the solution profiles at the initial and final times are shown in Figures~\ref{fig:SolutionAd3} and~\ref{fig:SolutionDd3}.
Movies illustrating the time evolution are found through
\cite{CM2024Fig}.
%%%%%%%%%%%%%%%%%%%%%%%%%%%%%%%%%%%%%%%%%%%%%%%%%%%%%%%%%%%%%%%%%%%%%%%%%%%%%%%%%%%%%%%%%%%%%%%%%%%%%%%%%%%%%%%%%%%%%%%%%%%%%%%%%%%%%%%%%%%%%%%%%% 
%%%%%%%%%%%%%%%%%%%%%%%%%%%%%%%%%%%%%%%%%%%%%%%%%%%%%%%%%%%%%%%%%%%%%%%%%%%%%%%%%%%%%%%%%%%%%%%%%%%%%%%%%%%%%%%%%%%%%%%%%%%%%%%%%%%%%%%%%%%%%%%%%% 
\section*{Acknowledgments}
%%%%%%%%%%%%%%%%%%%%%%%%%%%%%%%%%%%%%%%%%%%%%%%%%%%%%%%%%%%%%%%%%%%%%%%%%%%%%%%%%%%%%%%%%%%%%%%%%%%%%%%%%%%%%%%%%%%%% 
%%%%%%%%%%%%%%%%%%%%%%%%%%%%%%%%%%%%%%%%%%%%%%%%%%%%%%%%%%%%%%%%%%%%%%%%%%%%%%%%%%%%%%%%%%%%%%%%%%%%%%%%%%%%%%%%%%%%%
The research of JAC was supported by the Advanced Grant Nonlocal\--CPD (Non\-local PDEs for Complex Particle Dynamics: Phase Transitions, Patterns and Synchronization) of the European Research Council Executive Agency (ERC) under the European Union’s Horizon 2020 research and innovation programme (grant agreement No. 883363).
JAC was also partially supported by the “Maria de Maeztu” Excellence Unit IMAG, reference CEX2020-001105-M, funded by MCIN/AEI/10.13039/501100011033/.
MT has been supported by the Austrian Science Fund FWF through
a stand-alone project (grant-doi 10.55776/PAT1281625).
JAC and MT acknowledge support from the programme "UIBK Guest Professorship" of the University of Innsbruck. 

\bibliography{LS}
\bibliographystyle{abbrv}

%%%%%%%%%%%%%%%%%%%%%%%%%%%%%%%%%%%%%%%%%%%%%%%%%%%%%%%%%%%%%%%%%%%%%%%%%%%%%%%%%%%%%%%%%%%%%%%%%%%%%%%%%%%%%%%%%%%%%%%%%%%%%%%%%%%%%%%%%%%%%%%%%% 
\appendix
\section{Further approaches}
%%%%%%%%%%%%%%%%%%%%%%%%%%%%%%%%%%%%%%%%%%%%%%%%%%%%%%%%%%%%%%%%%%%%%%%%%%%%%%%%%%%%%%%%%%%%%%%%%%%%%%%%%%%%%%%%%%%%%%%%%%%%%%%%%%%%%%%%%%%%%%%%%% 
\subsection{Approaches~CCT1-2 (Constant kernel) and~CRT1-2 (Regular kernel)}
\label{sec:SectionCRT1-2}
%%%%%%%%%%%%%%%%%%%%%%%%%%%%%%%%%%%%%%%%%%%%%%%%%%%%%%%%%%%%%%%%%%%%%%%%%%%%%%%%%%%%%%%%%%%%%%%%%%%%%%%%%%%%%%%%%%%%%%%%%%%%%%%%%%%%%%%%%%%%%%%%%% 
\MyParagraph{Simplified approaches}
The consideration of the Maxwellian molecules case and test problems involving regular integral kernels and permits significant simplifications concerning the representations of the fundamental integrals.
On account of situations exemplified in~\eqref{eq:TestProblemB}, where the approach~CST2 fails due to the incorrect replacement of the shifted integral domain $v - \Omega$ by the original domain~$\Omega$, we first include the details for approaches~CRT1 and~CCT1 and then, for the sake of completeness, the corresponding relations for approaches~CRT2 and~CCT2.

\MyParagraph{Approach CRT1}
For a regular integral kernel with formal Fourier series expansion
%\begin{subequations}
\begin{equation*}
\varphi(v) = \sum_{\ell \in \ZZ^d} \varphi_{\ell} \, \nFb_{\ell}(v)\,, \quad v \in \Omega\,, 
\end{equation*}
the fundamental integrals are given by 
\begin{equation*}
%\label{eq:ApproachCRT1}
\begin{gathered}
\widetilde{I}(k, m - \ell) = \int_{\Omega} w^k \, \nFb_{m - \ell}(w) \; \dd w\,, \quad k \in \nM_3\,, \\
\widetilde{I}_{ij}^{(m - \ell)}(v) = v_i \, v_j \int_{\Omega} \nFb_{m - \ell}(w) \; \dd w + v_i \int_{\Omega} w_j \, \nFb_{m - \ell}(w) \; \dd w \\
\qquad\qquad\qquad\quad + \; v_j \int_{\Omega} w_i \, \nFb_{m - \ell}(w) \; \dd w + \int_{\Omega} w_i \, w_j \, \nFb_{m - \ell}(w) \; \dd w\,, \\
\nIm_{ij}(v) = \sum_{\ell \in \ZZ^d} \varphi_{\ell} \, \nEb_{- \ell}(b_1) \, \nFb_{\ell}(v) \, \widetilde{I}_{ij}^{(m - \ell)}(v)\,, \\
\nI_{ij}(f)(v) = \sum_{m \in \ZZ^d} f_m \; \nIm_{ij}(v)\,, \quad (i, j) \in \nM_1\,, \\
\nJ_{ijk}(f)(v) = \sum_{m \in \ZZ^d} \mu^{(b_{k1}, \, b_{k2})}_{m_k} \, f_m \; \nIm_{ij}(v)\,, \quad (i, j,k) \in \nM_2\,, \\
v \in \Omega\,.
\end{gathered}
\end{equation*}

\MyParagraph{Approach CCT1}
For a constant integral kernel, that is $\varphi = C$, we arrive at the further simplification
\begin{equation}
\label{eq:ApproachCCT1}
\begin{gathered}
\widetilde{I}(k, m) = \int_{\Omega} w^k \, \nFb_m(w) \; \dd w\,, \quad k \in \nM_3\,, \\
\nIm_{ij}(v) = C \, \bigg(v_i \, v_j \int_{\Omega} \nFb_m(w) \; \dd w + v_i \int_{\Omega} w_j \, \nFb_m(w) \; \dd w \\
\qquad\qquad\qquad\qquad\; + \; v_j \int_{\Omega} w_i \, \nFb_m(w) \; \dd w + \int_{\Omega} w_i \, w_j \, \nFb_m(w) \; \dd w\bigg)\,, \\
\nI_{ij}(f)(v) = \sum_{m \in \ZZ^d} f_m \; \nIm_{ij}(v)\,, \quad (i, j) \in \nM_1\,, \\
\nJ_{ijk}(f)(v) = \sum_{m \in \ZZ^d} \mu^{(b_{k1}, \, b_{k2})}_{m_k} \, f_m \; \nIm_{ij}(v)\,, \quad (i, j,k) \in \nM_2\,, \\
v \in \Omega\,.
\end{gathered}
\end{equation}

\MyParagraph{Approach CRT2}
On the other hand, for a regular integral kernel, we obtain 
\begin{equation}
\begin{gathered}
\nIm_{ij}(v) = \nEb_{-m}(b_1) \, \nFb_{m}(v) \sum_{\ell \in \ZZ^d} \varphi_{\ell} \int_{v - \Omega} u_i \, u_j \; \nFb_{\ell - m}(u) \; \dd u\,, \\
\nI_{ij}(f)(v) = \sum_{m \in \ZZ^d} f_m \; \nIm_{ij}(v)\,, \quad (i, j) \in \nM_1\,, \\
\nJ_{ijk}(f)(v) = \sum_{m \in \ZZ^d} \mu^{(b_{k1}, \, b_{k2})}_{m_k} \, f_m \; \nIm_{ij}(v)\,, \quad (i, j, k) \in \nM_2\,, \\
v \in \Omega\,.
\end{gathered}
\end{equation}

\MyParagraph{Approach CCT2}
In case of a constant integral kernel, the former approach reduces to 
\begin{equation}
\begin{gathered}
\nIm_{ij}(v) = C \, \nFb_m(v) \int_{v - \Omega} u_i \, u_j \, \nEb_{- m}(u) \; \dd u\,, \\
\nI_{ij}(f)(v) = \sum_{m \in \ZZ^d} f_m \; \nIm_{ij}(v)\,, \quad (i, j) \in \nM_1\,, \\
\nJ_{ijk}(f)(v) = \sum_{m \in \ZZ^d} \mu^{(b_{k1}, \, b_{k2})}_{m_k} \, f_m \; \nIm_{ij}(v)\,, \quad (i, j, k) \in \nM_2\,, \\
v \in \Omega\,.
\end{gathered}
\end{equation}
%%%%%%%%%%%%%%%%%%%%%%%%%%%%%%%%%%%%%%%%%%%%%%%%%%%%%%%%%%%%%%%%%%%%%%%%%%%%%%%%%%%%%%%%%%%%%%%%%%%%%%%%%%%%%%%%%%%%%%%%%%%%%%%%%%%%%%%%%%%%%%%%%% 
\subsection{Approaches~NST1 and~NRT1 (Non-Conservative form)}
\label{sec:SectionNST1NRT1}
%%%%%%%%%%%%%%%%%%%%%%%%%%%%%%%%%%%%%%%%%%%%%%%%%%%%%%%%%%%%%%%%%%%%%%%%%%%%%%%%%%%%%%%%%%%%%%%%%%%%%%%%%%%%%%%%%%%%%%%%%%%%%%%%%%%%%%%%%%%%%%%%%% 
\MyParagraph{Non-conservative formulation}
With regard to numerical comparisons for different test problems, it is worth mentioning that alternative approaches in the lines of CST2 and CST1, based on a non-conservative formulation of the Landau operator and hence avoiding numerical differentiation of the integral operator, are advantageous in certain situations, see for instance~\eqref{eq:TestProblemB}.
However, due to the fact that intrinsic properties such as the conservation of mass are lost, additional technicalities are needed, and the overall computation times are higher, we refrain from detailed descriptions and merely include the numerical results obtained for the general approach~NST1 and its simplification~NRT1 in the special case of a regular kernel. 
%%%%%%%%%%%%%%%%%%%%%%%%%%%%%%%%%%%%%%%%%%%%%%%%%%%%%%%%%%%%%%%%%%%%%%%%%%%%%%%%%%%%%%%%%%%%%%%%%%%%%%%%%%%%%%%%%%%%%%%%%%%%%%%%%%%%%%%%%%%%%%%%%% 
\section{Auxiliaries}
\label{sec:Auxiliaries}
%%%%%%%%%%%%%%%%%%%%%%%%%%%%%%%%%%%%%%%%%%%%%%%%%%%%%%%%%%%%%%%%%%%%%%%%%%%%%%%%%%%%%%%%%%%%%%%%%%%%%%%%%%%%%%%%%%%%%%%%%%%%%%%%%%%%%%%%%%%%%%%%%% 
In the following, we collect useful abbreviations and elementary results related to the Fourier spectral method.  
Its efficient implementation relies on fast Fourier techniques.

\MyParagraph{Short notation}
For the purpose of more compact representations, we introduce the common short notation 
\begin{equation}
\label{eq:ShortNotation}  
\begin{gathered}
v^j = v_1^{j_1} \cdots v_d^{j_d}\,, \quad \partial_v^j = \partial_{v_1}^{j_1} \cdots \partial_{v_d}^{j_d}\,, \\
v = (v_1, \dots, v_d) \in \RR^d\,, \quad j = (j_1, \dots, j_d) \in \NN_{\geq 0}^d\,.
\end{gathered}
\end{equation}

\MyParagraph{Truncated velocity domain}
With regard to the numerical approximation of the integral operator in~\eqref{eq:LandauOperator}, it is required to truncate the initially unbounded velocity domain $\Omega = \RR^d$.
Under the reasonable assumption that the density function can be approximated with high accuracy by a localised function, we may replace the original integral domain by a Cartesian product of sufficiently large intervals
\begin{equation}
\label{eq:Omb}
\Omb = [b_{11}, b_{12}] \times \cdots \times [b_{d1}, b_{d2}] \subset \Omega\,. 
\end{equation}

\MyParagraph{Fourier functions}
For Fourier functions depending on a single variable, we in addition specify the dependence on the canonical periodicity interval and the number of oscillations
\begin{subequations}
\label{eq:Fourier}  
\begin{equation}
\begin{gathered}
\nF_{\kappa}^{(\alpha)}(\xi) = \tfrac{1}{\sqrt{\alpha_2 - \alpha_1}} \; \ee^{\, \mu_{\kappa}^{(\alpha)} (\xi - \alpha_1)}\,, \quad
\mu_{\kappa}^{(\alpha)} = \tfrac{2 \, \pi \, \ii \, \kappa}{\alpha_2 - \alpha_1}\,, \\
\nF_{\kappa}^{(\alpha)}(\alpha_1) = \tfrac{1}{\sqrt{\alpha_2 - \alpha_1}} = \nF_{\kappa}^{(\alpha)}(\alpha_2)\,, \quad
\partial_{\xi} \, \nF_{\kappa}^{(\alpha)}(\xi) = \mu_{\kappa}^{(\alpha)} \, \nF_{\kappa}^{(\alpha)}(\xi)\,, \\
\kappa \in \ZZ\,, \quad \alpha = (\alpha_1, \alpha_2) \in \RR^2\,, \quad \alpha_1 < \alpha_2\,, \quad \xi \in \RR\,.
\end{gathered}
\end{equation}
Accordingly, for Fourier functions in several variables, we set 
\begin{equation}
\begin{gathered}
\nFb_m(v) = \nF_{m_1}^{(b_{11}, \, b_{12})}(v_1) \cdots \nF_{m_d}^{(b_{d1}, \, b_{d2})}(v_d)\,, \\
m = (m_1, \dots, m_d) \in \ZZ^d\,, \quad v = (v_1, \dots, v_d) \in \RR^d\,.
\end{gathered}
\end{equation}
Moreover, in view of compact representations, we employ the notation 
\begin{equation}
\begin{gathered}
\mub_m = \big(\mu_{m_1}^{(b_{11}, \, b_{12})}, \dots, \mu_{m_d}^{(b_{d1}, \, b_{d2})}\big) \in \CC^{1 \times d}\,, \\
\nEb_m(v) = \ee^{\, \mu_{m_1}^{(b_{11}, \, b_{12})} \, v_1} \cdots \ee^{\, \mu_{m_d}^{(b_{d1}, \, b_{d2})} \, v_d}\,, \\
\nEb_m(b_1) = \ee^{\, \mu_{m_1}^{(b_{11}, \, b_{12})} \, b_{11}} \cdots \ee^{\, \mu_{m_d}^{(b_{d1}, \, b_{d2})} \, b_{d1}}\,, \\
m \in \ZZ^d\,, \quad v \in \RR^d\,.
\end{gathered}
\end{equation}
\end{subequations}

\MyParagraph{Basic integrals and identities}
Elementary calculations permit to determine the basic integrals 
\begin{equation}
\label{eq:BasicIntegrals}  
\begin{gathered}
I_j(k_j, m_j) = \int_{b_{j1}}^{b_{j2}} w_j^{k_j} \, \nF_{m_j}^{(b_{j1}, \, b_{j2})}(w_j) \; \dd w_j\,, \\
I_j(k_j, m_j) = \begin{cases}
\sqrt{b_{j2} - b_{j1}}\,, & k_j = 0\,, \quad m_j = 0\,, \\
\tfrac{b_{j2}^2 - b_{j1}^2}{2 \, \sqrt{b_{j2} - b_{j1}}}\,, & k_j = 1\,, \quad m_j = 0\,, \\
\tfrac{b_{j2}^3 - b_{j1}^3}{3 \, \sqrt{b_{j2} - b_{j1}}}\,, & k_j = 2\,, \quad m_j = 0\,, \\
0\,, & k_j = 0\,, \quad m_j \in \ZZ \setminus \{0\}\,, \\
\tfrac{\sqrt{b_{j2} - b_{j1}}}{\mu^{(b_{j1}, \, b_{j2})}_{m_j}}\,, & k_j = 1\,, \quad m_j \in \ZZ \setminus \{0\}\,, \\
\tfrac{(b_{j2}^2 - b_{j1}^2) \, \mu^{(b_{j1}, \, b_{j2})}_{m_j} - 2 \, (b_{j2} - b_{j1})}{\sqrt{b_{j2} - b_{j1}} \, \big(\mu^{(b_{j1}, \, b_{j2})}_{m_j}\big)^2}\,, & k_j = 2\,, \quad m_j \in \ZZ \setminus \{0\}\,, 
\end{cases} \\
I(k, m) = \int_{\Omb} w^k \, \nF_m^{(b)}(w) \; \dd w = I_1(k_1, m_1) \cdots I_d(k_d, m_d)\,, \\
k \in \{0, 1, 2\}^d\,, \quad m \in \ZZ^d\,,
\end{gathered}
\end{equation}
see also~\eqref{eq:ShortNotation}.
Due to their close connection to the complex exponential, Fourier functions satisfy identities such as 
\begin{subequations}
\begin{gather}
\label{eq:FourierIdentity1}
\nFb_m(v - u) = \nFb_m(v) \, \nEb_{- m}(u)\,, \\
\label{eq:FourierIdentity2}
\nFb_{\ell}(u) \, \nFb_m(v - u) = \nEb_{-m}(b_1) \, \nFb_m(v) \, \nFb_{\ell - m}(u)\,, \\
\label{eq:FourierIdentity3}
\nFb_{\ell}(v - w) \, \nFb_m(w) = \nEb_{- \ell}(b_1) \, \nFb_{\ell}(v) \, \nFb_{m - \ell}(w)\,, \\
\notag
\ell, m \in \ZZ^d\,, \quad u, v, w \in \RR^d\,.
\end{gather}
\end{subequations}

\MyParagraph{Fourier series expansions}
Accordingly to the regularity of a real-valued function $g: \Omb \to \RR$ and its partial derivatives, Fourier series expansions 
\begin{equation*}
\begin{gathered}
g(v) \approx \sum_{m \in \nM_M} g_m \, \nFb_m(v)\,, \quad 
g_m = \int_{\Omb} g(v) \, \nFb_{-m}(v) \; \dd v\,, \quad m \in \nM_M\,, \\ 
\partial_{v_k} g(v) \approx \sum_{m \in \nM_M} \mu^{(b_{k1}, \, b_{k2})}_{m_k} \, g_m \, \nFb_m(v)\,, \quad k \in \{1, \dots, d\}\,, \\
v \in \Omb\,,
\end{gathered}
\end{equation*}
hold, where the index set
\begin{equation}
\label{eq:nM}
\nM_M = \big\{- \tfrac{M_1}{2}\,, \dots, \tfrac{M_1}{2} - 1\big\} \times \cdots \times \big\{- \tfrac{M_d}{2}\,, \dots, \tfrac{M_d}{2} - 1\big\} \subset \ZZ^d
\end{equation}
is characterised by $M = (M_1, \dots, M_d) \in \NN^d$ comprising even positive integers. 
For the sake of compact representations, we initially identify functions with their formal Fourier series 
\begin{equation*}
g(v) = \sum_{m \in \ZZ^d} g_m \, \nFb_m(v)\,, \quad v \in \Omega\,,
\end{equation*}
and afterwards comment on the practical implementation of the Fourier spectral method.

%%%%%%%%%%%%%%%%%%%%%%%%%%%%%%%%%%%%%%%%%%%%%%%%%%%%%%%%%%%%%%%%%%%%%%%%%%%%%%%%%%%%%%%%%%%%%%%%%%%%%%%%%%%%%%%%%%%%%
%%%%%%%%%%%%%%%%%%%%%%%%%%%%%%%%%%%%%%%%%%%%%%%%%%%%%%%%%%%%%%%%%%%%%%%%%%%%%%%%%%%%%%%%%%%%%%%%%%%%%%%%%%%%%%%%%%%%%

\begin{table}[t!]
\begin{center}
\begin{tabular}{|c|}\hline
\emph{Approach CST2} (Section~\ref{sec:SectionCST2}) \\\hline
Based on the conservative formulation of the Landau operator. \\\hline
Uses numerical differentiation of the integral operator. \\\hline
Adapted to kernels with an isolated singularity at the origin. \\\hline
The integral transform is applied to the singular kernel and its regularisation. \\\hline\hline
\emph{Approach CST1} (Section~\ref{sec:SectionCST1}) \\\hline
Based on the conservative formulation of the Landau operator. \\\hline
Uses numerical differentiation of the integral operator. \\\hline
Adapted to kernels with an isolated singularity at the origin. \\\hline
The integral transform is applied to the singular kernel. \\\hline\hline
\emph{Approaches CRT2 and CRT1} (Section~\ref{sec:SectionCRT1-2}) \\\hline
Simplification of CST2 and CST1 for the case of a regular kernel. \\\hline\hline
\emph{Approach NST1} (Section~\ref{sec:SectionNST1NRT1}) \\\hline
Based on the non-conservative formulation of the Landau operator. \\\hline
Avoids numerical differentiation of the integral operator. \\\hline
Adapted to kernels with an isolated singularity at the origin. \\\hline
The integral transform is applied to the singular kernel and its derivatives. \\\hline\hline
\emph{Approach NRT1} (Section~\ref{sec:SectionNST1NRT1}) \\\hline
Simplification of NST1 for the case of a regular kernel. \\\hline
\end{tabular} \\[2mm] 
\caption{Overview on different approaches for the numerical evaluation of the Landau collision operator~\eqref{eq:LandauOperatorxvt}.}
\label{tab:Table1}
\end{center}
\end{table}

%%%%%%%%%%%%%%%%%%%%%%%%%%%%%%%%%%%%%%%%%%%%%%%%%%%%%%%%%%%%%%%%%%%%%%%%%%%%%%%%%%%%%%%%%%%%%%%%%%%%%%%%%%%%%%%%%%%%% 
% Quadrature
%%%%%%%%%%%%%%%%%%%%%%%%%%%%%%%%%%%%%%%%%%%%%%%%%%%%%%%%%%%%%%%%%%%%%%%%%%%%%%%%%%%%%%%%%%%%%%%%%%%%%%%%%%%%%%%%%%%%% 

\begin{table}[t!]
\begin{center}
\begin{tabular}{|c||c|}\hline
Quadrature on a small neighbourhood & Precomputation time CT \\\hline
Quadrature on the whole domain & $88 \times$CT \\\hline
\end{tabular} \\[2mm] 
\caption{Test problem~$C$ (regular integral kernel, unbounded domain, known solution) in two dimensions.
Numerical evaluation of the Landau operator based on $256 \times 256$ uniform grid points covering the truncated velocity domain $[- \, 10, 10] \times [- \, 10, 10]$.
Precomputation times observed for a quadrature approximation based on $5 \times 5$ grid points versus a quadrature approximation on the whole domain based on $256 \times 256$ grid points.
In both cases, an overall relative accuracy of about $4 \cdot 10^{- \, 11}$ is obtained.}
\label{tab:Table2}
\end{center} 
\end{table}

%%%%%%%%%%%%%%%%%%%%%%%%%%%%%%%%%%%%%%%%%%%%%%%%%%%%%%%%%%%%%%%%%%%%%%%%%%%%%%%%%%%%%%%%%%%%%%%%%%%%%%%%%%%%%%%%%%%%%
% Singular kernel 
%%%%%%%%%%%%%%%%%%%%%%%%%%%%%%%%%%%%%%%%%%%%%%%%%%%%%%%%%%%%%%%%%%%%%%%%%%%%%%%%%%%%%%%%%%%%%%%%%%%%%%%%%%%%%%%%%%%%%

%%%%%%%%%%%%%%%%%%%%%%%%%%%%%%%%%%%%%%%%%%%%%%%%%%%%%%%%%%%%%%%%%%%%%%%%%%%%%%%%%%%%%%%%%%%%%%%%%%%%%%%%%%%%%%%%%%%%% 
% Comparisons (A / B / C / D)
%%%%%%%%%%%%%%%%%%%%%%%%%%%%%%%%%%%%%%%%%%%%%%%%%%%%%%%%%%%%%%%%%%%%%%%%%%%%%%%%%%%%%%%%%%%%%%%%%%%%%%%%%%%%%%%%%%%%% 

%%%%%%%%%%%%%%%%%%%%%%%%%%%%%%%%%%%%%%%%%%%%%%%%%%%%%%%%%%%%%%%%%%%%%%%%%%%%%%%%%%%%%%%%%%%%%%%%%%%%%%%%%%%%%%%%%%%%% 
% Error over time (A)
%%%%%%%%%%%%%%%%%%%%%%%%%%%%%%%%%%%%%%%%%%%%%%%%%%%%%%%%%%%%%%%%%%%%%%%%%%%%%%%%%%%%%%%%%%%%%%%%%%%%%%%%%%%%%%%%%%%%% 

%%%%%%%%%%%%%%%%%%%%%%%%%%%%%%%%%%%%%%%%%%%%%%%%%%%%%%%%%%%%%%%%%%%%%%%%%%%%%%%%%%%%%%%%%%%%%%%%%%%%%%%%%%%%%%%%%%%%% 
% Mass over time etc. (A / D)
%%%%%%%%%%%%%%%%%%%%%%%%%%%%%%%%%%%%%%%%%%%%%%%%%%%%%%%%%%%%%%%%%%%%%%%%%%%%%%%%%%%%%%%%%%%%%%%%%%%%%%%%%%%%%%%%%%%%% 

%%%%%%%%%%%%%%%%%%%%%%%%%%%%%%%%%%%%%%%%%%%%%%%%%%%%%%%%%%%%%%%%%%%%%%%%%%%%%%%%%%%%%%%%%%%%%%%%%%%%%%%%%%%%%%%%%%%%% 
% Solution profiles (A / D)
%%%%%%%%%%%%%%%%%%%%%%%%%%%%%%%%%%%%%%%%%%%%%%%%%%%%%%%%%%%%%%%%%%%%%%%%%%%%%%%%%%%%%%%%%%%%%%%%%%%%%%%%%%%%%%%%%%%%% 

%%%%%%%%%%%%%%%%%%%%%%%%%%%%%%%%%%%%%%%%%%%%%%%%%%%%%%%%%%%%%%%%%%%%%%%%%%%%%%%%%%%%%%%%%%%%%%%%%%%%%%%%%%%%%%%%%%%%% 
%%%%%%%%%%%%%%%%%%%%%%%%%%%%%%%%%%%%%%%%%%%%%%%%%%%%%%%%%%%%%%%%%%%%%%%%%%%%%%%%%%%%%%%%%%%%%%%%%%%%%%%%%%%%%%%%%%%%% 
\end{document}